\documentclass[A4paper,11pt,twoside,openright]{report}
\usepackage{amsfonts}
\usepackage{amsmath}
\usepackage{amssymb}
\usepackage{graphicx}
\usepackage{rlepsf}
\usepackage{amscd}
\begin{document}
\date{}

\newtheorem{theorem}{Theorem}
\newtheorem{acknowledgment}[theorem]{Acknowledgment}
\newtheorem{algorithm}[theorem]{Algorithm}
\newtheorem{axiom}[theorem]{axiom}
\newtheorem{case}[theorem]{Case}
\newtheorem{claim}[theorem]{Claim}
\newtheorem{conclusion}[theorem]{Conclusion}
\newtheorem{condition}[theorem]{Condition}
\newtheorem{conjecture}[theorem]{Conjecture}
\newtheorem{corollary}[theorem]{Corollary}
\newtheorem{criterion}[theorem]{Criterion}
\newtheorem{definition}{Def\mbox{}inition}
\newtheorem{example}{Example}
\newtheorem{exercise}[theorem]{Exercise}
\newtheorem{lemma}{\indent Lemma}
\newtheorem{notation}[theorem]{Notation}
\newtheorem{problem}[theorem]{Problem}
\newtheorem{proposition}{Proposition}
\newtheorem{remark}[theorem]{Remark}
\newtheorem{solution}[theorem]{Solution}
\newtheorem{summary}[theorem]{Summary}
\newcommand{\ud}{\mathrm{d}}

\def\gcd{\mathop{\rm gcd}}
\def\Ker{\mathop{\rm Ker}}
\def\max{\mathop{\rm max}}
\def\map{\mathop{\rm map}}
\def\lcm{\mathop{\rm lcm}}
\def\kraj{\hfill\rule{6pt}{6pt}}
\def\diag{\mathop{\rm diag}}
\def\span{\mathop{\rm span}}
\def\deg{\mathop{\rm deg}}
\def\rank{\mathop{\rm rank}}
\def\sgn{\mathop{\rm sgn}}
\def\kvn{\{n\}_q}
\def\F{\mathbb{F}}
\def\R{\mathbb{R}}
\def\C{\mathcal{C}}
\def\P{\mathcal{P}}
\def\A{\mathcal{A}}
\def\D{\mathcal{D}}
\def\H{\mathcal{H}}
\def\N{\mathbb{N}}
\def\K{\mathbb{K}}
\def\Z{\mathbb{Z}}
\def\Q{\mathbb{Q}}
\def\X{\qbezier(0.00,0.00)(0.50,1.00)(1.00,2.00)
\qbezier(1.00,0.00)(0.80,0.40)(0.60,0.80)
\qbezier(0.00,2.00)(0.20,1.60)(0.40,1.20)
}

\def\Y{\qbezier(1.00,0.00)(0.50,1.00)(0.00,2.00)
\qbezier(0.00,0.00)(0.20,0.40)(0.40,0.80)
\qbezier(1.00,2.00)(0.80,1.60)(0.60,1.20)
}

\def\O{\qbezier(0.00,0.00)(0.20,1.00)(0.00,2.00)
\qbezier(1.00,0.00)(0.80,1.00)(1.00,2.00)
}
\def\OF{\qbezier(0.00,0.00)(0.00,1.00)(0.00,2.00)
\qbezier(1.00,0.00)(1.00,1.00)(1.00,2.00)
}

\def\l{
\qbezier(0.00,0.00)(0.50,1.20)(1.00,0.00)
\qbezier(0.00,2.00)(0.50,0.80)(1.00,2.00)
}

\def\LPP{
\qbezier(0.00,0.00)(0.35,0.70)(0.70,1.40)
\qbezier(0.70,0.60)(0.65,0.70)(0.60,0.80)
\qbezier(0.00,2.00)(0.20,1.60)(0.40,1.20)
\qbezier(0.7,1.4)(1.3,2.6)(1.3,1)
\qbezier(0.7,0.6)(1.3,-0.6)(1.3,1)}

\def\LPM{\qbezier(0.00,0.00)(0.2,0.40)(0.40,0.80)
\qbezier(0.60,1.20)(0.65,1.30)(0.70,1.40)
\qbezier(0.00,2.00)(0.35,1.30)(0.70,0.60)
\qbezier(0.7,1.4)(1.3,2.6)(1.3,1)
\qbezier(0.7,0.6)(1.3,-0.6)(1.3,1)}

\def\XS{\qbezier(0.00,0.00)(0.50,1.00)(1.00,2.00)
\qbezier(1.00,0.00)(0.80,0.40)(0.60,0.80)
\qbezier(0.00,2.00)(0.20,1.60)(0.40,1.20)

\qbezier(0.0,2.00)(0.0,1.85)(0,1.70)
\qbezier(1.00,2.00)(1,1.85)(1.0,1.70)
\qbezier(0.00,2.00)(0.10,1.90)(0.20,1.80)
\qbezier(1.00,2.00)(0.9,1.9)(0.8,1.80)}

\def\YS{\qbezier(1.00,0.00)(0.50,1.00)(0.00,2.00)
\qbezier(0.00,0.00)(0.20,0.40)(0.40,0.80)
\qbezier(1.00,2.00)(0.80,1.60)(0.60,1.2)
\qbezier(0.0,2.00)(0.0,1.85)(0,1.70)
\qbezier(1.00,2.00)(1,1.85)(1.0,1.70)
\qbezier(0.00,2.00)(0.10,1.90)(0.20,1.80)
\qbezier(1.00,2.00)(0.9,1.9)(0.8,1.80)}

\def\GR{\qbezier(0.00,0.00)(0.25,0.25)(0.50,0.50)
\qbezier(1.00,0.00)(0.75,0.25)(0.50,0.50)
\qbezier(0.00,2.00)(0.25,1.75)(0.50,1.5)
\qbezier(1.0,2.00)(0.75,1.75)(0.5,1.50)

\qbezier(0.00,2.00)(0,1.9)(0,1.80)
\qbezier(0.00,2.00)(0.1,2)(0.20,2)
\qbezier(1.00,2.00)(1,1.9)(1,1.80)
\qbezier(1.00,2.00)(0.9,2)(0.8,2)
\linethickness{2.5pt}
\qbezier(0.5,0.5)(0.5,1)(0.50,1.5)}

\def\GRF{\qbezier(0.00,0.00)(0.25,0.25)(0.50,0.50)
\qbezier(1.00,0.00)(0.75,0.25)(0.50,0.50)
\qbezier(0.00,2.00)(0.25,1.75)(0.50,1.5)
\qbezier(1.0,2.00)(0.75,1.75)(0.5,1.50)

%\qbezier(0.00,2.00)(0,1.9)(0,1.80)
%\qbezier(0.00,2.00)(0.1,2)(0.20,2)
%\qbezier(1.00,2.00)(1,1.9)(1,1.80)
%\qbezier(1.00,2.00)(0.9,2)(0.8,2)
\linethickness{2.5pt}
\qbezier(0.5,0.5)(0.5,1)(0.50,1.5)}

\def\OS{\qbezier(0.00,0.00)(0.60,1.00)(0.00,2.00)
\qbezier(1.00,0.00)(0.40,1.00)(1.00,2.00)
\qbezier(0.0,2.00)(0.0,1.85)(0,1.70)
\qbezier(1.00,2.00)(1,1.85)(1.0,1.70)
\qbezier(0.00,2.00)(0.10,1.90)(0.20,1.80)
\qbezier(1.00,2.00)(0.9,1.9)(0.8,1.80)}

\pagenumbering{roman}

\thispagestyle{empty}
{\Large{
\begin{center}{\textbf{ UNIVERSIDADE T\'ECNICA DE LISBOA\\
 INSTITUTO SUPERIOR T\'ECNICO}}\end{center}
}}
\vspace{2.2cm}
\begin{center}{\LARGE{\textbf{Khovanov homology of links and graphs}}}\end{center}

\vspace{1.2cm}

{\Large{\begin{center} \textbf{Marko Sto\v si\'c}\end{center}
}}
\vspace{1cm}

{\large{\begin{center}Disserta\c c\~ao para obten\c c\~ao do Grau de\\ Doutor em Matem\'atica\end{center}}}

\vspace{1cm}
{\large {\noindent Orientador: Doutor Roger Francis Picken}}

\vspace{0.7cm}
{\large{
\begin{center}
\textbf{J\'uri\\}\end{center}}}

{\large{\noindent Presidente: Reitor da Universidade T\'ecnica de Lisboa}}\\

{\large{\noindent Vogais: \hspace{4pt} Doutor Mikhail Khovanov\\
\indent \indent \indent Doutor Rui Ant\'onio Loja Fernandes\\
\indent \indent \indent Doutor Roger Francis Picken\\
\indent \indent \indent Doutor Paul Turner\\
\indent \indent \indent Doutor Marco Arien Mackaay\\
\indent \indent \indent Doutor Gustavo Rui Gon\c calves Fernandes de Oliveira Granja}}

{\large
\vspace{1cm}
\begin{center}\textbf{Mar\c co de 2006}\end{center} 
}

%{\Large{
%\begin{center}{\textbf{ UNIVERSIDADE T\'ECNICA DE LISBOA\\
% INSTITUTO SUPERIOR T\'ECNICO}}\end{center}
%}}
%\vspace{4cm}
%\begin{center}{\LARGE{\textbf{Khovanov homology of links and %graphs}}}\end{center}
%\vspace{2cm}
%{\Large{\begin{center} \textbf{Marko Sto\v si\'c}\end{center}
%}}
%\vspace{1.5cm}
%{\Large{\begin{center}Disserta\c c\~ao para obten\c c\~ao do Grau de\\ %Doutor em Matem\'atica\end{center}}}
%\vspace{1cm}
%{\Large{
%\begin{center}
%DOCUMENTO PROVIS\'ORIO\end{center}
%{\large
%\vspace{1cm}
%\begin{center}Janeiro de 2006\end{center} 
%}}}

%{\Large{
%\begin{center} UNIVERSIDADE T\'ECNICA DE LISBOA\\
% INSTITUTO SUPERIOR T\'ECNICO\end{center}
%}}
%\vspace{4cm}
%\begin{center}{\Large{\textbf{Khovanov homology of links and %graphs}}}\end{center}
%
%\vspace{2cm}
%
%{\large{\begin{center} \textbf{Marko Sto\v si\'c}\\
%(Licenciado)\end{center}
%}}
%\vspace{0.5cm}
%
%{\large{\begin{center}Disserta\c c\~ao para obten\c c\~ao do Grau de\\ %Doutor em Matem\'atica\end{center}}}
%
%\vspace{2cm}
%{\large{
%\begin{center}
%DOCUMENTO PROVIS\'ORIO\end{center}
%
%\vspace{1cm}
%\begin{center}Janeiro 2006\end{center} 
%}}
%\newpage
%
%

\arraycolsep 6pt

\chapter*{}
\bigskip
\begin{flushright}
{\Large{
Mojoj \mbox{}\mbox{} Maj\mbox{}i
}}
\end{flushright}

\chapter*{Resumo}
\addcontentsline{toc}{chapter}{Resumo}

Nesta tese, trabalhamos com a homologia de Khovanov de enlaces e com a homologia de grafos. A homologia de Khovanov de enlaces consta dum complexo graduado de cadeias que \'e um invariante de enlace, cuja caracter\'istica de Euler \'e igual ao polin\'omio de Jones do enlace, e portanto pode ser vista como a ``categorif\mbox{}ica\c c\~ao" do polin\'omio de Jones. \\
\indent Provamos que o primeiro grupo de homologia dum n\'o de tran\c ca positiva ({\em positive braid knot}) \'e trivial. Depois, provamos que os n\'os de toro, n\~ao alternados s\~ao homologicamente espessos ({\em thick}). Tamb\'em, demonstramos que podemos diminuir o n\'umero de tor\c c\~oes completas de n\'os de toro sem mudar a homologia de baixo grau, implicando a exist\^encia de homologia est\'avel de n\'os de toro. Mais, provamos a maior parte das propriedades anteriores tamb\'em para a homologia de Khovanov-Rozansky. \\
\indent Na homologia de grafos, categorif\mbox{}icamos o polin\'omio dicrom\'atico (e em consequ\^encia, o polin\'omio de Tutte) para  grafos, atrav\'es da categorif\mbox{}ica\c c\~ao dum conjunto inf\mbox{}inito de especializa\c c\~oes duma vari\'avel. Tamb\'em, categori\-f\mbox{}icamos explicitamente a especializa\c c\~ao que \'e analoga ao polin\'omio de Jones dum enlace alternado correspondente ao grafo inicial. Mais, categorif\mbox{}icamos explicitamente o polin\'omio dicrom\'atico completo de duas vari\'aveis de grafos, usando os complexos de Koszul.\\
\\
\\
\textbf{Palavras-chave: } homologia de Khovanov, polin\'omio de Jones, enlace, n\'os de toro, grafo, polin\'omio dicrom\'atico

\chapter*{Abstract}
\addcontentsline{toc}{chapter}{Abstract}

In this thesis we work with Khovanov homology of links and its generalizations, as well as with the homology of graphs. Khovanov homology of links consists of graded chain complexes which are link invariants, up to chain homotopy, with graded Euler characteristic equal to the Jones polynomial of the link. Hence, it can be regarded as the ``categorif\mbox{}ication" of the Jones polynomial.\\
\indent We prove that the f\mbox{}irst homology group of positive braid knots is trivial. Futhermore, we prove that non-alternating torus knots are homologically thick. In addition, we show that we can decrease the number of full twists of torus knots without changing low-degree homology and consequently that there exists stable homology for torus knots. We also prove most of the above properties for Khovanov-Rozansky homology.\\
\indent Concerning graph homology, we categorify the dichromatic (and consequently Tutte) polynomial for graphs, by categorifying an inf\mbox{}inite set of its one-variable specializations. We categorify explicitly the one-variable specialization that is an analog of the Jones polynomial of an alternating link corresponding to the initial graph. Also, we categorify explicitly the whole two-variable dichromatic polynomial of graphs by using Koszul complexes.\\
\\
\\
\textbf{Key-words: }  Khovanov homology, Jones polynomial, link, torus knot, graph, dichromatic polynomial

\chapter*{Acknowledgements}
\addcontentsline{toc}{chapter}{Acknowledgements}

First of all, I want to thank my supervisor Roger Picken 
for allowing me to explore Mathematics during my PhD studies and to f\mbox{}ind the most appropriate topic to work on, and for all the help during the last four years. Also, I am thankful to him for giving me the opportunity to get to know Portugal and to live in this beautiful country. \\
\indent I specially thank Jacob Rasmussen, Jozef Przytycki, Yongwu Rong and the whole group from George Washington University and, most of all, Mikhail Khovanov for all the hospitality, the discussions, both of a mathematical and nonmathematical nature, and for the great time I spent during my visit to the United States in the Autumn of 2005.\\
\indent I am extremely grateful to the Funda\c c\~ao para a Ci\^encia e a Tecnologia (FCT) for the material backing during my PhD studies, as well as for the f\mbox{}inancial support for the trip to the USA.\\

\indent Finally, most of all, I thank my wife Marija Dodig for all the love and support that she is giving me. 

\tableofcontents

\chapter{Introduction}
\pagenumbering{arabic}

One of the most interesting and promising recent developments in mathematics in general is the latest construction in knot theory called the ``cate\-gorif\mbox{}ication" of link invariants, initiated by M. Khovanov in  \cite{kov}. This construction opened new horizons by providing tools for approaching and solving various problems in knot theory and low-dimensional topology. Since the topic is very new, the frontiers of its applications are still far away, and  its potential is  
enormous.

\section{Brief historical review of polynomial link invariants}

Mathematical properties of knots -- circles smoothly embedded in 3-space -- have been studied for a long time. Although the knots are 3-dimensional objects, there is a nice way to present them in two dimensions by a planar projection (also called a regular diagram). The f\mbox{}irst attempts to tabulate  knots date from the late XIX century. The knots and links (f\mbox{}inite disjoint collection of knots) are naturally identif\mbox{}ied by isotopy. Hence in order to distinguish two knots (or links) one needs an invariant, i.e. a function on the set of knots which has the same value on isotopic knots.\\
\indent A f\mbox{}irst step toward f\mbox{}inding  knot invariants was made by K. Reidemeister in \cite{rei}, who proved that two planar projections represent isotopic knots (or links) if and only if one can be obtained from the another by the f\mbox{}inite application of the famous Reidemeister moves. Hence, a function on the set of planar diagrams will be a link invariant if and only if it is invariant under Reidemeister moves. \\
\indent Even before the discovery of the Reidemeister moves, in 1923 in \cite{ax1}, Alexander def\mbox{}ined  geometrically (by working in 3-space) a polynomial link  invariant which was very good at distinguishing knots. Later, in 1970 in \cite{con}, Conway showed that the Alexander polynomial can also be def\mbox{}ined combinatorially by working with planar projections and that it satisf\mbox{}ies a certain skein relation. \\
\indent In 1984, V.F.R. Jones revolutionarized  knot theory by def\mbox{}ining  a  polynomial invariant which satisf\mbox{}ies a dif\mbox{}ferent skein relation to the Alexander-Conway one. The original Jones construction  uses von Neumann algebras and is rather complicated (\cite{jones}). However, in 1987, L. Kauf\mbox{}fman in \cite{kaufss} introduced a state-sum model  construction of the Jones polynomial which is purely combinatorial and remarkably simple. Although the def\mbox{}ining relations of the Alexander-Conway and Jones polynomials are very similar, their properties are quite dif\mbox{}ferent. For instance, the Jones polynomial can distinguish knots from their mirror images, and it was used to prove some old Tait conjectures on regular diagrams of  alternating links (\cite{mus}, \cite{kauf}).  \\
\indent  Soon after the discovery of the Jones polynomial there was an explosion of  generalizations, which culminated in the two-variable HOMFLYPT polynomial named after its many dicoverers (\cite{homfly}, \cite{pt}). It contains the Alexander-Conway and Jones polynomial as  specializations and it is the most general polynomial invariant that possesses a  skein relation. Alternatively, the two-variable HOMFLYPT polynomial can be given by an inf\mbox{}inite sequence of one-variable specializations, one for each positive integer $n$, which is enough to recover the whole HOMFLYPT polynomial. Finally, in \cite{moy}  a state-sum model was obtained for the $n$-specializations of the HOMFLYPT polynomial.\\  
\indent As was shown in the late 80's, \cite{restur}, \cite{tur}, the $n$-specializations of the HOMFLYPT polynomial can be obtained by using the representation theory of the quantum $sl(n)$ group. Hence, the $n$-specializations of the HOMFLYPT polynomial are also called the $sl(n)$ link polynomials. Also, the homology theories that categorify them (see the following section), i.e. Khovanov homology for the $sl(2)$ case (see Section \ref{kovnot}), and Khovanov-Rozansky homology for the general $sl(n)$ case (see Section \ref{kovroznot}), are also called $sl(2)$-link homology and $sl(n)$-link homology, respectively.

\section{Link homology}

\indent In 1999, M. Khovanov in \cite{kov} gave a completely new perspective on link polynomials. For each link $L$ in $S^3$ he def\mbox{}ined a graded chain complex, with grading preserving dif\mbox{}ferentials, whose graded Euler characteristic is equal to the Jones polynomial of the link $L$. This is done by starting from the state-sum expression for the Jones polynomial (which is written as an alternating sum), then constructing for each term a graded module whose graded dimension is equal to the value of that term, and f\mbox{}inally, def\mbox{}ining the dif\mbox{}ferentials as appropriate grading preserving maps, so that the complex obtained is a link invariant (up to chain homotopy). It has also been shown by M. Khovanov in \cite{tan} and in a dif\mbox{}ferent way by M. Jacobsson in \cite{jac}, that the homology groups are functorial, i.e. that there exists a functor from the category of links as objects and cobordisms between them (surfaces in 4-space), to the category of abelian groups, such that on links it is given by the homology of the Khovanov complex. Also, Khovanov homology is strictly stronger than the Jones polynomial, i.e. there exist knots with the same Jones polynomial and dif\mbox{}ferent Khovanov homology (\cite{bn}). The situation is similar to the Euler characteristic of a manifold: it can be seen as an alternating sum of the dimensions of the homology groups of a manifold, which are functorial, and which are stronger invariants than the Euler characteristic itself.\\
\indent Although the theory is rather new, it already has strong applications in low-dimensional topology. For instance, the short proof of the Milnor conjecture by Rasmussen in \cite{ras}, as well as the proof of the existence of exotic dif\mbox{}ferential structures on $\R^4$ (\cite{ras2}, \cite{xy}), which were previously accessible only by 
%the huge machinery of 
gauge theory (\cite{km}, \cite{don}).\\
\indent The advantage of Khovanov homology theory is that its def\mbox{}inition is combinatorial and since there is a straightforward algorithm for computing it, it is (theoretically) highly calculable. Nowadays there are several computer programs \cite{bnprog}, \cite{shum} that  can calculate ef\mbox{}fectively Khovanov homology of links with up to 50 crossings.\\
\indent Based on the calculations there are many conjectures about the properties of link homology, see e.g. \cite{bn}, \cite{pat}, \cite{dgr}. Some of the properties have been proved by now (see \cite{lee}, \cite{lee2}), but many of them are still open. \\

\indent In this thesis we solve some  conjectures concerning  Khovanov homology. First, in Chapter \ref{braid}, we prove the following conjecture from \cite{pat}:\\

\textbf{Theorem \ref{braidgl}:} \textit{The f\mbox{}irst homology group of a  positive braid knot is trivial.} \\

\indent The Khovanov homology of a link is naturally bigraded, and it is common practice to represent it in a planar array, by writing at the position $(i,j)$ the rank of the homology group $\H^{i,j}(L)$. As shown by calculations, the nonzero entries are grouped along the diagonals. The homology of every knot occupies at least two diagonals, and the knots whose homology occupy exactly two diagonals are called homologically thin. All other knots are called homologically thick.\\
\indent The torus knots -- the knots that can be drawn on a surface of a  solid torus without intersection -- constitute an important class of links. They are classif\mbox{}ied up to isotopy by two integers, $p$ and $q$. Every torus knot is a positive braid knot.  Hence, their f\mbox{}irst homology group is trivial by  Theorem \ref{braidgl}. In Chapter \ref{torus} we also obtain the following:\\

\textbf{Corollary \ref{C1}:} \textit{The torus knots $T_{p,q}$, for $3\le p \le q$ are homologically thick.}  \\

Moreover, in the course of the proof, we obtain  relations between the homology groups of the torus knot $T_{p,q}$ and the torus knot $T_{p,q-1}$ (see Theorem \ref{T3} in Chapter \ref{torus}). As a consequence, this gives low-degree homology groups of torus knots and also a stability property for the Khovanov homology of torus knots, conjectured by Dunf\mbox{}ield, Gukov and Rasmussen in \cite{dgr}.\\

\indent \textbf{Theorem \ref{stablet}:} \textit{There exists stable Khovanov homology for torus knots.}\\

\indent An analogous categorif\mbox{}ication of the $n$-specializations of the HOMFLYPT polynomial was carried out by M. Khovanov and L. Rozansky in 2004 (\cite{kovroz}). For each positive integer $n$ and link $L$ in $S^3$ they def\mbox{}ined a graded chain complex, with grading preserving dif\mbox{}ferentials, whose graded Euler characteristic is equal to the $n$-specialization of the HOMFLYPT polynomial of the link $L$. The construction uses the state-sum model for the HOMFLYPT polynomial (\cite{moy}) and is analogous to the categorif\mbox{}ication of the Jones polynomial: it uses the same cubic complex construction, and there exists a  similar long exact sequence in homology. However, since the state-sum model for the HOMFLYPT polynomial is much more complicated than  Kauf\mbox{}fman's state-sum model, the explicit calculation of the homology groups is practically impossible. Consequently, the values of the $sl(n)$-link homology are known only for a very small class of knots -- the two-bridge knots (see \cite{ras2b}).\\
\indent On the other hand, our approach to the properties of Khovanov homology for positive braid knots and torus knots, essentially relies on the form of the construction (cubic complex) and on the existence of a  long exact sequence in homology. Hence, with some minor modif\mbox{}ications of the proofs in the $sl(2)$ case we obtain.\\

\textbf{Theorem \ref{slngl} (Chapter \ref{braid}):} \textit{The f\mbox{}irst $sl(n)$-homology group of a positive braid knot is trivial for every positive integer $n$}.\\

%Also, from the formula (\ref{spustsln}) in Chapter \ref{torus} we obtain

\textbf{Theorem \ref{torussln} (Chapter \ref{torus}):} \textit{There exists stable $sl(n)$-homology for torus knots.}  

\section{Graph homology}

Apart from the categorif\mbox{}ications of the Jones and of the HOMFLYPT polynomial mentioned in the previous section, there also exist categorif\mbox{}ications of various other polynomial link invariants, see \cite{kov3}, \cite{kol}, \cite{bn1}. All these constructions start from a state-sum expression of the polynomial invariant, and then ``lift" it to a chain complex.\\
\indent Since the chromatic polynomial for graphs, \cite{bol}, also possess a state-sum expression, it means that one can try to categorify it as well. This was done  by L. Helme-Guizon and Y. Rong in 2004, \cite{gr}.   For each graph $G$, they def\mbox{}ined a graded chain complex whose graded Euler characteristic is equal to the chromatic polynomial of $G$. Since for graphs there is no isotopy invariance, and therefore no analogs of Reidemeister moves, the homology groups are automatically invariants. Hence, there are no restrictions on the algebraic structure in the categorif\mbox{}ication of graph invariants and consequently L. Helme-Guizon and Y. Rong, together with J. Przytycki have generalized their construction for an arbitrary algebra, \cite{gr2}, \cite{grp}. On the other hand, since the structure of the homology for graphs is simpler than that for links, it can be seen as a step towards a  better understanding of  link homology (\cite{grp}, \cite{pr}).\\
\indent A natural generalization of the chromatic polynomial, is the two-variable dichromatic polynomial (see \cite{kauf}), which is up to a change of variables equal to the important Tutte polynomial of graphs (\cite{bol}, \cite{jaeger}).  Since the dichromatic polynomial is a two-variable polynomial which contains the chromatic polynomial as an one-variable specialization, we have a similar situation to that for the  HOMFLYPT polynomial and the Jones polynomial.\\
\indent In this thesis we categorify the dichromatic polynomial (and consequently the Tutte polynomial). In Chapter \ref{graf},  we def\mbox{}ine an inf\mbox{}inite set of one-variable specializations of the dichromatic polynomial and then we categorify each of these polynomials. We also give an explicit computation of the homology groups of polygon graphs, where we use the known relation between  graph homology and  Hochschild homology, see \cite{pr}.   \\
\indent To each planar graph $G$ there corresponds bijectively, up to  mirror image, an alternating link $L^G$. Then the Jones polynomial of an alternating link $L^G$ corresponds to a certain one-variable specialization of the dichromatic polynomial, $J_G(q)$, of an appropriate planar graph  $G$ (see Appendix \ref{a1}). However, the set of one-variable specializations that we def\mbox{}ined previously ``misses"   the polynomial $J_G(q)$, so we also def\mbox{}ine a second set of one-variable specializations of the dichromatic polynomial that contains $J_G(q)$, and then categorify each of these polynomials.\\

\indent In 2005 there was a new development, when M. Khovanov and L. Rozansky in \cite{KR2} def\mbox{}ined triply-graded link homology, i.e. to each link $L$ they associated a doubly-graded chain complex, whose doubly-graded Euler characteristic is equal to the full HOMFLYPT polynomial of $L$. The method used in the construction, i.e. a cubic complex, is the same as in all previous categorif\mbox{}ications -- the only dif\mbox{}ference is that in the explicit def\mbox{}inition of chain groups they use Koszul complexes (Koszul complexes are also used in \cite{kovroz} in the categorif\mbox{}ication of the $n$-specializations of the HOMFLYPT polynomial). In Section \ref{coment} of Chapter \ref{Koszul} we give a slightly simplif\mbox{}ied version of this construction. Then, in Section \ref{newgraf} of Chapter \ref{Koszul} we categorify the whole two-variable dichromatic polynomial for graphs. Specif\mbox{}ically, for each graph $G$ we construct  a doubly-graded chain complex whose doubly-graded Euler characteristic is equal to the dichromatic polynomial. The construction is in the spirit of the construction of the triply-graded theory for links, described in Section \ref{coment}. \\ 

\section{Organization of the thesis}

\indent The organization of the thesis is as follows: in Chapter \ref{not} we establish the notation that we will use throughout the thesis and  recall brief\mbox{}ly the basic def\mbox{}initions (Jones and HOMFLYPT polynomial and $sl(2)$ and $sl(n)$ link homology).\\
\indent  In Chapter \ref{graf} we  def\mbox{}ine two dif\mbox{}ferent sets of one-variable specializations 
of the dichromatic polynomial,  and then categorify each of the  specializations. Also we present some calculations for polygon graphs and give an explicit categorif\mbox{}ication of the analog of the Jones polynomial of the corresponding alternating link.\\
\indent  In the f\mbox{}irst section of  Chapter \ref{Koszul} we give a  simplif\mbox{}ied def\mbox{}inition of the triply-graded homology theory for links. Then, in  the second section of the same chapter we def\mbox{}ine an analogous triply-graded theory for graphs, i.e. we def\mbox{}ine a doubly-graded chain complex whose doubly-graded Euler characteristic is equal to the two-variable dichromatic polynomial.\\
\indent In the next two chapters we deal with the properties of link homology. In  Chapter \ref{braid} we prove that the f\mbox{}irst (Khovanov) homology of an arbitrary positive braid knot is trivial. We also extend the proof to show that the f\mbox{}irst $sl(n)$ homology group of a positive braid knot is trivial for every positive integer $n$.\\
\indent Finally, in Chapter \ref{torus} we prove that all non-alternating torus knots are homologically thick. In the course of the proof we also obtain a stability property of  Khovanov homology for torus knots and  explicit values for the homology of torus knots for low homological degrees. Also, we  show that an analogous stability property holds for the $sl(n)$ homology of torus knots.\\
\indent At the end of the thesis we have added two Appendices, since we wanted the thesis to be as self-contained as possible. In the f\mbox{}irst one we explain the relation between planar graphs and alternating links, as well as between the Jones polynomial of an alternating link and a certain specialization of the dichromatic polynomial of the corresponding graph. In Appendix \ref{ap2} we discuss the most general diagrammatics based on the state-sum relation and we show how the models of the HOMFLYPT polynomial of Sections \ref{homflynot} and \ref{coment} naturally arise.\\

\indent Most of this thesis, specif\mbox{}ically Chapters \ref{graf}-\ref{torus}, is based on papers by the author, see \cite{moj1}, \cite{moj2}, \cite{moj3}, \cite{moj4}.

\chapter{Notation and Def\mbox{}initions}\label{not}

\section{Knots and links}

By a knot we mean a smoothly embedded circle $S^1$ in  3-space $\R^3$ (or $S^3$). By a link we mean a f\mbox{}inite disjoint collection of knots. We will identify  isotopic knots (and links), i.e. we will identify those that  can be continuously transformed  into each other. Every knot that is isotopic to the standard circle $S^1$ we will call the unknot.\\
\indent By a knot invariant we mean any function on the set of knots which has the same value on isotopic knots, and analogously for a link invariant.\\
\indent To every link $L$ we can associate (inf\mbox{}initely many, of course) planar projections $D$, which are also called regular diagrams (or just diagrams). In fact, regular diagrams are planar 4-valent graphs with slightly modif\mbox{}ied vertices such that one can distinguish between the over- and undercrossing. These modif\mbox{}ied vertices we call the crossings of the diagrams $D$ and we denote the set of all crossings of $D$ by $c(D)$.\\
\indent The famous theorem of Reidemeister \cite{rei}, states that two diagrams $D_1$ and $D_2$ represent two isotopic links if and only if one can be transformed into the other by a f\mbox{}inite sequence of Reidemeister moves R1, R2 and R3:

\begin{center}
\setlength{\unitlength}{4mm}
\begin{picture}(25,2)
\linethickness{0.7pt}
\put(0,0.8){R1:}
\put(-5,0){
\put(8,0){\LPP}
\put(10.8,0.8){$\longleftrightarrow$}
\put(14.5,0){\qbezier(0,0)(0.5,1)(0,2)}
} 
\put(8,0){
\put(4,0.8){and}
\put(8,0){\LPM}
\put(10.8,0.8){$\longleftrightarrow$}
\put(14.5,0){\qbezier(0,0)(0.5,1)(0,2)}
}

\end{picture}
\end{center}

\begin{center}
\setlength{\unitlength}{4mm}
\begin{picture}(25,4)
\linethickness{0.7pt}
\put(0,1.8){R2:}

\put(8,0){\X}
\put(8,2){\Y}
\put(10.5,1.8){ {$\longleftrightarrow$}}
\put(14.5,0){\qbezier(0,0)(0.5,2)(0,4)\qbezier(1,0)(0.5,2)(1,4)}

\end{picture}
\end{center}

\begin{center}
\setlength{\unitlength}{4mm}
\begin{picture}(25,4)
\linethickness{0.7pt}
\put(0,1.8){R3:}

\put(5,0){
\qbezier(0,0)(2,2)(4,4)
\qbezier(4,0)(3.1,0.9)(2.2,1.8)
\qbezier(1.2,2.8)(1.5,2.5)(1.8,2.2)
\qbezier(0,4)(0.4,3.6)(0.8,3.2)
\qbezier(2,4)(1.4,3.4)(0.8,2.8)
\qbezier(0.8,2.8)(0,2)(0.8,1.2)
\qbezier(2,0)(1.6,0.4)(1.2,0.8)
}
\put(10.5,1.8){ {$\longleftrightarrow$}}
\put(14,0){
\bezier{100}(0,0)(2,2)(4,4)
%\linethickness{1pt}
%\put(0,0){\line(1,1){4}}
\linethickness{0.7pt}
\bezier{100}(4,0)(3.6,0.4)(3.2,0.8)
\bezier{100}(2.8,1.2)(2.5,1.5)(2.2,1.8)
\bezier{100}(0,4)(0.9,3.1)(1.8,2.2)
\bezier{100}(2,0)(2.6,0.6)(3.2,1.2)
\bezier{100}(3.2,2.8)(4,2)(3.2,1.2)
\bezier{100}(2,4)(2.4,3.6)(2.8,3.2)
}

\end{picture}
\end{center}

Here we mean that two diagrams are related by the Reidemeister move if they look the same except locally where they are like the left and right picture, respectively, of any of the above moves. Hence, a function def\mbox{}ined on the set of planar diagrams will be a link invariant if and only if it is invariant under the Reidemeister moves.\\
\indent We can also introduce  an orientation on  knots and links. Then locally each crossing of $D$ will look like one of the following two pictures
{{\begin{center}
\setlength{\unitlength}{4mm}
\begin{picture}(25,6)
\linethickness{0.8pt}
\qbezier(5.00,5.00)(2.50,2.50)(0.00,0.00)
\qbezier(5.00,0.00)(4.00,1.00)(3.00,2.00)
\qbezier(0.00,5.00)(1.00,4.00)(2.00,3.00)
\qbezier(0.00,5.00)(0.25, 5.00)(0.50,5.00)
\qbezier(0.00,5.00)(0.00,4.75)(0.00,4.50)
\qbezier(5.00,5.00)(4.75, 5.00)(4.50,5.00)
\qbezier(5.00,5.00)(5.00,4.75)(5.00,4.50)
\qbezier(22.00,2.00)(21.00,1.00)(20.00,0.00)
\qbezier(25.00,5.00)(24.00,4.00)(23.00,3.00)
\qbezier(20.00,5.00)(22.50,2.50)(25.00,0.00)
\qbezier(20.00,5.00)(20.25,5.00)(20.50,5.00)
\qbezier(20.00,5.00)(20.00,4.75)(20.00,4.50)
\qbezier(25.00,5.00)(24.75,5.00)(24.50,5.00)
\qbezier(25.00,5.00)(25.00,4.75)(25.00,4.50)
\put(1.00,-1.00){{\small\textrm{positive}}}
\put(21.00,-1.00){{\small\textrm{negative}}}
%\put(10.00,-2.00){{\small\textrm{f\mbox{}igure 1}}}
\end{picture}
\end{center}
}}
\bigskip
If it looks like the picture on the left we call it a positive crossing and if it looks like the picture on the right we call it a negative crossing. The number of positive (respectively negative) crossings of the diagram $D$ we denote by $n_+(D)$ (respectively $n_-(D)$). We say that a link is positive if it has a planar projection whose crossings are all positive. \\
\indent In this thesis we will deal with knot invariants, $I$, that satisfy a skein relation. A skein relation is a linear relation between the values of the invariant $I$ on  three oriented links $L_+$, $L_-$ and $L_0$ which have identical planar projections except near one crossing where they are as in the following picture:
\begin{center}
\setlength{\unitlength}{3.5mm}
\begin{picture}(30,5.5)
\linethickness{0.8pt}
\qbezier(5.00,5.00)(2.50,2.50)(0.00,0.00)
\qbezier(5.00,0.00)(4.00,1.00)(3.00,2.00)
\qbezier(0.00,5.00)(1.00,4.00)(2.00,3.00)
\qbezier(0.00,5.00)(0.25, 5.00)(0.50,5.00)
\qbezier(0.00,5.00)(0.00,4.75)(0.00,4.50)
\qbezier(5.00,5.00)(4.75, 5.00)(4.50,5.00)
\qbezier(5.00,5.00)(5.00,4.75)(5.00,4.50)

\put(12,0){\qbezier(2.00,2.00)(1.00,1.00)(0.00,0.00)
\qbezier(5.00,5.00)(4.00,4.00)(3.00,3.00)
\qbezier(0.00,5.00)(2.50,2.50)(5.00,0.00)
\qbezier(0.00,5.00)(0.25,5.00)(0.50,5.00)
\qbezier(0.00,5.00)(0.00,4.75)(0.00,4.50)
\qbezier(5.00,5.00)(4.75,5.00)(4.50,5.00)
\qbezier(5.00,5.00)(5.00,4.75)(5.00,4.50)}

\put(24,0)
{\qbezier(0,0)(2,2.5)(0,5)
\qbezier(5,0)(3,2.5)(5,5)
\qbezier(0.00,5.00)(0.25,5.00)(0.50,5.00)
\qbezier(0.00,5.00)(0.00,4.75)(0.00,4.50)
\qbezier(5.00,5.00)(4.75,5.00)(4.50,5.00)
\qbezier(5.00,5.00)(5.00,4.75)(5.00,4.50)}

\put(2.00,-2.00){$L_+$}
\put(14.00,-2.00){$L_-$}
\put(26.00,-2.00){$L_0$}
\end{picture}
\end{center}

\bigskip

\section{Braids}\label{braidnot}
 
Braids are a special type of open links. An $n$ strand braid is  obtained by f\mbox{}ixing $n$ points in the plane and the same $n$ points in a parallel copy of the same plane below the f\mbox{}irst one. A braid is obtained by attaching the top $n$ points to the bottom $n$ points by strings (later called strands) in such a way that the strings never intersect each other and such that the tangent vector to a string is never parallel to the planes. We identify  braids up to isotopy: we only assume that the $2n$ endpoints are f\mbox{}ixed and the tangent vectors satisfy the above properties through the isotopy.\\ 
% intermediate steps we have $n$ braids. \\
\indent The most useful thing about braids is that they form groups under obvious concatenation. We denote the group of $n$-strand braids by $B_n$.\\
\indent The braid group may be converted into a purely algebraic object by means of the  Artin presentation (see \cite{ar}). The generators of the group $B_n$ are $\sigma_1,\ldots,\sigma_{n-1}$  where $\sigma_i$ is:    

{{
\begin{center}
\setlength{\unitlength}{4mm}
\begin{picture}(30,6) 
\linethickness{0.8pt}
%\line(1.00,2.00,1.00,3.00)
\qbezier(5.00,1.00)(5.00,3.00)(5.00,5.00)
\qbezier(9.00,1.00)(9.00,3.00)(9.00,5.00)

\qbezier(15.00,5.00)(13.00,3.00)(11.00,1.00)
\qbezier(15.00,1.00)(14.30,1.70)(13.60,2.40)
\qbezier(11.00,5.00)(11.70,4.30)(12.40,3.60)

\qbezier(17.00,1.00)(17.00,3.00)(17.00,5.00)
\qbezier(21.00,1.00)(21.00,3.00)(21.00,5.00)

\put(6.00,3.00){$\cdot$}
\put(7.00,3.00){$\cdot$}
\put(8.00,3.00){$\cdot$}
\put(18.00,3.00){$\cdot$}
\put(19.00,3.00){$\cdot$}
\put(20.00,3.00){$\cdot$}

\put(5.00,5.50){1}
\put(11.00,5.50){$i$}
\put(14.50,5.50){$i+1$}
\put(21.00,5.50){$p$}
\put(12.80,0.00){{\Large $\sigma_i$}}

\end{picture}
\end{center}
}}

Then the presentation of the braid group is given by the relations
\begin{eqnarray*}
\sigma_i\sigma_j&=&\sigma_j\sigma_i,\quad |i-j|>1\\
\sigma_i\sigma_{i+1}\sigma_i&=&\sigma_{i+1}\sigma_i\sigma_{i+1}
\end{eqnarray*}

\indent We will identify the braid itself and the braid word (word in the alphabet $\sigma_1,\sigma^{-1}_1,\ldots,\sigma_{n-1},\sigma_{n-1}^{-1}$) that represents it.\\
\indent  We will assume that all strands of a braid are oriented from bottom to top. From any braid, $B$, we can obtain a knot (or link) by closing the strands i.e. by connecting the initial points to the endpoints by means of a collection of parallel strands:

{{
\begin{center}
\setlength{\unitlength}{4mm}
\begin{picture}(25,6)
\linethickness{0.7pt}
\put(7,1.8){\qbezier(0,0)(0,1.2)(0,2.4)
\qbezier(2.4,0)(2.4,1.2)(2.4,2.4)
\qbezier(0,0)(1.2,0)(2.4,0)
\qbezier(0,2.4)(1.2,2.4)(2.4,2.4)
\put(0.7,0.7){{\Large{$B$}}}

\qbezier(0.6,0)(0.6,-0.3)(0.6,-0.6)
\put(0.6,0){\qbezier(0.6,0)(0.6,-0.3)(0.6,-0.6)}
\put(1.2,0){\qbezier(0.6,0)(0.6,-0.3)(0.6,-0.6)}
\put(0,3){\qbezier(0.6,0)(0.6,-0.3)(0.6,-0.6)}
\put(1.2,3){\qbezier(0.6,0)(0.6,-0.3)(0.6,-0.6)}
\put(0.6,3){\qbezier(0.6,0)(0.6,-0.3)(0.6,-0.6)}

%\put(3,-0.6){\qbezier(0,0)(0,1.8)(0,3.6)}
%\put(3.6,-0.6){\qbezier(0,0)(0,1.8)(0,3.6)}
%\put(4.2,-0.6){\qbezier(0,0)(0,1.8)(0,3.6)}

%\qbezier(0.6,-0.6)(2.4,-3.5)(4.2,-0.6)
%\qbezier(1.2,-0.6)(2.4,-2.5)(3.6,-0.6)
%\qbezier(1.8,-0.6)(2.4,-1.5)(3,-0.6)

%\qbezier(0.6,3)(2.4,5.9)(4.2,3)
%\qbezier(1.2,3)(2.4,4.9)(3.6,3)
%\qbezier(1.8,3)(2.4,3.9)(3,3)
}

\put(11,3){$\longrightarrow$}

\put(14,1.8){\qbezier(0,0)(0,1.2)(0,2.4)
\qbezier(2.4,0)(2.4,1.2)(2.4,2.4)
\qbezier(0,0)(1.2,0)(2.4,0)
\qbezier(0,2.4)(1.2,2.4)(2.4,2.4)
\put(0.7,0.7){{\Large{$B$}}}

\qbezier(0.6,0)(0.6,-0.3)(0.6,-0.6)
\put(0.6,0){\qbezier(0.6,0)(0.6,-0.3)(0.6,-0.6)}
\put(1.2,0){\qbezier(0.6,0)(0.6,-0.3)(0.6,-0.6)}
\put(0,3){\qbezier(0.6,0)(0.6,-0.3)(0.6,-0.6)}
\put(1.2,3){\qbezier(0.6,0)(0.6,-0.3)(0.6,-0.6)}
\put(0.6,3){\qbezier(0.6,0)(0.6,-0.3)(0.6,-0.6)}

\put(3,-0.6){\qbezier(0,0)(0,1.8)(0,3.6)}
\put(3.6,-0.6){\qbezier(0,0)(0,1.8)(0,3.6)}
\put(4.2,-0.6){\qbezier(0,0)(0,1.8)(0,3.6)}

\qbezier(0.6,-0.6)(2.4,-3.5)(4.2,-0.6)
\qbezier(1.2,-0.6)(2.4,-2.5)(3.6,-0.6)
\qbezier(1.8,-0.6)(2.4,-1.5)(3,-0.6)

\qbezier(0.6,3)(2.4,5.9)(4.2,3)
\qbezier(1.2,3)(2.4,4.9)(3.6,3)
\qbezier(1.8,3)(2.4,3.9)(3,3)
}
\end{picture}
\end{center}
}}

Conversely, it was shown by Alexander that every link or knot can be put in the form of a closed braid, i.e. it is isotopic to a link which is the closure of a braid (\cite{ax2}). We say that a knot (or link) is a positive braid knot if it has a planar projection which is the closure of a positive braid, i.e. a braid whose crossings are all positive.  

\section{The Jones polynomial}\label{jonesnot}

In this section we recall the def\mbox{}inition of the Jones polynomial. It was def\mbox{}ined by V.F.R. Jones in 1985 in \cite{jones}. He obtained the polynomial by introducing a certain tower of von Neumann algebras and by constructing a certain representation into the Artin braid group. The Jones polynomial $J(L)$ is a Laurent polynomial with values in $\Z[t^{1/2},t^{-1/2}]$ satisfying the following skein relation
\begin{equation}
tJ(L_+)-t^{-1}J(L_-)=(t^{1/2}-t^{-1/2})J(L_0),
\label{jonsk}
\end{equation}
and such that the value of the unknot is 1.\\
\indent In 1987, L. Kauf\mbox{}fman in \cite{kaufss} introduced a state-sum model of the link $L$, and obtained an alternative def\mbox{}inition of the Jones polynomial, which is remarkably simple and purely combinatorial. We will describe Kauf\mbox{}fman's state-sum model. In order to do that, we will use slightly dif\mbox{}ferent variables than the original ones of Kauf\mbox{}fman. Specif\mbox{}ically, we will use the ones introduced by M.Khovanov in \cite{kov} which are now widely used in  link homology theory.\\
\indent Let $L$ be a link, and let $D$ be a planar projection of $L$. Then for each crossing $c\in D$ we will def\mbox{}ine two resolutions (or states) of $c$ according to the following picture:
{{
\begin{center}
\setlength{\unitlength}{4mm}
\begin{picture}(25,6)
\linethickness{0.8pt}
\qbezier(1.00,1.00)(1.00,3.00)(1.00,5.00)
\qbezier(3.00,1.00)(3.00,3.00)(3.00,5.00)

\qbezier(5.00,3.00)(7.00,3.00)(9.00,3.00)
\qbezier(5.00,3.00)(5.25,3.25)(5.50,3.50)
\qbezier(5.00,3.00)(5.25,2.75)(5.50,2.50)

\qbezier(15.00,5.00)(13.00,3.00)(11.00,1.00)
\qbezier(15.00,1.00)(14.30,1.70)(13.60,2.40)
\qbezier(11.00,5.00)(11.70,4.30)(12.40,3.60)

\qbezier(17.00,3.00)(19.00,3.00)(21.00,3.00)
\qbezier(21.00,3.00)(20.75,3.25)(20.50,3.50)
\qbezier(21.00,3.00)(20.75,2.75)(20.50,2.50)

\qbezier(23.00,5.00)(24.00,2.00)(25.00,5.00)
\qbezier(23.00,1.00)(24.00,4.00)(25.00,1.00)

\put(5.00,3.75){{\small\textrm{0-resolution}}}
\put(17.00,3.75){{\small\textrm{1-resolution}}}

\end{picture}
\end{center}
}}

We call a collection of circles that are obtained by resolving all crossings of $D$, a total resolution of the diagram $D$. We order the set of crossings of $D$, say from 1 to $n$, where $n$ is the number of crossings of $D$. Then, 
obviously, we have a bijection between the set of the total resolutions of the diagram $D$ and the set $\{0,1\}^{n}$, i.e. the set of all $n$-tuples of zeros and ones.
To each element $\epsilon=(\epsilon_1,\ldots,\epsilon_n)\in\{0,1\}^{n}$ we assign the resolution $D_{\epsilon}$ which is obtained by resolving the $i$-th crossing as an $\epsilon_i$-resolution. We denote by $|\epsilon|$ the sum of elements of $\epsilon$, and by $c(\epsilon)$ the number of disjoint circles of the resolution $D_{\epsilon}$.  Now we def\mbox{}ine the Kauf\mbox{}fman bracket of the diagram $D$ by:
\begin{equation}\langle D\rangle=\sum_{\epsilon\in\{0,1\}^n} {(-1)^{|\epsilon|}q^{|\epsilon|}(q+q^{-1})^{c(\epsilon)}}.\label{kaufbr}\end{equation}

Alternatively, we can describe the Kauf\mbox{}fman bracket by the following recursive axioms 
\begin{center}
\setlength{\unitlength}{4mm}
\begin{picture}(25,2)
\linethickness{0.7pt}
\put(7.4,0.8){$\langle$}
\put(8,0){\X}
\put(9.4,0.8){$\rangle$}
\put(10,0.8){=}
\put(11.2,0.8){$\langle$}
\put(12,0){\O}
\put(13.4,0.8){$\rangle$}
\put(14,0.8){$-q$}
\put(15.4,0.8){$\langle$}
\put(16,0){\l}
\put(17.4,0.8){$\rangle$}
\end{picture}
\end{center}
%\begin{eqnarray*}
%<D>&=&<D_0>-q<D_1>\\
$$\langle U_k \rangle=(q+q^{-1})^k,$$
%\end{eqnarray*}
where $U_k$ is the collection of $k$ disjoint circles in the plane. Here in the f\mbox{}irst formula we assume that the relation is valid for any three diagrams that are the same except near one crossing where they look like the picture. 
It can be easily verif\mbox{}ied that the Kauf\mbox{}fman bracket is invariant under the R3 move, and  that under the R1 and R2 moves it has the following simple behaviour

\begin{center}
\setlength{\unitlength}{4mm}
\begin{picture}(25,2)
\linethickness{0.7pt}
\put(7.4,0.8){$\langle$}
\put(8,0){\LPP}
\put(9.8,0.8){$\rangle$}
\put(10.8,0.8){=}
\put(12,0.8){$q^{-1}$}
\put(14,0.8){$\langle$}
\put(15,0){\qbezier(0,0)(0.5,1)(0,2)}
\put(16,0.8){$\rangle$}
\end{picture}
\end{center}

\begin{center}
\setlength{\unitlength}{4mm}
\begin{picture}(25,2)
\linethickness{0.7pt}
\put(7.4,0.8){$\langle$}
\put(8,0){\LPM}
\put(9.8,0.8){$\rangle$}
\put(10.8,0.8){=}
\put(12,0.8){$-q^{2}$}
\put(14,0.8){$\langle$}
\put(15,0){\qbezier(0,0)(0.5,1)(0,2)}
\put(16,0.8){$\rangle$}
\end{picture}
\end{center}

\begin{center}
\setlength{\unitlength}{4mm}
\begin{picture}(25,4)
\linethickness{0.7pt}
\put(7.2,1.8){$\langle$}
\put(8,0){\X}
\put(8,2){\Y}
\put(9.6,1.8){$\rangle$}
\put(10.5,1.8){ {$=$}}
\put(12,1.8){$-q$}
\put(13.8,1.8){$\langle$}
\put(14.5,0){\qbezier(0,0)(0.5,2)(0,4)\qbezier(1,0)(0.5,2)(1,4)}
\put(16,1.8){$\rangle$}

\end{picture}
\end{center}

%$$<L_{+}>=q^{-1}<L>$$
%$$<L_{-}>=-q^{2}<L>$$
%$$<L_{00}>=-q<L>.$$
%Hence, if by $n_{+}$ and $n_-$ we denote the number of positive and %negative crossings of diagram $D$, respectively, 
Then the polynomial $\hat{J}(D)$ def\mbox{}ined by 
$$ \hat{J}(D)=(-1)^{n_-}q^{n_+-2n_-}\langle D \rangle,$$
is obviously invariant under all Reidemeister moves. Hence, we can write $\hat{J}(L)$, and we call it the unnormalized Jones polynomial. From this def\mbox{}inition, we easily obtain the skein relation:
$$q^{-2}\hat{J}(L_+)-q^{2}\hat{J}(L_-)=(q^{-1}-q)\hat{J}(L_0),$$
and the value of the unknot is $q+q^{-1}$. Thus, if we def\mbox{}ine $J(L)=\hat{J}(L)/(q+q^{-1})$ then $J(L)$ is the Jones polynomial (with the variable $t$ from (\ref{jonsk}) replaced by $q^{-2}$).  
 
\section{Khovanov homology ($sl(2)$ link homology)}\label{kovnot}

Here we recall  the def\mbox{}inition of Khovanov homology for links (see \cite{bn},\cite{kov}). The idea of Khovanov homology is to def\mbox{}ine a graded chain complex whose graded Euler characteristic is equal to the Jones polynomial. This is done by following the lines of the Kauf\mbox{}fman bracket construction and by ``categorifying" the polynomials involved, i.e. by replacing the polynomials with graded modules (whose graded dimensions are equal to the polynomials), sums with direct sums, products with tensor products, etc. We f\mbox{}irst recall the basic properties of graded $\Z$-modules:
\subsection{Graded dimension of a graded $\Z$-module}

\begin{definition}
Let $M=\oplus_{k}M_k$ for $k \in \Z$ be a graded $\Z$-module where $\{M_k\}$ denotes the $k$-th graded component of $M$. The graded (or quantum) dimension of $M$ is the power series
$$q\dim M:=\sum_k{q^k\rank(M_k)},$$
where $\rank(M_k)=\dim_{\Q}(M_k\otimes\Q)$.
\end{definition}

\indent The direct sum and tensor product can be def\mbox{}ined in the graded category in an obvious way. The following proposition is straightforward.
\begin{proposition}
Let $M$ and $N$ be graded $\Z$-modules. Then
\begin{eqnarray*}
q\dim (M\oplus N)&=&q\dim(M)+q\dim(N) \\
q\dim(M\otimes N)&=&q\dim(M) q\dim(N).
\end{eqnarray*}
\end{proposition}
\begin{definition}
Let $\{ l\}$, $l \in \Z$, be the ``degree shift" operation on graded $\Z$-modules. That is, if $M=\oplus_k M_k$ is a graded $\Z$-module where $M_k$ denotes the $k$-th graded component of $M$, we set $M\{l\}_k:=M_{k-l}$. Then we have 
$q\dim M\{l\}=q^{l}q\dim M.$ In other words, all the degrees are increased by $l$.
\end{definition}

\indent Now  go back to the construction of Khovanov homology. Let $L$ be a link, and $D$ its planar projection. Take an ordering of the crossings of $D$. 
Denote by $n$ the number of crossings of  $D$. Then, as before, there is a  bijective correspondence between the total resolutions of $D$ and the set $\{0,1\}^n$.\\
\indent
Now we organize the Kauf\mbox{}fman bracket state-sum expression (\ref{kaufbr}) in the following way: 
\begin{equation}\langle D\rangle=\sum_{i=0}^n{(-1)^iq^i}\sum_{|\epsilon|=i} {(q+q^{-1})^{c(\epsilon)}}.\label{}\end{equation}
\indent To each total resolution $D_{\epsilon}$ we want to assign a graded module whose quantum dimension is equal to $q^{|\epsilon|}(q+q^{-1})^{c(\epsilon)}$. Let $V$ be a  graded $\Z$-module $V$, which is freely generated by two basis vectors $1$ and $X$, with $\deg 1 = 1$ and $\deg X=-1$. Then to $D_{\epsilon}$ we assign  a graded module $V_{\epsilon}$, which is the tensor product of $V$'s over all circles in the resolution.  We draw all the resolutions as the vertices of a (skewed) $n$-dimensional cube such that in the $i$-th column we have the resolutions $D_{\epsilon}$ with $|\epsilon|=i$. We def\mbox{}ine the $i$-th chain group $C^i(D)$ by:
$$C^i(D)=\oplus_{|\epsilon|=i}V_{\epsilon}\{i\}.$$
\indent We obviously have that the quantum dimension of $C^i(D)$ is equal to $q^i\sum_{|\epsilon|=i}{(q+q^{-1})^{c(\epsilon)}}.$\\
\indent Now, we want to obtain the whole graded  chain   
complex, i.e. we want to def\mbox{}ine the dif\mbox{}ferential. First, we give the general def\mbox{}initions for the graded chain complexes.

\subsection{Graded chain complexes and the graded Euler characteristic}
\label{subs1}

\begin{definition}
Let  $M=\oplus_j M_j$ and $N=\oplus_j N_j$ be graded $\Z$-modules where $M_j$ and $N_j$ denote the $j$-th graded component of $M$ and $N$, respectively. A $\Z$-module map $\alpha:M\rightarrow N$ is said to be graded with degree $d$ if $\alpha(M_j)\subset N_{j+d}$, i.e. elements of degree $j$ are mapped to elements of degree $j+d$. A $\Z$-module map is called degree preserving (or grading preserving) if it is graded of degree zero.\\
A graded chain complex is a chain complex for which the chain groups are graded $\Z$-modules and the dif\mbox{}ferentials are graded.
\end{definition}

\begin{definition}
Let $[s]$, $s\in \Z$, be the ``height shift" operation on chain complexes. That is, if $C$ is a chain complex, $\cdots\to C^r\to C^{r+1}\cdots$ of modules (or vector spaces), and if $\bar{C}=C[s]$, then $\bar{C}^r=C^{r-s}$, with all dif\mbox{}ferentials shifted accordingly.
\end{definition}

\begin{definition}
The graded Euler characteristic $\chi_q(C)$ of a graded chain complex $C$ is the alternating sum of the graded dimensions of its homology groups, i.e. $\chi_q(C)=\sum_i{(-1)^i q\dim(H^i)}$.
\end{definition}

\begin{proposition}
(\cite{bn}, \cite{gr}) If the dif\mbox{}ferentials are degree preserving and all %$C^{i,j}$'s (the sets of f\mbox{}ixed degree of each chain group)
chain groups are f\mbox{}inite-dimensional, the graded Euler characteristic is also equal to the alternating sum of the graded dimensions of its chain groups, i.e.
$$\chi_q(C)=\sum_{i} (-1)^iq\dim(H^i)=\sum_{i}(-1)^iq\dim (C^i).$$
\label{p1}
\end{proposition}
\begin{remark} From the last Proposition we have that if we def\mbox{}ine a degree preserving dif\mbox{}ferential between the groups $C^i(D)$, the chain complex $C(D)$ obtained in this way will have as Euler characteristic the  Kauf\mbox{}fman bracket $\langle D\rangle$. Hence, we are left with def\mbox{}ining such a dif\mbox{}ferential.
\label{rem1}
\end{remark}
\indent The dif\mbox{}ferentials $d^i:C^i(D)\to C^{i+1}(D)$ are def\mbox{}ined as the (signed) sum of  ``per-edge" dif\mbox{}ferentials. Specif\mbox{}ically, the only nonzero maps are from $D_{\epsilon}$ to $D_{\epsilon'}$, where $\epsilon=(\epsilon_1,\ldots,\epsilon_n)$, $\epsilon_i \in \{0,1\}$, if and only if $\epsilon'$ has all entries the same as $\epsilon$ except one $\epsilon_j$, for some $j=1,\ldots,m$, which is changed from 0 to 1. We denote these dif\mbox{}ferentials by $d_{\nu}: V_{\epsilon}\{|\epsilon|\}\to V_{\epsilon'}\{|\epsilon'|\}$, where $\nu$ is the $n$-tuple which consists of the label $\ast$ at the position $j$ and of $n-1$ 0's and 1's (the same as the remaining entries of $\epsilon$).\\
\indent  
Note that in these cases, either two circles of $D_{\epsilon}$ merge into one circle of $D_{\epsilon'}$ or one circle of $D_{\epsilon}$ splits into two circles of $D_{\epsilon'}$, while all other circles remain the same. \\
\indent In the f\mbox{}irst case, we def\mbox{}ine the map $d_{\nu}$ as the identity on the tensor factors ($V$) that correspond to the unchanged circles, and on the remaining factors we def\mbox{}ine it to be the (grading preserving) multiplication map $m:V\otimes V \to V\{1\}$, which is given on basis vectors by:
$$m(1\otimes 1)=1,\quad m(1\otimes X)=m(X\otimes 1)=X,\quad m(X\otimes X)=0.$$
\indent In the second case, we def\mbox{}ine the map $d_{\nu}$ as the identity on the tensor factors ($V$) that correspond to the unchanged circles,
and on the remaining factors we def\mbox{}ine it to be the (grading preserving) comultiplication map $\Delta:V\to V \otimes V\{1\}$, which is given on basis vectors by:
$$\Delta(1)=1\otimes X+X\otimes 1,\quad \Delta(X)=X\otimes X.$$
\indent Finally, to obtain the dif\mbox{}ferential $d^i$ of the chain complex $C(D)$, we sum all contributions $d_{\nu}$ with $|\nu|=i$, multiplied by the sign $(-1)^{f(\nu)}$, where $f(\nu)$ is equal to the number of 1's ordered before $\ast$ in $\nu$. This makes every square of our cubic complex anticommutative, and hence we obtain a genuine dif\mbox{}ferential (i.e. $d^{i+1}d^i=0$).\\
\indent We denote the homology groups of the complex $(C(D),d)$ obtained
in this way by $H^i(D)$ and call them the \textit{unnormalized homology groups of $D$}. Since the chain groups are naturally graded by the construction of the  chain complex, so are the homology groups, and we write $H^{i,j}(D)$ for the summand of degree $j$ in $H^i(D)$. In order to obtain link invariants (i.e. independence of the chosen projection), we have to shift the chain complex (and hence the homology groups) by:
\begin{equation}
\C(D)=C(D)[-n_{-}]\{n_+-2n_-\},
\label{eq1}
\end{equation}
where $n_+$ and $n_-$ are the number of positive and negative crossings, respectively, of the diagram $D$.  The homology groups of the shifted complex $\C(D)$ we denote by $\H^i(D)$. Hence, we have $\H^{i,j}(D)=H^{i+n_-,j-n_++2n_-}(D)$. 
\begin{theorem}(\cite{kov},\cite{bn})
The homology groups $\H^i(D)$ are independent of the choice of the planar projection $D$. Furthermore, the graded Euler characteristic of the complex $\C(D)$ is equal to the Jones polynomial of the link $L$.
\end{theorem}
\indent Hence, we can write $\H(L)$, and we call $\H^i(L)$ \textit{the homology groups of the link $L$}. 

\subsection{Basic properties of Khovanov homology}

 \indent If the diagram $D$ has only positive crossings then we do not have an overall shift in the homology degrees, and so we have that, for example, $\H^1(L)$ is trivial if and only if $H^{1}(D)$ is trivial. Also, in the general case if we have a positive knot $K$ (a knot that has a planar  projection with only positive crossings) then we have that $\H^i(K)$ is trivial for all $i<0$. Furthermore, if $D$ is a planar projection of a positive knot $K$ with $n_-$ negative crossings then $H^i(D)$ is trivial for $i<n_-$. \\

\indent Usually, the homology groups of the link $K$ are represented as a planar array in such a way that  $\rank \H^{i,j}(K)=\dim (\H^{i,j}(K)\otimes \Q)$ (or the whole group $\H^{i,j}(K)$, if we want to keep track of the torsions) is specif\mbox{}ied in the position $(i,j)$. As can be easily seen, the $q$-gradings ($j$) of the generators of nontrivial $\H^{i,j}(K)$ and the number of components of $K$ are of the same parity (either all are even or all are odd). Hence, by a \textit{diagonal} of the homology of the link $K$, we mean a line $j-2i=a=const$, when there exist integers $i$ and $j$ such that $j-2i=a$ and $\rank \H^{i,j}(K)>0$. If $a_{\max}$ and $a_{\min}$ are the maximal and minimal value of $a$ such that the line $j-2i=a$ is a diagonal of the homology of the link $K$, then we def\mbox{}ine the homological width of the link $K$ to be  $h(K)=(a_{\max}-a_{\min})/2+1$.\\
\indent Every knot (link) occupies at least two diagonals (i.e. $h(K) \ge 2$ for every link $K$), and the ones that occupy exactly two diagonals are called H-thin, or homologically thin. For example all alternating knots are H-thin (\cite{lee}), and the free part of the homology of any H-thin knot is determined by its Jones polynomial and the signature. A knot that is not H-thin is called H-thick or homologically thick.\\
\indent Furthermore, for an element $x\in \H^{i,j}(K)$, we denote its homological grading -- $i$ -- by $t(x)$, and its $q$-grading (also called the quantum grading) -- $j$ -- by $q(x)$. We also introduce a third grading $\delta(x)$ by $\delta(x)=q(x)-2t(x)$. Hence, we have that the knot $K$ is H-thick, if there exist three generators of $\H(K)$ with dif\mbox{}ferent values of the  $\delta$-grading.\\
\indent An alternative way of presenting the homology of the knot is by means of the two-variable Poincar\'e polynomial $P(K)(t,q)$ of the chain complex $\C(D)$, i.e.:
$$P(K)(t,q)=\sum_{i,j\in \Z}{t^iq^j \rank \H^{i,j}(K)}.$$
%\indent We also def\mbox{}ine the unnormalized Poincare polynomial $P(D)(t,q)$ %of the diagram $D$ as the Poincare polynomial of the chain complex %$C(D)$, i.e.
%$$P(K)(t,q)=t^{n_-}q^{-n_++2n_-}\P(K)(t,q),$$
%where $n_+$ and $n_-$ are the numbers of positive and negative %crossings, respectively, of the diagram $D$.\\

\indent Let $D$ be a diagram of a link $L$ and let $c$ be one of its  crossings. Denote by $D_i$, $i=0,1$ the diagram that is obtained after performing an $i$-resolution of the crossing $c$. Then there exists a  long exact sequence of (unnormalized) homology groups (see e.g. \cite{v} and Remark \ref{rem2} in Chapter \ref{graf}):
{\small{
\begin{equation}\cdots\rightarrow H^{i-1,j-1}(D_1)\to H^{i,j}(D)\to H^{i,j}(D_0)\to H^{i,j-1}(D_1)\to H^{i+1,j}(D)\to \cdots\label{eq2}\end{equation}
}}

\section{The HOMFLYPT polynomial}\label{homflynot}

The HOMFLYPT polynomial $P(L)(a,q)$, named after its many discoverers (\cite{homfly}, \cite{pt}), is a two-variable generalization of the Jones polynomial. It is given by the following skein relation:
\begin{equation}
aP(L_+)-a^{-1}P(L_-)=(q-q^{-1})P(L_0),
\label{homflysk}
\end{equation}
and by prescribing the value of the unknot to be equal to 1. As can be easily shown (\cite{mus}), the skein relation (\ref{homflysk}) is the most general skein relation of this form. \\
\indent There is an alternative description of the HOMFLYPT polynomial, given by def\mbox{}ining an inf\mbox{}inite set of one-variable specializations, by $P_n(L)(q)$ $=$\- $P(L)(q^n,q)$, one for each positive integer $n$, which is enough to recover the whole HOMFLYPT polynomial. Obviously, the case $n=2$, i.e. the polynomial $P_2(L)(q)$ is equal to the Jones polynomial $J(L)(q^{-1})$.

\subsection{A state-sum model}

In order to categorify the HOMFLYPT polynomial, we need an analogous state-sum model for its $n$-specializations so that we can proceed as in the case of $sl(2)$-link homology. Such a model was introduced by Murakami, Ohtsuki and Yamada in \cite{moy} (in greater generality), and also previously by Kauf\mbox{}fman (\cite{kauf}, section I.11). \\
\indent Let $L$ be an oriented link, and let $D$ be a planar projection of $L$. Then to each crossing $c$ of $D$ we associate 0- and 1-resolutions, depending on the sign of the crossing $c$, according to the following picture:

{{
\begin{equation}
\setlength{\unitlength}{4mm}
\begin{picture}(25,4.5) 
\linethickness{0.6pt}
%\line(1.00,2.00,1.00,3.00)
\put(0,-1){
\qbezier(4.00,1.00)(4.00,3.00)(4.00,5.00)
\qbezier(3.80,4.60)(3.90,4.80)(4.00,5.00)
\qbezier(4.20,4.60)(4.10,4.80)(4.00,5.00) \qbezier(2.00,1.00)(2.00,3.00)(2.00,5.00)
\qbezier(1.80,4.60)(1.90,4.80)(2.00,5.00)
\qbezier(2.20,4.60)(2.10,4.80)(2.00,5.00) 
\qbezier(9.5,3)(10.5,3)(11.5,3)
\qbezier(11.5,3)(11.3,3.1)(11.1,3.2)
\qbezier(11.5,3)(11.3,2.9)(11.1,2.8)

\qbezier(4.5,3)(5.5,3)(6.5,3)
\qbezier(4.5,3)(4.7,3.1)(4.9,3.2)
\qbezier(4.5,3)(4.7,2.9)(4.9,2.8)

\put(10,0){
\qbezier(9,1)(8,3)(7,5)
\qbezier(7,1)(7.4,1.8)(7.8,2.6)
\qbezier(9,5)(8.6,4.2)(8.2,3.4)
}

\linethickness{2.5pt}
\qbezier(13.00,2.00)(13.00,3.00)(13.00,4.00)
\linethickness{0.6pt}
\qbezier(12.00,1.00)(12.50,1.50)(13.00,2.00)
\qbezier(14.00,1.00)(13.50,1.50)(13.00,2.00)
\qbezier(12.00,5.00)(12.50,4.50)(13.00,4.00)
\qbezier(14.00,5.00)(13.50,4.50)(13.00,4.00)

\qbezier(14.00,5.00)(13.80,4.90)(13.60,4.80)
\qbezier(14.00,5.00)(13.90,4.80)(13.80,4.60)

\qbezier(12.00,5.00)(12.20,4.90)(12.40,4.80)
\qbezier(12.00,5.00)(12.10,4.80)(12.20,4.60)

\qbezier(12.50,1.50)(12.30,1.40)(12.10,1.30)\qbezier(12.50,1.50)(12.40,1.30)(12.30,1.10)

\qbezier(13.50,1.50)(13.70,1.40)(13.90,1.30)\qbezier(13.50,1.50)(13.60,1.30)(13.70,1.10)
\qbezier(14.5,3)(15.5,3)(16.5,3)
\qbezier(14.5,3)(14.7,3.1)(14.9,3.2)
\qbezier(14.5,3)(14.7,2.9)(14.9,2.8)

\qbezier(19.5,3)(20.5,3)(21.5,3)
\qbezier(21.5,3)(21.3,3.1)(21.1,3.2)
\qbezier(21.5,3)(21.3,2.9)(21.1,2.8)
\put(-10,0){
\qbezier(19,5)(18,3)(17,1)
\qbezier(17,5)(17.4,4.2)(17.8,3.4)
\qbezier(19,1)(18.6,1.8)(18.2,2.6)

\qbezier(19,5)(19,4.75)(19,4.5)
\qbezier(19,5)(18.8,4.8)(18.6,4.6)
\qbezier(17,5)(17,4.75)(17,4.5)
\qbezier(17,5)(17.2,4.8)(17.4,4.6)
}
\put(10,0){
\qbezier(9,5)(9,4.75)(9,4.5)
\qbezier(9,5)(8.8,4.8)(8.6,4.6)
\qbezier(7,5)(7,4.75)(7,4.5)
\qbezier(7,5)(7.2,4.8)(7.4,4.6)
}

\qbezier(22.00,1.00)(22.00,3.00)(22.00,5.00)
\qbezier(21.80,4.60)(21.90,4.80)(22.00,5.00)
\qbezier(22.20,4.60)(22.10,4.80)(22.00,5.00)

\qbezier(24.00,1.00)(24.00,3.00)(24.00,5.00)
\qbezier(23.80,4.60)(23.90,4.80)(24.00,5.00)
\qbezier(24.20,4.60)(24.10,4.80)(24.00,5.00)

\put(5.40,3.50){0}
\put(10.3,3.50){1}
\put(15.40,3.50){0}
\put(20.3,3.50){1}
}
\end{picture}
\label{rezsln}
\end{equation}
}}

\indent Again we have a bijection between the set of all $n$-tuples $\epsilon$  of 0's and 1's and the set of total resolutions $D_{\epsilon}$. 
The total resolutions  in this case consist of the planar trivalent graphs $\Gamma$, such that at each vertex we have exactly one (unoriented) wide edge and two oriented thin edges, and such that at one end of each wide edge both thin edges are incoming and at the other end both thin edges are outgoing. \\
\indent Now, introduce the function $\hat{P}_n(\Gamma)(q)\in \Z[q,q^{-1}]$ on the set of such planar trivalent graphs, by the following f\mbox{}ive axioms:

{\Large{
\begin{center}
\setlength{\unitlength}{4mm}
\begin{picture}(25,3) 
\linethickness{0.6pt}
\put(10,1.5){\oval(4,3)}
\put(10,3){\qbezier(0,0)(0.15,0.15)(0.3,0.3)\qbezier(0,0)(0.15,-0.15)(0.3,-0.3)}
\put(13,1.1){$=\quad[n]$}

\end{picture}
\end{center}
}}
{\Large{
\begin{center}
\setlength{\unitlength}{5mm}
\begin{picture}(35,3) 
\linethickness{0.6pt}
\put(2,0){
\linethickness{2.5pt}\qbezier(0,0.5)(0,0.75)(0,1)\qbezier(0,3)(0,3.25)(0,3.5)
\linethickness{0.6pt}

\put(0,0.1){\qbezier(-0.5,2)(-0.6,1.85)(-0.7,1.7)
\qbezier(-0.5,2)(-0.4,1.85)(-0.3,1.7)}
\put(0,0.1){\qbezier(0.5,2)(0.6,1.85)(0.7,1.7)
\qbezier(0.5,2)(0.4,1.85)(0.3,1.7)}

\qbezier(0,1)(-0.5,1.6)(-0.5,2)
\qbezier(0,3)(-0.5,2.4)(-0.5,2)
\qbezier(0,1)(0.5,1.6)(0.5,2)
\qbezier(0,3)(0.5,2.4)(0.5,2)

%\qbezier(0,3)(0.15,2.95)(0.3,2.9)
%\qbezier(0,3)(-0.15,2.95)(-0.3,2.9)

%\qbezier(0,3)(0.05,2.85)(0.1,2.7)
%\qbezier(0,3)(-0.05,2.85)(-0.1,2.7)

}

\put(4.2,1.6){$=\quad[2]$}
\put(8.2,0){\linethickness{2.5pt}\qbezier(0,0.5)(0,2)(0,3.5)}

\linethickness{0.6pt}
\put(13.5,0){
\put(1.8,2){\oval(1.6,2.5)}
\qbezier(1,1.5)(1,0.8)(0.5,0.2)
\qbezier(1,2.5)(1,3.2)(0.5,3.8)
\linethickness{2.5pt}
\qbezier(1,1.5)(1,2)(1,2.5)
\linethickness{0.6pt}
\put(2.6,2)
{\qbezier(0,0)(-0.12,0.12)(-0.24,0.24)
\qbezier(0,0)(0.12,0.12)(0.24,0.24)
}
\put(0.9,0.9)
{\qbezier(0,0)(0.03,-0.15)(0.06,-0.3)
\qbezier(0,0)(-0.13,-0.1)(-0.26,-0.2)
}
\put(0.8,3.32)
{\qbezier(0,0)(0,-0.16)(0,-0.32)
\qbezier(0,0)(0.14,-0.1)(0.28,-0.2)
}
\put(4.5,1.8){$=\quad [n-1]$}
\put(10.5,0){\qbezier(0,0.2)(0,2)(0,3.8)
\qbezier(0,2)(-0.12,1.88)(-0.24,1.76)
\qbezier(0,2)(0.12,1.88)(0.24,1.76)
}
}

\end{picture}
\end{center}
}}

{\Large{
\begin{center}
\setlength{\unitlength}{5mm}
\begin{picture}(35,3.5) 
\linethickness{0.6pt}
\put(3,0)
{
\put(2,2){\oval(2,3)}
\qbezier(1,1.5)(1,0.8)(0.5,0.2)
\qbezier(1,2.5)(1,3.2)(0.5,3.8)
\qbezier(3,1.5)(3,0.8)(3.5,0.2)
\qbezier(3,2.5)(3,3.2)(3.5,3.8)
\linethickness{2.5pt}
\qbezier(1,1.5)(1,2)(1,2.5)
\qbezier(3,1.5)(3,2)(3,2.5)
\linethickness{0.6pt}
%\put(2.6,2)
%{\qbezier(0,0)(-0.12,0.12)(-0.24,0.24)
%\qbezier(0,0)(0.12,0.12)(0.24,0.24)
%}
\put(0.9,0.9)
{\qbezier(0,0)(0.03,-0.15)(0.06,-0.3)
\qbezier(0,0)(-0.13,-0.1)(-0.26,-0.2)
}
\put(3.1,3.1)
{\qbezier(0,0)(-0.03,0.15)(-0.06,0.3)
\qbezier(0,0)(0.13,0.1)(0.26,0.2)
}
\put(0.8,3.32)
{\qbezier(0,0)(0,-0.16)(0,-0.32)
\qbezier(0,0)(0.14,-0.1)(0.28,-0.2)
}
\put(3.2,0.68)
{\qbezier(0,0)(0,0.16)(0,0.32)
\qbezier(0,0)(-0.14,0.1)(-0.28,0.2)
}
\put(1.9,0.5)
{\qbezier(0,0)(0.14,0.1)(0.28,0.2)
\qbezier(0,0)(0.14,-0.1)(0.28,-0.2)
}
\put(2.1,3.5){
\qbezier(0,0)(-0.14,0.1)(-0.28,0.2)
\qbezier(0,0)(-0.14,-0.1)(-0.28,-0.2)
}
}
\put(8.5,1.8){$=$}
\put(10,0)
{\qbezier(0.5,0.2)(2,2.5)(3.5,0.2)
\qbezier(0.5,3.8)(2,1.5)(3.5,3.8)
\put(2.1,1.35)
{\qbezier(0,0)(-0.14,0.1)(-0.28,0.2)
\qbezier(0,0)(-0.14,-0.1)(-0.28,-0.2)
}
\put(1.9,2.65)
{\qbezier(0,0)(0.14,0.1)(0.28,0.2)
\qbezier(0,0)(0.14,-0.1)(0.28,-0.2)
}
}
\put(14.5,1.8){$+\,\,[n-2]$}
\put(18.5,0)
{\qbezier(0.5,0.2)(2.2,2)(0.5,3.8)
\qbezier(3.5,0.2)(1.8,2)(3.5,3.8)
\put(1.35,2.1)
{\qbezier(0,0)(-0.1,-0.14)(-0.2,-0.28)
\qbezier(0,0)(0.1,-0.14)(0.2,-0.28)
}
\put(2.65,1.9)
{\qbezier(0,0)(-0.1,0.14)(-0.2,0.28)
\qbezier(0,0)(0.1,0.14)(0.2,0.28)
}
}
\end{picture}
\end{center}
}}
{\Large{
\begin{center}
\setlength{\unitlength}{5mm}
\begin{picture}(22,6) 
\linethickness{0.6pt}
\put(3,0){
\put(0,0){\GRF}
\put(0,4){\GR}
\put(1,2){\GRF}
\qbezier(2,0)(2,1)(2,2)
\qbezier(0,4)(0,3)(0,2)
\qbezier(2,6)(2,5)(2,4)
\put(2,6)
{\qbezier(0,0)(-0.08,-0.08)(-0.16,-0.16)
\qbezier(0,0)(0.08,-0.08)(0.16,-0.16)
}
\put(0,3)
{\qbezier(0,0)(-0.08,-0.08)(-0.16,-0.16)
\qbezier(0,0)(0.08,-0.08)(0.16,-0.16)
}
\put(2,1)
{\qbezier(0,0)(-0.08,-0.08)(-0.16,-0.16)
\qbezier(0,0)(0.08,-0.08)(0.16,-0.16)
}
\put(1,2)
{\qbezier(0.00,2.00)(0,1.90)(0,1.80)
\qbezier(0.00,2.00)(0.1,2)(0.20,2)
}
\put(0.7,-1.7)
{\qbezier(0.00,2.00)(0,1.90)(0,1.80)
\qbezier(0.00,2.00)(0.1,2)(0.20,2)
}
\put(-0.7,-1.7)
{\qbezier(1.00,2.00)(1,1.9)(1,1.80)
\qbezier(1.00,2.00)(0.9,2)(0.8,2)
}
\put(0,0)
{\qbezier(1.00,2.00)(1,1.9)(1,1.80)
\qbezier(1.00,2.00)(0.9,2)(0.8,2)
}
}
\put(6,2.9){$-$}

\put(8,0){
\put(0,2){\GRF}
\qbezier(1,0)(1,1)(1,2)
\qbezier(0,0)(0,1)(0,2)
\qbezier(0,6)(0,5)(0,4)
\qbezier(1,6)(1,5)(1,4)
\qbezier(2,0)(2,3)(2,6)
\put(0,6)
{\qbezier(0,0)(-0.08,-0.08)(-0.16,-0.16)
\qbezier(0,0)(0.08,-0.08)(0.16,-0.16)
}
\put(1,6)
{\qbezier(0,0)(-0.08,-0.08)(-0.16,-0.16)
\qbezier(0,0)(0.08,-0.08)(0.16,-0.16)
}
\put(2,6)
{\qbezier(0,0)(-0.08,-0.08)(-0.16,-0.16)
\qbezier(0,0)(0.08,-0.08)(0.16,-0.16)
}
\put(0,1)
{\qbezier(0,0)(-0.08,-0.08)(-0.16,-0.16)
\qbezier(0,0)(0.08,-0.08)(0.16,-0.16)
}
\put(1,1)
{\qbezier(0,0)(-0.08,-0.08)(-0.16,-0.16)
\qbezier(0,0)(0.08,-0.08)(0.16,-0.16)
}

}
\put(11,2.85){$=$}
\put(13,0){
\put(0,2){\GRF}
\put(1,4){\GR}
\put(1,0){\GRF}
\qbezier(0,0)(0,1)(0,2)
\qbezier(2,4)(2,3)(2,2)
\qbezier(0,6)(0,5)(0,4)
\put(0,6)
{\qbezier(0,0)(-0.08,-0.08)(-0.16,-0.16)
\qbezier(0,0)(0.08,-0.08)(0.16,-0.16)
}
\put(0,1)
{\qbezier(0,0)(-0.08,-0.08)(-0.16,-0.16)
\qbezier(0,0)(0.08,-0.08)(0.16,-0.16)
}
\put(2,3)
{\qbezier(0,0)(-0.08,-0.08)(-0.16,-0.16)
\qbezier(0,0)(0.08,-0.08)(0.16,-0.16)
}
\put(1,0)
{\qbezier(0.00,2.00)(0,1.90)(0,1.80)
\qbezier(0.00,2.00)(0.1,2)(0.20,2)
}
\put(1.7,-1.7)
{\qbezier(0.00,2.00)(0,1.90)(0,1.80)
\qbezier(0.00,2.00)(0.1,2)(0.20,2)
}
\put(0.3,-1.7)
{\qbezier(1.00,2.00)(1,1.9)(1,1.80)
\qbezier(1.00,2.00)(0.9,2)(0.8,2)
}
\put(0,2)
{\qbezier(1.00,2.00)(1,1.9)(1,1.80)
\qbezier(1.00,2.00)(0.9,2)(0.8,2)
}

}

\put(16,2.9){$-$}

\put(18,0){
\put(1,2){\GRF}
\qbezier(1,0)(1,1)(1,2)
\qbezier(2,0)(2,1)(2,2)
\qbezier(2,6)(2,5)(2,4)
\qbezier(1,6)(1,5)(1,4)
\qbezier(0,0)(0,3)(0,6)
\put(0,6)
{\qbezier(0,0)(-0.08,-0.08)(-0.16,-0.16)
\qbezier(0,0)(0.08,-0.08)(0.16,-0.16)
}
\put(1,6)
{\qbezier(0,0)(-0.08,-0.08)(-0.16,-0.16)
\qbezier(0,0)(0.08,-0.08)(0.16,-0.16)
}
\put(2,6)
{\qbezier(0,0)(-0.08,-0.08)(-0.16,-0.16)
\qbezier(0,0)(0.08,-0.08)(0.16,-0.16)
}
\put(1,1)
{\qbezier(0,0)(-0.08,-0.08)(-0.16,-0.16)
\qbezier(0,0)(0.08,-0.08)(0.16,-0.16)
}
\put(2,1)
{\qbezier(0,0)(-0.08,-0.08)(-0.16,-0.16)
\qbezier(0,0)(0.08,-0.08)(0.16,-0.16)
}

}

\end{picture}
\end{center}
}}

Again by the above relations we mean the respective equality of the value of the function $\hat{P}_n$ on the trivalent graphs which look the same except one small part where they are like the left and right side, respectively, of the above equalities. We denote by $[k]$ the quantum integer $k$, given by the expression:
$$[k]=\frac{q^{k}-q^{-k}}{q-q^{-1}}=q^{k-1}+q^{k-3}+\ldots+q^{1-k}.$$
The consistency of the above axioms is shown in \cite{moy}. The last relation can be rewritten by introducing four-valent vertices and a new type of (thick) edges, that we label with the number 3:

\begin{center}
\setlength{\unitlength}{5mm}
\begin{picture}(22,5) 
\linethickness{0.6pt}
\put(9,0)
{\qbezier(0,0)(0.75,0.75)(1.5,1.5)
\qbezier(1.5,0)(1.5,0.75)(1.5,1.5)
\qbezier(3,0)(2.25,0.75)(1.5,1.5)
\qbezier(0,5)(0.75,4.25)(1.5,3.5)
\qbezier(1.5,5)(1.5,4.25)(1.5,3.5)
\qbezier(3,5)(2.25,4.25)(1.5,3.5)
\put(1,2.2){3}
\linethickness{3pt}
\qbezier(1.5,1.5)(1.5,2.5)(1.5,3.5)
\linethickness{0.6pt}
\put(1.5,0.7){
\qbezier(0,0)(-0.08,-0.08)(-0.16,-0.16)
\qbezier(0,0)(0.08,-0.08)(0.16,-0.16)
}
\put(1.5,4.5){
\qbezier(0,0)(-0.08,-0.08)(-0.16,-0.16)
\qbezier(0,0)(0.08,-0.08)(0.16,-0.16)
}
\put(0.7,0.7)
{\qbezier(0,0)(0,-0.1)(0,-0.2)
\qbezier(0,0)(-0.1,0)(-0.2,0)
}
\put(2.5,4.5)
{\qbezier(0,0)(0,-0.1)(0,-0.2)
\qbezier(0,0)(-0.1,0)(-0.2,0)
}
\put(2.3,0.7)
{\qbezier(0,0)(0,-0.1)(0,-0.2)
\qbezier(0,0)(0.1,0)(0.2,0)
}
\put(0.5,4.5)
{\qbezier(0,0)(0,-0.1)(0,-0.2)
\qbezier(0,0)(0.1,0)(0.2,0)
}

}
\end{picture}
\end{center}
and we def\mbox{}ine the value of $\hat{P}_n$ on this graph as equal to the expression on either side of the last of the above axioms. In other words, we can replace the f\mbox{}ifth axiom by the following two relations (see \cite{moy} and Section 1 of \cite{kovroz}):
\begin{equation}
%\begin{center}
\setlength{\unitlength}{5mm}
\begin{picture}(22,6) 
\linethickness{0.6pt}
\put(3,0)
{
\put(0,0){\GRF}
\put(0,4){\GR}
\put(1,2){\GRF}
\qbezier(2,0)(2,1)(2,2)
\qbezier(0,4)(0,3)(0,2)
\qbezier(2,6)(2,5)(2,4)
\put(2,6)
{\qbezier(0,0)(-0.08,-0.08)(-0.16,-0.16)
\qbezier(0,0)(0.08,-0.08)(0.16,-0.16)
}
\put(0,3)
{\qbezier(0,0)(-0.08,-0.08)(-0.16,-0.16)
\qbezier(0,0)(0.08,-0.08)(0.16,-0.16)
}
\put(2,1)
{\qbezier(0,0)(-0.08,-0.08)(-0.16,-0.16)
\qbezier(0,0)(0.08,-0.08)(0.16,-0.16)
}
\put(1,2)
{\qbezier(0.00,2.00)(0,1.90)(0,1.80)
\qbezier(0.00,2.00)(0.1,2)(0.20,2)
}
\put(0.7,-1.7)
{\qbezier(0.00,2.00)(0,1.90)(0,1.80)
\qbezier(0.00,2.00)(0.1,2)(0.20,2)
}
\put(-0.7,-1.7)
{\qbezier(1.00,2.00)(1,1.9)(1,1.80)
\qbezier(1.00,2.00)(0.9,2)(0.8,2)
}
\put(0,0)
{\qbezier(1.00,2.00)(1,1.9)(1,1.80)
\qbezier(1.00,2.00)(0.9,2)(0.8,2)
}
}
\put(6,2.9){$=$}

\put(8,0){
\put(0,2){\GRF}
\qbezier(1,0)(1,1)(1,2)
\qbezier(0,0)(0,1)(0,2)
\qbezier(0,6)(0,5)(0,4)
\qbezier(1,6)(1,5)(1,4)
\qbezier(2,0)(2,3)(2,6)
\put(0,6)
{\qbezier(0,0)(-0.08,-0.08)(-0.16,-0.16)
\qbezier(0,0)(0.08,-0.08)(0.16,-0.16)
}
\put(1,6)
{\qbezier(0,0)(-0.08,-0.08)(-0.16,-0.16)
\qbezier(0,0)(0.08,-0.08)(0.16,-0.16)
}
\put(2,6)
{\qbezier(0,0)(-0.08,-0.08)(-0.16,-0.16)
\qbezier(0,0)(0.08,-0.08)(0.16,-0.16)
}
\put(0,1)
{\qbezier(0,0)(-0.08,-0.08)(-0.16,-0.16)
\qbezier(0,0)(0.08,-0.08)(0.16,-0.16)
}
\put(1,1)
{\qbezier(0,0)(-0.08,-0.08)(-0.16,-0.16)
\qbezier(0,0)(0.08,-0.08)(0.16,-0.16)
}

}
\put(11,2.85){$+$}

\put(13.5,0.5)
{
\qbezier(-0.5,-0.5)(0.5,0.5)(1.5,1.5)
\qbezier(1.5,-0.5)(1.5,0.5)(1.5,1.5)
\qbezier(3.5,-0.5)(2.5,0.5)(1.5,1.5)
\qbezier(-0.5,5.5)(0.5,4.5)(1.5,3.5)
\qbezier(1.5,5.5)(1.5,4.5)(1.5,3.5)
\qbezier(3.5,5.5)(2.5,4.5)(1.5,3.5)
\put(1,2.2){3}
\linethickness{3pt}
\qbezier(1.5,1.5)(1.5,2.5)(1.5,3.5)
\linethickness{0.6pt}
\put(1.5,0.7){
\qbezier(0,0)(-0.08,-0.08)(-0.16,-0.16)
\qbezier(0,0)(0.08,-0.08)(0.16,-0.16)
}
\put(1.5,4.5){
\qbezier(0,0)(-0.08,-0.08)(-0.16,-0.16)
\qbezier(0,0)(0.08,-0.08)(0.16,-0.16)
}
\put(0.7,0.7)
{\qbezier(0,0)(0,-0.1)(0,-0.2)
\qbezier(0,0)(-0.1,0)(-0.2,0)
}
\put(2.5,4.5)
{\qbezier(0,0)(0,-0.1)(0,-0.2)
\qbezier(0,0)(-0.1,0)(-0.2,0)
}
\put(2.3,0.7)
{\qbezier(0,0)(0,-0.1)(0,-0.2)
\qbezier(0,0)(0.1,0)(0.2,0)
}
\put(0.5,4.5)
{\qbezier(0,0)(0,-0.1)(0,-0.2)
\qbezier(0,0)(0.1,0)(0.2,0)
}
}
\end{picture}
%\end{center}
\label{axioma51}
\end{equation}
\begin{equation}
%\begin{center}
\setlength{\unitlength}{5mm}
\begin{picture}(22,6) 
\linethickness{0.6pt}
\put(3,0){
\put(0,2){\GRF}
\put(1,4){\GR}
\put(1,0){\GRF}
\qbezier(0,0)(0,1)(0,2)
\qbezier(2,4)(2,3)(2,2)
\qbezier(0,6)(0,5)(0,4)
\put(0,6)
{\qbezier(0,0)(-0.08,-0.08)(-0.16,-0.16)
\qbezier(0,0)(0.08,-0.08)(0.16,-0.16)
}
\put(0,1)
{\qbezier(0,0)(-0.08,-0.08)(-0.16,-0.16)
\qbezier(0,0)(0.08,-0.08)(0.16,-0.16)
}
\put(2,3)
{\qbezier(0,0)(-0.08,-0.08)(-0.16,-0.16)
\qbezier(0,0)(0.08,-0.08)(0.16,-0.16)
}
\put(1,0)
{\qbezier(0.00,2.00)(0,1.90)(0,1.80)
\qbezier(0.00,2.00)(0.1,2)(0.20,2)
}
\put(1.7,-1.7)
{\qbezier(0.00,2.00)(0,1.90)(0,1.80)
\qbezier(0.00,2.00)(0.1,2)(0.20,2)
}
\put(0.3,-1.7)
{\qbezier(1.00,2.00)(1,1.9)(1,1.80)
\qbezier(1.00,2.00)(0.9,2)(0.8,2)
}
\put(0,2)
{\qbezier(1.00,2.00)(1,1.9)(1,1.80)
\qbezier(1.00,2.00)(0.9,2)(0.8,2)
}

}

\put(6,2.9){$=$}

\put(8,0){
\put(1,2){\GRF}
\qbezier(1,0)(1,1)(1,2)
\qbezier(2,0)(2,1)(2,2)
\qbezier(2,6)(2,5)(2,4)
\qbezier(1,6)(1,5)(1,4)
\qbezier(0,0)(0,3)(0,6)
\put(0,6)
{\qbezier(0,0)(-0.08,-0.08)(-0.16,-0.16)
\qbezier(0,0)(0.08,-0.08)(0.16,-0.16)
}
\put(1,6)
{\qbezier(0,0)(-0.08,-0.08)(-0.16,-0.16)
\qbezier(0,0)(0.08,-0.08)(0.16,-0.16)
}
\put(2,6)
{\qbezier(0,0)(-0.08,-0.08)(-0.16,-0.16)
\qbezier(0,0)(0.08,-0.08)(0.16,-0.16)
}
\put(1,1)
{\qbezier(0,0)(-0.08,-0.08)(-0.16,-0.16)
\qbezier(0,0)(0.08,-0.08)(0.16,-0.16)
}
\put(2,1)
{\qbezier(0,0)(-0.08,-0.08)(-0.16,-0.16)
\qbezier(0,0)(0.08,-0.08)(0.16,-0.16)
}

}
\put(11,2.85){$+$}
\put(13.5,0.5)
{
\qbezier(-0.5,-0.5)(0.5,0.5)(1.5,1.5)
\qbezier(1.5,-0.5)(1.5,0.5)(1.5,1.5)
\qbezier(3.5,-0.5)(2.5,0.5)(1.5,1.5)
\qbezier(-0.5,5.5)(0.5,4.5)(1.5,3.5)
\qbezier(1.5,5.5)(1.5,4.5)(1.5,3.5)
\qbezier(3.5,5.5)(2.5,4.5)(1.5,3.5)
\put(1,2.2){3}
\linethickness{3pt}
\qbezier(1.5,1.5)(1.5,2.5)(1.5,3.5)
\linethickness{0.6pt}
\put(1.5,0.7){
\qbezier(0,0)(-0.08,-0.08)(-0.16,-0.16)
\qbezier(0,0)(0.08,-0.08)(0.16,-0.16)
}
\put(1.5,4.5){
\qbezier(0,0)(-0.08,-0.08)(-0.16,-0.16)
\qbezier(0,0)(0.08,-0.08)(0.16,-0.16)
}
\put(0.7,0.7)
{\qbezier(0,0)(0,-0.1)(0,-0.2)
\qbezier(0,0)(-0.1,0)(-0.2,0)
}
\put(2.5,4.5)
{\qbezier(0,0)(0,-0.1)(0,-0.2)
\qbezier(0,0)(-0.1,0)(-0.2,0)
}
\put(2.3,0.7)
{\qbezier(0,0)(0,-0.1)(0,-0.2)
\qbezier(0,0)(0.1,0)(0.2,0)
}
\put(0.5,4.5)
{\qbezier(0,0)(0,-0.1)(0,-0.2)
\qbezier(0,0)(0.1,0)(0.2,0)
}
}
\end{picture}
%\end{center}
\label{axioma52}
\end{equation}

Also, from \cite{moy} we know that the value of $\hat{P}_n(\Gamma)$ on every trivalent graph $\Gamma$ with thick edges is a Laurent polynomial in $q$ and $q^{-1}$ whose coef\mbox{}f\mbox{}icients are nonnegative integers.\\ 
\indent Like in the case of the Jones polynomial, we will extend the def\mbox{}inition of the function $\hat{P}_n(\Gamma)$ to oriented links by the following recursive relations:
\begin{eqnarray}
\hat{P}_n(D_+)&=&q^{1-n}\hat{P}_n(D_0)-q^{-n}\hat{P}_n(D_{\Gamma}) \label{slnr1}\\
\hat{P}_n(D_-)&=&q^{n-1}\hat{P}_n(D_0)-q^{n}\hat{P}_n(D_{\Gamma}). \label{slnr2}
\end{eqnarray} 
Here we denote by $D_+$, $D_-$, $D_0$  and $D_{\Gamma}$  any four diagrams that look the same except near one crossing where they look like:

\begin{center}
\setlength{\unitlength}{5mm}
\begin{picture}(25,3)
\linethickness{0.6pt}
\put(5,1){\XS}
\put(10,1){\YS}
\put(15,1){\OS}
\put(20,1){\GR}
\put(5.2,-0.2){$D_{+}$}
\put(10.2,-0.2){$D_{-}$}
\put(15.2,-0.2){$D_{0}$}
\put(20.2,-0.2){$D_{\Gamma}$}
\end{picture}
\end{center}

Then it can be easily verif\mbox{}ied that $\hat{P}_n(D)$ is invariant under the Reidemeister moves, and hence that it is a link invariant. Thus, we can write $\hat{P}_n(L)$.  Furthermore, from (\ref{slnr1}) and (\ref{slnr2}) we obtain the following skein relation for $\hat{P}_n$:
$$q^n\hat{P}_n(L_+)-q^{-n}\hat{P}_n(L_-)=(q-q^{-1})\hat{P}_n(L_0).$$
Hence, if we def\mbox{}ine $P_n(L)=\hat{P}_n(L)/[n]$, we obtain the $n$-specialization of the HOMFLYPT polynomial.
  
\section{Khovanov-Rozansky homology ($sl(n)$ link homology)}\label{kovroznot}

\subsection{The construction of the chain complex}

The categorif\mbox{}ication of the $sl(n)$ link invariant, analogous to the categorif\mbox{}ica\-tion of the Jones polynomial, was def\mbox{}ined by Khovanov and Rozansky in \cite{kovroz}. The construction follows the same lines as the cubic complex construction of Section \ref{kovnot}, and uses the state-sum model for the $sl(n)$ link invariant from the previous section.
The main dif\mbox{}ference is that the contributions from the states (i.e. the total resolutions) given by $\hat{P}_n(\Gamma)$ are much more complicated expressions (so that they satisfy the f\mbox{}ive axioms from the previous section) than in the Jones polynomial case (where the contributions are just $(q+q^{-1})^k$). Fortunately, because of the positivity of the calculus (the coef\mbox{}f\mbox{}icients of $\hat{P}_n(\Gamma)$ are nonnegative integers), it was possible to assign graded vector space to any state $\Gamma$, with graded dimension equal to $\hat{P}_n(\Gamma)$ (see below).\\
\indent Again for every regular diagram $D$, they also def\mbox{}ined a graded chain complex whose graded Euler characteristic is equal to the unnormalized $sl(n)$ link invariant - $\hat{P}_n(D)$. However, we will use a slightly dif\mbox{}ferent, intermediate, normalization $Q_n(D)$ instead of $\hat{P}_n(D)$ from the previous section. On the level of trivalent graphs $Q_n$ and $\hat{P}_n$ coincide. However, for the recursive relations we take:
\begin{eqnarray}
Q_n(D_+)&=&Q_n(D_0)-q^{-1}Q_n(D_{\Gamma}) \label{qslnr1}\\
Q_n(D_-)&=&Q_n(D_{\Gamma})-q^{-1}Q_n(D_0). \label{qslnr2}
\end{eqnarray} 
As a result, we obviously have $\hat{P}_n(D)=(-1)^{-n_-}q^{-(n-1)\cdot n_+ + n \cdot n_-}Q_n(D)$.\\
\indent To every total resolution $D_{\epsilon}$ they assigned a graded  $\Q$-vector space $\bar{V}_{\epsilon}$ which is the cohomology of a certain 2-periodic complex $C(D_{\epsilon})$ (for details see \cite{kovroz} and the following subsection). It can be shown that $q\dim \bar{V}_{\epsilon}=\hat{P}_{n}(D_{\epsilon})=Q_{n}(D_{\epsilon})$. Indeed, this is proved by showing that the complex $C(D_{\epsilon})$ categorif\mbox{}ies the axioms for $\hat{P}_{n}(D_{\epsilon})$, in the sense that, for example, the second axiom becomes:
$$C(\Gamma_1) \cong C(\Gamma_2)\{1\} \oplus C(\Gamma_2)\{-1\},$$ 
where by $\Gamma_1$ and $\Gamma_2$ we have denoted the graphs on the left and on the righthand side, respectively, of the second axiom.\\
\indent Again, the resolutions are organized in the same cubic complex as in Section \ref{kovnot} (just bearing in mind the dif\mbox{}ferent 0- and 1-resolutions, see (\ref{rezsln})) and we obtain the chain groups by summing along columns:  \begin{equation}
C_n^i(D)=\bigoplus_{|\epsilon|=i}{\bar{V}_{\epsilon}}\{-i\}.
\label{form1}
\end{equation}
\indent Furthermore, we def\mbox{}ine the grading-preserving dif\mbox{}ferential $d_n^{i}:C_{n}^i(D)\to C_{n}^{i+1}(D)$ as a signed sum of per-edge maps (the signs are the same as in the $sl(2)$-case, Section \ref{kovnot}).  For the explicit def\mbox{}inition of the per-edge maps see \cite{kovroz} and the following subsection. For our purposes, we want to point out that the complex associated to a (generalized) diagram $D$, with a specif\mbox{}ied crossing $c$, is again the mapping cone of a certain homomorphism $f:C_n(D_0) \to C_n(D_1)\{-1\}$, where by $D_i$, $i=0,1$, we have denoted the diagram obtained from $D$ by performing the $i$-resolution of the crossing $c$ (according to (\ref{rezsln})). Furthermore, a complex $C_n(D')$ can be associated to any generalized diagram $D'$, i.e. a regular diagram with thick edges.\\
\indent We def\mbox{}ine the homology groups of the complex $(C_n(D),d_n)$ by $H_n^i(D)=H^i(C_n(D))$.
Like in the $sl(2)$ case, the homology groups $H_n^i(D)$ inherit a natural grading from the internal gradings of the chain groups and we write $H_n^{i,j}(D)$ for the subset of elements of $H_n^i(D)$ of the degree $j$.\\
\indent  Finally we def\mbox{}ine a graded chain complex $\C_n(D)$ by 
$$\C_n(D)=C_n(D)[-n_-]\{(1-n)n_+ +n\cdot n_-\},$$
where $n_-$ and $n_+$ are the numbers of negative and positive crossings of $D$, respectively. \\
\indent Since the dif\mbox{}ferential is grading preserving we have that the graded Euler characteristic of the chain complex $\C_n(D)$ is equal to $\hat{P}_n(D)$. In fact, we have 
\begin{theorem}{\cite{kovroz}}
The homology groups $\H^i_n(D)=H^i(\C_n(D))$ are independent of the planar projection $D$ of the link $L$. Hence, the homology groups are link invariants and we can write $\H^i_n(L)$.
\label{krgl}
\end{theorem} 
\begin{remark}
If we put  $n=2$ in the $sl(n)$ theory \cite{kovroz} we obtain the standard Khovanov homology \cite{kov} with the $q$-gradings inverted, i.e. we have an isomorphism $\H_2^{i,j}(L)\cong\H^{i,-j}(L)\otimes_{\Z}\Q$.  We would obtain the same convention for the standard Khovanov homology (Section \ref{kovnot}) if we def\mbox{}ined $\deg 1=-1$, $\deg X=1$ and in all later shifts in $q$-gradings we replaced $\{i\}$  by $\{-i\}$. These conventions for the $sl(2)$ homology theory are used in \cite{tan} and \cite{pat}.
\end{remark}

\subsection{Some details concerning the construction} \label{kovroznot2}

As we mentioned above, the explicit construction of the complexes $C(\Gamma)$ and $C_n(D)$ from the previous subsection is rather complicated, and hence we omit the technical details (for these, see \cite{kovroz}). However, we want to point some features of the construction and extract some elementary properties that we will use in the thesis.\\
\indent In the construction, one starts with a regular diagram $D$ of a link $L$, and assigns a dif\mbox{}ferent variable $x_i$ to every arc between two crossings. We def\mbox{}ine the ring $R$ as the ring of polynomials over $\Q$ in all variables $x_i$. We introduce a grading on $R$ by def\mbox{}ining $\deg x_i=2$. Then every resolution naturally acquires a labelling, according to the following picture:
\bigskip 
\begin{center}
\setlength{\unitlength}{5mm}
\begin{picture}(25,5) 
\linethickness{0.6pt}
%\line(1.00,2.00,1.00,3.00)
\put(-0.5,0.5){
\qbezier(4.00,1.00)(4.00,3.00)(4.00,5.00)
\qbezier(3.80,4.60)(3.90,4.80)(4.00,5.00)
\qbezier(4.20,4.60)(4.10,4.80)(4.00,5.00) \qbezier(2.00,1.00)(2.00,3.00)(2.00,5.00)
\qbezier(1.80,4.60)(1.90,4.80)(2.00,5.00)
\qbezier(2.20,4.60)(2.10,4.80)(2.00,5.00) 
\qbezier(9.5,3)(10.5,3)(11.5,3)
\qbezier(11.5,3)(11.3,3.1)(11.1,3.2)
\qbezier(11.5,3)(11.3,2.9)(11.1,2.8)

\qbezier(4.5,3)(5.5,3)(6.5,3)
\qbezier(4.5,3)(4.7,3.1)(4.9,3.2)
\qbezier(4.5,3)(4.7,2.9)(4.9,2.8)

\put(10,0){
\qbezier(9,1)(8,3)(7,5)
\qbezier(7,1)(7.4,1.8)(7.8,2.6)
\qbezier(9,5)(8.6,4.2)(8.2,3.4)
}

\linethickness{2.5pt}
\qbezier(13.00,2.00)(13.00,3.00)(13.00,4.00)
\linethickness{0.6pt}
\qbezier(12.00,1.00)(12.50,1.50)(13.00,2.00)
\qbezier(14.00,1.00)(13.50,1.50)(13.00,2.00)
\qbezier(12.00,5.00)(12.50,4.50)(13.00,4.00)
\qbezier(14.00,5.00)(13.50,4.50)(13.00,4.00)

\qbezier(14.00,5.00)(13.80,4.90)(13.60,4.80)
\qbezier(14.00,5.00)(13.90,4.80)(13.80,4.60)

\qbezier(12.00,5.00)(12.20,4.90)(12.40,4.80)
\qbezier(12.00,5.00)(12.10,4.80)(12.20,4.60)

\qbezier(12.50,1.50)(12.30,1.40)(12.10,1.30)\qbezier(12.50,1.50)(12.40,1.30)(12.30,1.10)

\qbezier(13.50,1.50)(13.70,1.40)(13.90,1.30)\qbezier(13.50,1.50)(13.60,1.30)(13.70,1.10)
\qbezier(14.5,3)(15.5,3)(16.5,3)
\qbezier(14.5,3)(14.7,3.1)(14.9,3.2)
\qbezier(14.5,3)(14.7,2.9)(14.9,2.8)

\qbezier(19.5,3)(20.5,3)(21.5,3)
\qbezier(21.5,3)(21.3,3.1)(21.1,3.2)
\qbezier(21.5,3)(21.3,2.9)(21.1,2.8)
\put(-10,0){
\qbezier(19,5)(18,3)(17,1)
\qbezier(17,5)(17.4,4.2)(17.8,3.4)
\qbezier(19,1)(18.6,1.8)(18.2,2.6)

\qbezier(19,5)(19,4.75)(19,4.5)
\qbezier(19,5)(18.8,4.8)(18.6,4.6)
\qbezier(17,5)(17,4.75)(17,4.5)
\qbezier(17,5)(17.2,4.8)(17.4,4.6)
}
\put(10,0){
\qbezier(9,5)(9,4.75)(9,4.5)
\qbezier(9,5)(8.8,4.8)(8.6,4.6)
\qbezier(7,5)(7,4.75)(7,4.5)
\qbezier(7,5)(7.2,4.8)(7.4,4.6)
}

\qbezier(22.00,1.00)(22.00,3.00)(22.00,5.00)
\qbezier(21.80,4.60)(21.90,4.80)(22.00,5.00)
\qbezier(22.20,4.60)(22.10,4.80)(22.00,5.00)

\qbezier(24.00,1.00)(24.00,3.00)(24.00,5.00)
\qbezier(23.80,4.60)(23.90,4.80)(24.00,5.00)
\qbezier(24.20,4.60)(24.10,4.80)(24.00,5.00)
\small{
\put(4.70,3.50){0-res.}
\put(9.6,3.50){1-res.}
\put(14.7,3.50){0-res.}
\put(19.6,3.50){1-res.}
}
\footnotesize{
\put(1.2,5.2){$x_1$}
\put(6.2,5.2){$x_1$}
\put(11.2,5.2){$x_1$}
\put(16.2,5.2){$x_1$}
\put(21.2,5.2){$x_1$}

\put(1.2,0.6){$x_4$}
\put(6.2,0.6){$x_4$}
\put(11.2,0.6){$x_4$}
\put(16.2,0.6){$x_4$}
\put(21.2,0.6){$x_4$}

\put(4.2,5.2){$x_2$}
\put(9.2,5.2){$x_2$}
\put(14.2,5.2){$x_2$}
\put(19.2,5.2){$x_2$}
\put(24.2,5.2){$x_2$}

\put(4.2,0.6){$x_3$}
\put(9.2,0.6){$x_3$}
\put(14.2,0.6){$x_3$}
\put(19.2,0.6){$x_3$}
\put(24.2,0.6){$x_3$}

}
}
\put(2.3,-0.3){$\Gamma_0$}
\put(22.3,-0.3){$\Gamma_0$}
\put(12.3,-0.3){$\Gamma_1$}

\end{picture}
\end{center}
\bigskip
\bigskip
Now, to every such resolution one assigns a certain matrix factorization (i.e. a generalization of a 2-periodic complex -- see \cite{kovroz}), and to each total resolution $D_{\epsilon}$ one assigns the tensor product over the crossings of $D$ of the corresponding matrix factorizations. In this way, one obtains a (genuine) 2-periodic complex $C(D_{\epsilon})$, whose homology $H(D_{\epsilon})$ is the vector space $\bar{V}_{\epsilon}$ that we assigned to the resolution $D_{\epsilon}$ in the previous subsection. Thus, we are left with def\mbox{}ining the particular matrix factorizations. \\
\indent In the case when we have a thin edge with two variables: 
 
\begin{center}
\setlength{\unitlength}{5mm}
\begin{picture}(25,2) 
\linethickness{0.6pt}
%\put(7,0.9){{\large $L_j^i:$}}
\put(0,0){\qbezier(8,0.8)(13,0.8)(18,0.8)
\qbezier(18,0.8)(17.85,0.95)(17.7,1.1)
\qbezier(18,0.8)(17.85,0.65)(17.7,0.5)
\put(9,1.3){$x_j$}
\put(16,1.3){$x_i$}}
\end{picture}
\end{center}
we assign to it the matrix factorization $L_j^i$, def\mbox{}ined by:
\begin{equation*}
\begin{CD}
R @>\pi_{i,j}^n>> R\{1-n\} @>x_i-x_j>> R
\end{CD}
\end{equation*}
where $\pi_{i,j}^n$ is given by:
$$\pi_{i,j}^n=\frac{x_i^{n+1}-x_j^{n+1}}{x_i-x_j}=\sum_{k=0}^n{x_i^k x_j^{n-k}}.$$
Finally, to the resolution $\Gamma_0$, we assign the tensor product of matrix factorizations $L_4^1$ and $L_3^2$, corresponding to the thin edges that form it.\\ 
\indent In the case of the resolution with a thick edge, $\Gamma_1$, the corresponding matrix factorization is def\mbox{}ined as the tensor product of the factorizations:
\begin{equation*}
\begin{CD}
R @>u_1>> R\{1-n\} @>x_1+x_2-x_3-x_4>> R
\end{CD}
\end{equation*}
and
\begin{equation*}
\begin{CD}
R @>u_2>> R\{3-n\} @>x_1x_2-x_3x_4>> R
\end{CD}
\end{equation*}
with the grading shifted down by 1.
Here by $u_1$ and $u_2$ we have denoted the following polynomials
\begin{eqnarray*}
u_1&=&u_1(x_1,x_2,x_3,x_4)=\frac{g(x_1+x_2,x_1x_2)-g(x_3+x_4,x_1x_2)}{x_1+x_2-x_3-x_4},\\
u_2&=&u_2(x_1,x_2,x_3,x_4)=\frac{g(x_3+x_4,x_1x_2)-g(x_3+x_4,x_3x_4)}{x_1x_2-x_3x_4},\\
\end{eqnarray*}
where $g(x,y)$ is the unique polynomial such that $g(x+y,xy)=x^{n+1}+y^{n+1}$. Although $u_1$ and $u_2$ have in general a rather complicated explicit expressions, in the case when $x_1+x_2=x_3+x_4$ and $x_1x_2=x_3x_4$, we have that $u_1=\pi_{1,2}^n$ and $u_2=\pi_{1,2}^{n-1}$.\\
\indent The main problem with the def\mbox{}initions from above, is the fact that the homology of the complex that corresponds to the total resolution is incalculable in almost all cases. However, we know that the homology of the total complex is an $R/I$-module, where $I$ is the ideal generated by all polynomials that appear in the def\mbox{}initions of the matrix factorizations (see Proposition 2 of \cite{kovroz} and also Section 2.3 of \cite{bojan}). Furthermore, in some very simple cases like all the total resolutions $D_{\epsilon}$ of a standard diagram $D$ of a torus $(2,n)$-knot (see Chapters \ref{braid} and \ref{torus}), the homology $H(D_{\epsilon})$ is isomorphic to the above quotient ring of polynomials (up to a certain shift in grading). \\

\indent Now, let us move on to the dif\mbox{}ferential. As can be shown, the homology $H(D_{\epsilon})$ is concentrated in only one homological degree of the 2-periodic complex $C(D_{\epsilon})$ (the one with  parity equal to the parity of the number of circles in the total resolution of $D$ when we resolve all the crossing into the resolution $\Gamma_0$). Since the dif\mbox{}ferential should be a map between two homology groups $H(D_{\epsilon})$ and $H(D_{\epsilon'})$ ($\epsilon$ and $\epsilon'$ are related like in the $sl(2)$ case), we will def\mbox{}ine it as the map induced on homology by a homomorphism $f_{\nu}:C(D_{\epsilon}) \to C(D_{\epsilon'})$, def\mbox{}ined as follows:\\
 \indent Obviously, the only nontrivial dif\mbox{}ferentials are between two complexes associated to diagrams $D_{\epsilon}$ and $D_{\epsilon'}$ which look the same except near one crossing where one is of the form $\Gamma_0$ and the another one is of the form $\Gamma_1$. Since the complex $C(D_{\epsilon})$ is def\mbox{}ined as the tensor product of the matrix factorizations assigned to $\Gamma_0$ and $\Gamma_1$, we def\mbox{}ine the homomorphism $f_{\nu}$ as the tensor product of the identity maps over all crossings where we have the same resolutions, and at the crossing where the resolutions dif\mbox{}fer, we def\mbox{}ine it as a homomorphism $\chi_0:C(\Gamma_0)\to C(\Gamma_1)$ or $\chi_1:C(\Gamma_1) \to  C(\Gamma_0)$. Since the chain groups of the matrix factorizations $C(\Gamma_i)$, $i=0,1$, are of rank 2, both maps $\chi_i$ can be described by a pair of $2\times 2$ matrices. Explicitly, $\chi_0$ is given by a pair $(U_0, U_1)$ def\mbox{}ined by
\begin{eqnarray*}
U_0&=&\left(
\begin{array}{cc}
x_1-x_3 & 0 \\
a_1 & 1 
\end{array}
\right),\\
U_1&=&\left(
\begin{array}{cc}
x_1 & -x_3 \\
-1 & 1
\end{array}
\right),
\end{eqnarray*}
where 
$$a_1=\frac{u_1+x_1u_2-\pi_{2,3}^n}{x_1-x_4}.$$
The map $\chi_1$ is given by a pair of matrices $(V_0,V_1)$ def\mbox{}ined by
\begin{eqnarray*}
V_0&=&\left(
\begin{array}{cc}
1 & 0 \\
-a_1 & x_1-x_3 
\end{array}
\right),\\
V_1&=&\left(
\begin{array}{cc}
1 & x_3 \\
1 & x_1
\end{array}
\right).
\end{eqnarray*}
It can easily be seen that both maps $\chi_0$ and $\chi_1$ have degree 1. However, the dif\mbox{}ferentials they induce are grading preserving because of the grading shift in the def\mbox{}inition of the chain groups (\ref{form1}).

\chapter{Categorif\mbox{}ication of the dichromatic polynomial for graphs}
\label{graf}

\section{Introduction}
\indent In this chapter we categorify
the dichromatic polynomial of graphs, and consequently the Tutte polynomial. Since
the dichromatic polynomial, $P_G(q,v)$, is a two-variable polynomial and the standard techniques
of categorif\mbox{}ication work with one-variable polynomials, we will def\mbox{}ine a
class of one-variable specializations $P_n(G)$, one for every positive integer $n$ -- 
a family which is enough to recover the original dichromatic
polynomial -- and
then for each of them we construct a chain complex whose graded Euler characteristic
is equal to $P_n(G)$. \\
\indent As is well-known, to each alternating link $L$ there
corresponds bijectively a planar graph
 $G(L)$, 
and the value of the Jones polynomial of the link $L$ corresponds to
the specialization of the
Tutte polynomial $T(x,1/x)$ or analogously to the specialization of the dichromatic polynomial
$J(q)=P(q,q^2/(q-1))$ (see Appendix \ref{a1} and  \cite{bol} or \cite{kauf}). However, the set of specializations from the previous paragraph ``misses" the polynomial $J(q)$. \\
\indent In order to resolve this, we def\mbox{}ine a dif\mbox{}ferent set of one-variable specializations $Q_n(G)$, one for every integer $n \le 2$, and then categorify each of these one-variable polynomials. Although the chain complexes that categorify this set of specializations have inf\mbox{}inite-dimensional chain groups, for $n=2$ we obtain the categorif\mbox{}ication of the analog $J(q)$ of the Jones polynomial.
Furthermore, in the last case ($n=2$) we also give an alternative description in terms of the enhanced states. \\

\section{Preliminaries} \label{l2}

Let  $G$ be a graph specif\mbox{}ied by a set of vertices $V(G)$ and a set of
edges $E(G)$. If $e\in E(G)$ is an arbitrary edge of the graph $G$,
then by $G-e$ we denote the graph $G$ with the edge $e$ deleted, and
by $G/e$ the graph obtained by contracting edge $e$ (i.e. by
identifying the vertices incident to $e$ and deleting $e$). The
dichromatic polynomial of the graph, $P_G(q,v)$, is a two-variable
generalization of the
chromatic polynomial of the graph. It is given by the following axioms:
\begin{eqnarray*}
(A1)&\quad P_G=P_{G-e}-qP_{G/e},\\
(A2)&\quad P_{N_k}=v^k,
\end{eqnarray*}
where $N_k$ is the graph with $k$ vertices and no edges.\\
\indent These two axioms determine the polynomial uniquely. Obviously, if we put $q=1$, we obtain the usual chromatic polynomial. Furthermore, from $(A1)$ we have a recursive expression for the dichromatic polynomial in terms of the value of the polynomial on graphs with a smaller number of edges. By repeated use of $(A1)$ we will obtain the value of the dichromatic polynomial as a sum of  contributions from all spanning subgraphs of $G$ (subgraphs that contain all vertices of $G$), which we will call states.
Furthermore, if for each subset $s\subset E(G)$ we denote by $[G:s]$ the graph whose set of vertices is $V(G)$ and set of edges is $s$, then the contribution of the graph $[G:s]$ is $(-1)^{|s|}q^{|s|}v^{k(s)}$, where $|s|$ is the number of elements of $s$ and $k(s)$ is the number of connected components of $[G:s]$. Hence, we obtain the expression:
$$P_G(q,v)=\sum_{s\subset E(G)}{(-1)^{|s|}q^{|s|}v^{k(s)}}=\sum_{i\ge 0}{(-1)^i q^i\sum_{s\subset E(G),|s|=i}{v^{k(s)}}},$$
which is called the state-sum expansion of the polynomial $P_G$.\\
\indent The important two-variable Tutte polynomial of a graph $G$, $T_G(x,y)$ is given by a similar recursive relation to the one for the dichromatic polynomial. Namely, it is given by 
$$T_G(x,y)=\left\{ 
\begin{array}{lr}
xT_{G/e}(x,y)&\textrm{  if }e \textrm{ is a bridge}\\
yT_{G-e}(x,y)&\textrm{  if }e \textrm{ is a loop}\\
T_{G-e}(x,y)+T_{G/e}(x,y)&\textrm{  otherwise}\end{array}\right.$$

\indent For the Tutte polynomial there exists a state-sum expression given by:
$$T_G(x,y)=\sum_{s\subset E(G)}{(x-1)^{k(s)-k(E(G))}(y-1)^{|s|-N+k(s)}},$$
where by $N$ we denote the number of vertices of $G$. It is obviously related to the dichromatic polynomial by:
{\small
\begin{eqnarray}
T_G(x,y)&=&\sum_{s\subset E(G)}{(x-1)^{k(s)-k(E(G))}(y-1)^{|s|-N+k(s)}}= \nonumber\\
        &=&(x-1)^{-k(E(G))}(y-1)^{-N}\sum_{s\subset E(G)}{(y-1)^{|s|}{((x-1)(y-1))}^{k(s)}}= \nonumber\\
        &=&(x-1)^{-k(E(G))}(y-1)^{-N}P_G(1-y,(x-1)(y-1)). \label{tutdich}
\end{eqnarray}
}
Hence, we have that the Tutte polynomial $T_G(x,y)$ is a multiple of the dichromatic polynomial $P_G(q,v)$, when we take $q=1-y$ and $v=(x-1)(y-1)$.\\ 

\indent Our aim is to def\mbox{}ine a graded chain complex whose graded Euler characteristic is equal to the dichromatic polynomial. However, since the dichromatic polynomial is a two-variable polynomial, we will f\mbox{}irst def\mbox{}ine an inf\mbox{}inite set of one-variable specializations and then ``categorify" each of the specializations. This situation is very similar to the two-variable HOMFLYPT polynomial for knots and its inf\mbox{}inite set of one-variable specializations (one for each positive integer $n$, with $n=2$ being the Jones polynomial). Note that in this way we will also categorify the Tutte polynomial.\\

\indent In Section \ref{l3}, we will observe the f\mbox{}irst set of specializations of the dichromatic polynomial. Namely, for each positive integer $n$ we def\mbox{}ine:
$$P_{G,n}(q)=P_G(q^n,1+q+\ldots+q^n).$$
\indent In other words, for each positive integer $n$, we replaced in the axioms $(A1)$ and $(A2)$ $q$ by $q^n$, and $v$ by the expression $1+q+\ldots+q^n$, which we will denote by $\kvn$. Hence we have the following state-sum expression for the specialization $P_{G,n}$:
$$P_{G,n}(q)=\sum_{i\ge 0}{(-1)^i \sum_{s\subset E(G),|s|=i}{q^{n|s|}\kvn^{k(s)}}}$$

\begin{figure}[h]
\centerline{\relabelbox
\epsfysize 10cm
\epsfbox{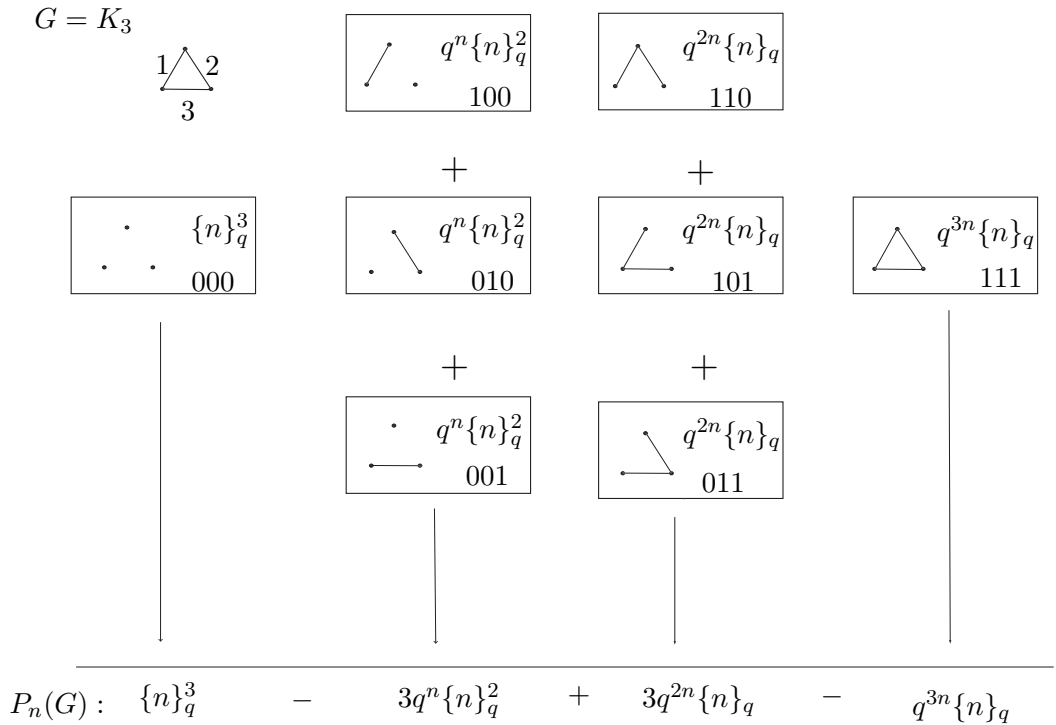}
\adjustrelabel <-15pt,20pt> {G=K3}{$G=K_3$}
\adjustrelabel <-4pt,0pt> {1}{$1$}
\relabel {2}{$2$}
\adjustrelabel <-3pt,-6pt> {3}{$3$}
\relabel {000}{$000$}
\relabel {001}{$001$}
\relabel {010}{$010$}
\relabel {100}{$100$}
\relabel {011}{$011$}
\relabel {101}{$101$}
\relabel {110}{$110$}
\relabel {111}{$111$}
\adjustrelabel <-5pt,-7pt> {n3}{$\kvn^3$}
\adjustrelabel <-13pt,-7pt> {qn1}{$q^n\kvn^2$}
\adjustrelabel <-13pt,-7pt> {qn2}{$q^n\kvn^2$}
\adjustrelabel <-13pt,-7pt> {qn3}{$q^n\kvn^2$}
\adjustrelabel <-13pt,-7pt> {q1n}{$q^{2n}\kvn$}
\adjustrelabel <-14pt,-7pt> {q2n}{$q^{2n}\kvn$}
\adjustrelabel <-8pt,-7pt> {q3n}{$q^{2n}\kvn$}
\adjustrelabel <-16pt,-7pt> {q3n1}{$q^{3n}\kvn$}
\relabel {+1}{{\Large $+$}}
\relabel {+2}{{\Large $+$}}
\relabel {+3}{{\Large $+$}}
\relabel {+4}{{\Large $+$}}
\relabel {pnG:}{$P_n(G):$}
\relabel {v1}{$\kvn^3$}
\adjustrelabel <-12pt,0pt> {v2}{$3q^n\kvn^2$}
\adjustrelabel <-14pt,0pt> {v3}{$3q^{2n}\kvn$}
\adjustrelabel <-15pt,0pt> {v4}{$q^{3n}\kvn$}
\relabel {-}{$-$}
\relabel {+}{$+$}
\relabel {-1}{$-$}
\endrelabelbox}
\caption{State-sum expression for the $n$-specialization of the dichromatic polynomial}
\label{sl1}
\end{figure}

\indent For every positive integer $n$ and graph $G$, we will def\mbox{}ine a graded chain complex whose graded Euler characteristic is equal to $P_{G,n}$. Denote by $m$ the number of edges of the graph $G$. The construction will depend on the ordering of the edges of $G$ ($e_1,\ldots,e_m$). Then each subset $s_{\epsilon}\subset E(G)$ will be uniquely determined by an element $\epsilon=(\epsilon_1,\ldots,\epsilon_m)\in \{0,1\}^m$, 
where $\epsilon_i=1$ if $e_i \in s$ and $\epsilon_i=0$ if $e_i \notin s$. Obviously, each such $m$-tuple is uniquely determined by the subset $s$ of $E(G)$, and so we have a bijective correspondence between $\{0,1\}^m$ and the set of all $s\subset E(G)$. Thus, every element $\epsilon \in \{0,1\}^m$ determines a set $s_{\epsilon}$ of edges of $G$, and hence a graph $[G:s_{\epsilon}]$ which we will denote by $G_{\epsilon}$.  \\
\indent We will construct our complex by ``f\mbox{}lattening" the cube of states (see f\mbox{}igure \ref{sl1}).

To each vertex $\epsilon$ of the cube $\{0,1\}^m$ (which is skewed 
such that in the $i$-th column, $0\le i \le m$, are all the $\epsilon$'s 
such that $|\epsilon|=i$), we will assign a graded $\Z$-module 
whose graded dimension is equal to the contribution of 
the subgraph $G_{\epsilon}$. In other words, a module whose 
graded dimension is $q^{n|\epsilon|}\kvn^{k_{\epsilon}}$, where
$k_{\epsilon}$ is the number of connected components of $G_{\epsilon}$.
We will build chain groups by taking direct sums along the columns.
Finally, we will def\mbox{}ine dif\mbox{}ferentials such that they are grading preserving,
and in this way we will obtain a chain complex whose graded Euler characteristic 
is equal to the $n$-th specialization $P_{G,n}$ of the dichromatic polynomial.\\
\indent In the next section we will def\mbox{}ine the chain complex more precisely.   \\

\indent In order to def\mbox{}ine the second set of specializations of the dichromatic polynomial, we introduce a new variable $a$ as $a=v(q-1)$, i.e. $v=a/(q-1)$. In this way, we obtain a two-variable polynomial $\tilde{P}_G(q,a)$, and the one-variable specializations we will def\mbox{}ine to be $Q_{G,n}(q)=\tilde{P}_{G}(q,q^n)$, for every integer $n \le 2$. Note that we then have $v=q^n/(q-1)=q^{n-1}/(1-q^{-1})$. We will denote the expression $q^n/(q-1)=q^{n-1} \sum_{i \ge 0}{q^{-i}}$ by $q_n$.\\
\indent We also obtain an analogous state-sum expression like for $P_{G,n}$, and then, in Section \ref{l4}, we def\mbox{}ine a cubic complex along the same lines as outlined before (for the f\mbox{}irst set of specializations).  

\section{The cubic complex construction of the chain complex} \label{l3}

\indent The construction of the complex (for both sets of specializations) will be similar to the case of the Jones polynomial for links (see Section \ref{kovnot}). So, we will have the same cubic complex and  we will use the same def\mbox{}initions of graded dimensions and graded chain complexes as in Section \ref{kovnot}. \\
\indent We are going to assign to each state a graded $\Z$-module whose graded dimension is $\kvn^{k}$ for the f\mbox{}irst set of specializations, and $q_n^k$ for the second set. The reader should note the following examples of graded modules:
\begin{example} \label{e1}
Let $V$ be the graded free $\Z$-module with $n+1$ basis elements: $X^0(=1),X,X^2,\ldots,X^n$ such that the degree of $X^i$ is equal to $n-i$, for $i=0,\ldots,n$. Then we have $q\dim V=1+q+\ldots+q^n=\kvn$. Furthermore, we have that $V^{\otimes k}\{l\}=q^l \kvn^k.$
Notice that we can also describe $V$ as the quotient of the ring of polynomials by $V=\Z[X]/(X^{n+1})\{n\}$, with $\deg X=-1$.
\end{example}
\begin{example} \label{e2}
Let $M=\Z[x_i]\{n-1\}$ be the ring of polynomials in one variable. If we def\mbox{}ine the degree of the variable $x_i$ to be $-1$, then $M$ becomes a graded free $\Z$-module, with quantum dimension $q\dim M =q^{n-1}\sum_{i\ge 0}{q^{-i}}=q_n$. Furthermore, if $k>0$, then $M^{\otimes k}$ is isomorphic to the ring of polynomials in $k$ variables $M_k=\Z[x_1,\ldots,x_k]\{k(n-1)\}$, and we have that $q\dim M^{\otimes k} = q\dim M_k = q_n^k$.  
\end{example}

\indent We are now ready to explain our construction. Let $n$ be a positive integer. Let $G$ be a graph with $m$ ordered edges and let $V$ be as in Example \ref{e1}. To each vertex $\epsilon=(\epsilon_1,\ldots,\epsilon_m)$ of the cube $\{0,1\}^m$, we associate a graded free $\Z$-module $V_{\epsilon}$ given by 
$V_{\epsilon}(G)=V^{\otimes k(\epsilon)}\{n|\epsilon|\}$, with $V$ being the graded $\Z$-module def\mbox{}ined in Example \ref{e1}. In other words, we assign one copy of $V$ to each connected component of the resolution of the graph, and then shift the degree up by $n|\epsilon|$. From the def\mbox{}inition we have that $q\dim V_{\epsilon}$ is (up to the sign $(-1)^{|\epsilon|}$) the contribution of $P_{G_{\epsilon}}$ to the dichromatic polynomial.\\

\begin{figure}[h]
\centerline
\relabelbox
\epsfysize 10cm
\epsfbox{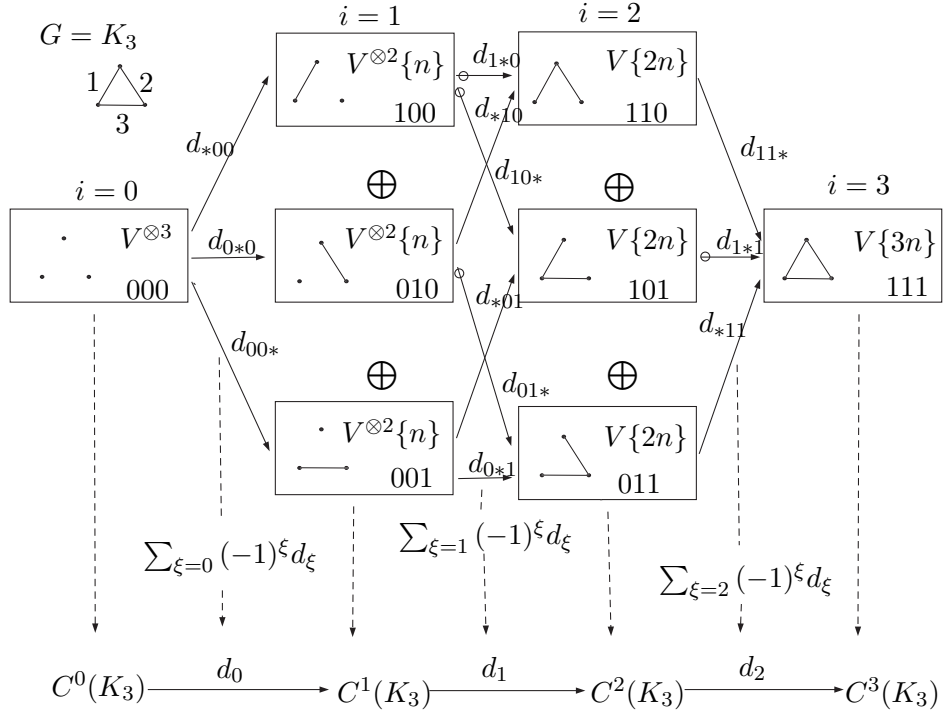}
\adjustrelabel <10pt,20pt> {G=K3}{$G=K_3$}
\adjustrelabel <-6pt,0pt> {1}{$1$}
\relabel {2}{$2$}
\adjustrelabel <-3pt,-5pt> {3}{$3$}
\relabel {000}{$000$}
\relabel {001}{$001$}
\relabel {010}{$010$}
\relabel {100}{$100$}
\relabel {011}{$011$}
\relabel {101}{$101$}
\relabel {110}{$110$}
\relabel {111}{$111$}
\adjustrelabel <-7pt,-7pt> {M3}{$V^{\otimes 3}$}
\adjustrelabel <-20pt,-7pt> {M21}{$V^{\otimes 2}\{n\}$}
\adjustrelabel <-20pt,-7pt> {M22}{$V^{\otimes 2}\{n\}$}
\adjustrelabel <-20pt,-7pt> {M23}{$V^{\otimes 2}\{n\}$}
\adjustrelabel <-10pt,-6pt> {M1}{$V\{2n\}$}
\adjustrelabel <-11pt,-7pt> {M2}{$V\{2n\}$}
\adjustrelabel <-6pt,-7pt> {M4}{$V\{2n\}$}
\adjustrelabel <-12pt,-7pt> {M5}{$V\{3n\}$}
\relabel {+1}{$\bigoplus$}
\relabel {+2}{$\bigoplus$}
\relabel {+3}{$\bigoplus$}
\relabel {+4}{$\bigoplus$}
\relabel {i=0}{$i=0$}
\relabel {i=1}{$i=1$}
\relabel {i=2}{$i=2$}
\relabel {i=3}{$i=3$}
\adjustrelabel <-8pt,0pt> {C0(K3)}{$C^0(K_3)$}
\adjustrelabel <-6pt,0pt> {C1(K3)}{$C^1(K_3)$}
\adjustrelabel <-6pt,0pt> {C2(K3)}{$C^2(K_3)$}
\adjustrelabel <-5pt,0pt> {C3(K3)}{$C^3(K_3)$}
\relabel {d0}{$d_0$}
\relabel {d1}{$d_1$}
\relabel {d2}{$d_2$}
\adjustrelabel <-5pt,0pt> {d*00}{$d_{\ast 00}$}
\adjustrelabel <0pt,2pt> {d0*0}{$d_{0\ast 0}$}
\relabel {d00*}{$d_{00\ast}$}
\relabel {d1*0}{$d_{1\ast 0}$}
\relabel {d*10}{$d_{\ast 10}$}
\relabel {d10*}{$d_{10\ast}$}
\relabel {d*01}{$d_{\ast 01}$}
\relabel {d01*}{$d_{01\ast}$}
\relabel {d1*1}{$d_{1\ast 1}$}
\relabel {d0*1}{$d_{0\ast 1}$}
\relabel {d11*}{$d_{11\ast}$}
\adjustrelabel <-3pt,0pt> {d*11}{$d_{\ast 11}$}
\adjustrelabel <-30pt,0pt> {S1}{$\sum_{\xi=0}{(-1)^{\xi}d_{\xi}}$}
\adjustrelabel <-30pt,0pt> {S2}{$\sum_{\xi=1}{(-1)^{\xi}d_{\xi}}$}
\adjustrelabel <-30pt,0pt> {S3}{$\sum_{\xi=2}{(-1)^{\xi}d_{\xi}}$}
\endrelabelbox
\caption{The chain complex and the dif\mbox{}ferentials}
\label{sl2}
\end{figure}

\indent To get the chain groups we ``f\mbox{}latten" the cube by taking direct sums along the columns. More precisely:
\begin{definition}
We set the $i$-th chain group $C_n^i(G)$ of the chain complex $\C_n(G)$ to be the direct sum of all $\Z$-modules at height $i$, i.e. $C_n^i(G)=\oplus_{|\epsilon|=i}{V_{\epsilon}(G)}$.\\
The grading is given by the degree of the elements and we can write the $i$-th chain group as $C_n^i(G)=\sum_{j\ge 0}C^{i,j}_n(G)$, where $C^{i,j}_n(G)$ denotes the set of elements of degree $j$ of $C_n^i(G)$.
\end{definition}

So, we are left with def\mbox{}ining the grading preserving dif\mbox{}ferential, since then according to Proposition \ref{p1} and Remark \ref{rem1} (see subsection \ref{subs1}), the chain complex obtained will have the $n$-specialization of the dichromatic polynomial as its Euler characteristic.

\subsection{The dif\mbox{}ferential}

\indent We def\mbox{}ine the dif\mbox{}ferential in the, now, standard way. We are f\mbox{}irst going to def\mbox{}ine per-edge maps between some vertices of the cube $\{0,1\}^m$ -- the maps that correspond to the edges of the cube. We def\mbox{}ine them as  linear, degree preserving maps and such that the cube is commutative, i.e.  every square is commutative. Then we build the dif\mbox{}ferential by summing with appropriate signs along columns, and hence obtain a map whose square is zero. This is presented graphically in Figure \ref{sl2}.\\
\indent So, let us f\mbox{}irst def\mbox{}ine per-edge maps. Each vertex of the cube $\{0,1\}^m$ is labeled with some $\epsilon=(\epsilon_1,\ldots,\epsilon_{m})\in\{0,1\}^m$. There are maps between two vertices only if one of the markers $\epsilon_i$ is changed from $0$ to $1$ when one goes from the f\mbox{}irst vertex to the second vertex and all the other $\epsilon_i$ are unchanged.\\
\indent Denote by $\epsilon$ the label of the f\mbox{}irst vertex. If the marker which is changed from $0$ to $1$ has index $j$ then the map will be denoted by $d_{\epsilon'}$, where $\epsilon'=(\epsilon'_1,\ldots,\epsilon'_{m})$ with $\epsilon'_i=\epsilon_i$ if $i\ne j$ and $\epsilon'_i=\ast$ if  $i=j$. Denote by $|\epsilon|$ the sum of components of $\epsilon$, and by $l$ the number of connected components of $G_{\epsilon}$.
Changing exactly one marker from 0 to 1 corresponds to adding an edge. That edge can either connect two dif\mbox{}ferent connected components, or it can connect one component to itself. \\
\indent In the f\mbox{}irst case, we have that two components are merging into one and the remaining $l-2$ components are left the same. We have to def\mbox{}ine a degree preserving map from $V^{\otimes l}\{n|\epsilon|\}$ to  $V^{\otimes (l-1)}\{n(|\epsilon|+1)\}$. We def\mbox{}ine it to be the identity on the $l-2$ tensor factors corresponding to the connected components which are left unchanged and on the remaining tensor product of  two copies of $V$ we have to def\mbox{}ine a map from $V\otimes V$ to $V\{n\}$. We def\mbox{}ine it on the basis vectors as: 
$$X^i \otimes X^j \mapsto X^{i+j},\quad i,j=0,\ldots,n,$$ 
where we assume $X^i=0$, for $i>n$.\\
\indent In the second case, the added edge belongs to one connected component, say the $p$-th one, and the remaining components will remain unchanged. In this case, we have to def\mbox{}ine a degree preserving map from $V^{\otimes l}\{n|\epsilon|\}$ to $V^{\otimes l}\{n(|\epsilon|+1)\}$. We def\mbox{}ine the map as the identity on those $l-1$ copies of $V$ which correspond to the connected components which are not af\mbox{}fected by adding the new edge, while on the remaining copy of $V$ (the one which corresponds to the $p$-th connected component), we have to def\mbox{}ine a linear map from $V$ to $V\{n\}$. We def\mbox{}ine it on the basis vectors by  $X^i \mapsto 0$, for $i=1,\ldots,n$, and $1 \mapsto X^{n}$.\\

\begin{remark} 
We could also def\mbox{}ine the dif\mbox{}ferential in the second case as the zero map (i.e. $1 \mapsto 0$ as well). This would also def\mbox{}ine a degree preserving map, and the theory would not be trivial since the dif\mbox{}ferential in the f\mbox{}irst case (the added edge connects two dif\mbox{}ferent components) is nontrivial.
\label{r1}
\end{remark}

 \indent The per-edge dif\mbox{}ferentials obtained in this way obviously make the cube commutative. We will sprinkle signs around on the edges of the cube and thus obtain a map whose square is zero. Namely, we def\mbox{}ine 
the dif\mbox{}ferential $d^i:C_n^i(G)\rightarrow C^{i+1}_n(G)$ on the chain complex $\C_n(G)$ by $$d^i:=\sum_{|\epsilon|=i}{(-1)^{\epsilon}d_{\epsilon}},$$ 
where we sum over all $m$-tuples of $0$'s, $1$'s, and exactly one $\ast$. Here we denoted the number of $1$'s in $\epsilon$ by $|\epsilon|$, and $(-1)^{\epsilon}$ equals $-1$ if there is odd number of $1$'s before $\ast$ in $\epsilon$, and $+1$ if there is even number of them. In Figure \ref{sl2} we indicated the dif\mbox{}ferentials with minus sign by adding a small circle at the tail of the arrow.\\

\indent A straightforward calculation implies:

\begin{proposition}
This def\mbox{}ines a dif\mbox{}ferential, that is, $d^{i+1}d^i=0$.
\end{proposition}

\indent Now we really have a chain complex $\C_n(G)$ where the chain groups and the dif\mbox{}ferential are def\mbox{}ined as in the previous two paragraphs, and according to Proposition \ref{p1}, we have

\begin{theorem}
The Euler characteristic of the chain complex $\C_n(G)$ is equal to the $n$-specialization $P_{G,n}$ of the dichromatic polynomial of the graph $G$.
\label{t1}
\end{theorem}

\indent Even though our construction depends on the ordering of the edges of the graph, in exactly the same way as in \cite{gr}, section 2.2.3, we obtain the following:

\begin{theorem}
The homology groups of the chain complex $\C_n(G)$ are graph invariants.
\label{t2}
\end{theorem}

\indent For each graph $G$ and positive integer $n$, we def\mbox{}ine the (two-variable) Poincar\'e polynomial of the corresponding chain complex
$\C_n(G)$, in the standard way, by $R_{G,n}(q,t)=\sum_{i\in \Z}{t^i q\dim (H^i)}$. Then from Theorems \ref{t1} and 
\ref{t2} we obtain:

\begin{proposition}
a) For each positive integer $n$, the polynomial $R_{G,n}$ is a graph invariant.\\
b) The $n$-specialization of the dichromatic polynomial $P_{G,n}(q)$ is equal to $R_{G,n}(q,-1)$. 
\end{proposition}
 
\subsection{Some calculations}

In this subsection we give the values of our homology theory for the polygon graph $\P_k$ with $k$ vertices and $k$ edges. The results from the following theorem are for the homology theory where the dif\mbox{}ferentials in the case when the added edge connects one connected component to itself is the zero map (see Remark \ref{r1}). In the general case, one can easily obtain an analogous result.

\begin{theorem}
{\em Odd case:} If $k=2g+1$, the free part of the homology is given by the Poincar\'e polynomial:
\begin{eqnarray*}
&(1+q+\ldots+q^{n-1})^{2g+1}+   \\
&(1+q+\ldots+q^{n-1})(\sum_{i=1}^{g-1}{(t^{2i-1}+t^{2i})q^{gn-g+i(n+1)}}+t^{2g-1}q^{2gn})+\\
& t^{2g}q^{2gn}(1+q+\ldots+q^n)+t^{2g+1}q^{(2g+1)n}(1+q+\ldots+q^n). 
\end{eqnarray*}
The torsion part is given by:
$$Tor(H^{\ast,\ast}_n(\P_{2g+1}))=\bigoplus_{i=1}^g {H^{2i-1,gn-g-1+i(n+1)}_n(\P_{2g+1})},$$
with each summand being isomorphic to $\Z_{n+1}$.\\
{\em Even case:} If $k=2g+2$, the free part of the homology is given by the Poincar\'e polynomial:
\begin{eqnarray*}
&(1+q+\ldots+q^{n-1})^{2g+2}+\\
&(1+q+\ldots+q^{n-1})(\sum_{i=0}^{g-1}{(t^{2i}+t^{2i+1})q^{gn+n-g+i(n+1)}}+t^{2g}q^{(2g+1)n})+\\
&t^{2g+1}q^{(2g+1)n}(1+q+\ldots+q^n)+t^{2g+2}q^{(2g+2)n}(1+q+\ldots+q^n). 
\end{eqnarray*}
The torsion part is given by:
$$Tor(H^{\ast,\ast}_n(\P_{2g+2}))=\bigoplus_{i=1}^g {H^{2i,(g+1)(n-1)+i(n+1)}_n(\P_{2g+2})},$$
with each summand being isomorphic to $\Z_{n+1}$.
\end{theorem}
\textbf{Proof:}\\
\indent First of all, note that up to the homological degree $k-2$, all the added edges connect two dif\mbox{}ferent components, and hence their contribution to the dif\mbox{}ferential is only of that type (algebra multiplication). Also, that part of the construction coincides with the construction \cite{gr2} and \cite{grp} for the algebra $\A_{n+1}$, and the only dif\mbox{}ference is in the grading. It is easy to see that the element of $q$-degree $i$ in \cite{gr2} becomes the element of $q$-degree $(2g+1)n-i$ in our theory, and so by applying the result from \cite{pr}, Corollary 4.3, we obtain the homology groups $H^i(\P_k)$, for $i=0,\ldots,k-2$. \\
\indent Since we took the case when the second type of the dif\mbox{}ferentials are zero maps, we have that $H^k(\P_k)=V\{kn\}$. Finally, the remaining cohomology group $H^{k-1}(P_k)$ we determine by using the fact that the restriction $t=-1$ (graded Euler characteristic) is equal to the $n$-specialization $P_{G,n}$ of the dichromatic polynomial. \kraj

\begin{remark}
We could also obtain the result of the previous theorem by using the Hochschild homology approach directly, like in \cite{pr}.
\end{remark}

\section{An ``inf\mbox{}inite-dimensional" set of specializations} \label{l4}

The set of specializations of the dichromatic polynomial that we have def\mbox{}ined in the previous section is enough to recover the original, two-variable dichromatic polynomial, but it ``misses" the specialization that corresponds to the Jones polynomial. In this section we def\mbox{}ine an alternative set of one-variable specializations (one polynomial for every integer $n \le 2$) and the corresponding set of chain complexes that categorify them. Although the categorif\mbox{}ication of this set of specializations has chain groups of inf\mbox{}inite rank, the advantage is that the case $n=2$  is equal (up to a multiple) to the Jones polynomial of an alternating link that corresponds to the planar graph.\\
\indent For every integer $n \le 2$, we def\mbox{}ine the polynomial $Q_{G,n}$ by
$$Q_{G,n}(q)=P_G(q,q^n/(q-1)).$$
Note that in this case we have $v=q^n/(q-1)=q^{n-1} \sum_{i \ge 0}{q^{-i}}=q_n$ in axiom $(A2)$.\\
\indent The categorif\mbox{}ication of $Q_n$ follows the same lines as the categorif\mbox{}ication in Section \ref{l3}. Namely, we organize the summands of the state-sum expression for $Q_n$ in the same cube as in  f\mbox{}igure \ref{sl1}, with the replacement of $\kvn$ by $q_n$ and of $q^n$ by $q$. Furthermore, to each $\epsilon \in \{0,1\}^m$ (i.e. to each vertex of the cube of resolutions) we assign a $\Z$-module $M_{\epsilon}$ in the following way: First of all, we order the vertices of $G$, and to the $i$-th one, we assign the index $i$. Let $K_1$, $K_2$,$\ldots$, $K_l$ be the set of connected components of $G_{\epsilon}$. Then to the component $K_j$ we will assign the variable $x_{i_j}$, where $i_j$ is the smallest index among the vertices $v$ that belong to $K_j$. Finally, we def\mbox{}ine $M_{\epsilon}(G)=\Z[x_{i_1},\ldots,x_{i_l}]\{l(n-1)+|\epsilon|\}$.\\
\indent In other words, to each connected component of $G_{\epsilon}$ we assign a copy of $M=\Z[x]\{n-1\}$, then take the tensor product and f\mbox{}inally increase the degree by $|\epsilon|$. From the def\mbox{}inition we have that $q\dim M_{\epsilon}(G)$ is (up to the sign $(-1)^{|\epsilon|}$) the contribution of $Q_{G_{\epsilon}}$ to the dichromatic polynomial (see Example \ref{e2} in Section \ref{l3}).\\
\indent In terms of pictures, we make the same picture as in f\mbox{}igure \ref{sl2}, only replacing $V$ by $M$ and all shifts $\{ni\}$ by $\{i\}$. Finally the chain groups $D^i_n(G)$ are obtained by summing along the columns, i.e. $D^i_n(G)=\oplus_{|\epsilon|=i}M_{\epsilon}(G)$.\\
\indent Now, we move to the dif\mbox{}ferential. The prescription is the same as in the previous section. We only have to def\mbox{}ine per-edge maps and make the cube commutative, and then the signs (def\mbox{}ined in the same way as previously) will make the cube anti-commutative, which after summing over columns def\mbox{}ines the map whose square is equal to zero, i.e. the dif\mbox{}ferential.\\
\indent Again, the edges of the cube of resolutions correspond to adding an edge of the graph $G$. This can af\mbox{}fect the connected components of the resolution in two ways: either the added edge connects two components (and consequently decreases the number of components by 1) or the added edge belongs to one of the components (and hence preserves the number of components).\\
\indent If adding the edge decreases the number of components by one, then we have that two components (say the $p$-th and the $q$-th one) merge into one, and the remaining $l-2$ components are left the same. Suppose that $i_p<i_q$. We have to def\mbox{}ine a degree preserving map from $M_{\epsilon}=\Z[x_{i_1},\ldots,x_{i_l}]\{l(n-1)\}$ to $M_{\epsilon'}=\Z[x_{i_1},\ldots,\hat{x}_{i_q},\ldots,x_{i_l}]\{(l-1)(n-1)+1\}$, where by $\hat{x}_{i_q}$ we have indicated that $x_{i_q}$ is omitted. We def\mbox{}ine it on the basis monomials by:
$$\prod_{j=1,\ldots,l}{x_{i_j}^{a_j}} \mapsto \prod_{{j=1,\ldots,l} \atop {j\ne q}}{x_{i_j}^{b_j}},$$
where $b_i=a_i$, for $i\ne p$, and $b_p=a_p+a_q+2-n$. In other words, we map $x_{i_q}$ to $x_{i_p}$ and multiply the result by $x_{i_p}^{2-n}$.\\
\indent If adding that edge does not af\mbox{}fect the number of components, then it means that the added edge will belong to some component, say the $p$-th one, and that the remaining components will remain unchanged. In this case, we have to def\mbox{}ine a grading preserving map from $M_{\epsilon}=\Z[x_{i_1},\ldots,x_{i_l}]\{l(n-1)\}$ to $M_{\epsilon'}=\Z[x_{i_1},\ldots,x_{i_l}]\{l(n-1)+1\}$. We def\mbox{}ine this map as the multiplication by $x_{i_p}$. \\

\begin{remark}
Like in the case of the f\mbox{}irst set of specializations from the previous Section, we could def\mbox{}ine the dif\mbox{}ferential in the second case (number of connected components is preserved) as the zero map. Since in the f\mbox{}irst case the dif\mbox{}ferential is nontrivial, this again def\mbox{}ines a chain complex with nontrivial degree-preserving dif\mbox{}ferentials.
\end{remark}
\indent These per-edge maps obviously make the cube commutative and after def\mbox{}ining signs in completely the same way as before, and adding over the columns (as in f\mbox{}igure \ref{sl2}), we obtain the dif\mbox{}ferential $d$ of our chain complex $D_n(G)$. Again, as before we have that the homology groups $H^i(D_n(G))$ do not depend on the ordering of edges of the graph $G$, and since the dif\mbox{}ferential is degree preserving, we have:

\begin{theorem}
The homology groups of the chain complex $D_n(G)$ are graph invariants. The graded Euler characteristic of the complex $D_n(G)$ is equal to the $n$-th specialization $Q_n(G)$ of the dichromatic polynomial of the graph $G$.
\end{theorem}

\begin{remark}
We could also def\mbox{}ine per-edge dif\mbox{}ferentials in the f\mbox{}inite\-dimensional case in Section \ref{l3} in the same manner as above, by using the alternative description of $V$ as a quotient of the ring of polynomials. Namely if the added edge connects two components, then the dif\mbox{}ferential is given on the corresponding tensor product of two copies of $V$ by multiplication of polynomials (and setting $X_q$ to be $X_p$), and in the other case is given by multiplication by $X_p^n$ (or as the zero map).
\end{remark}

\begin{remark}
Since both our constructions of categorif\mbox{}ication follow the same lines of the cubic complex construction, as it is used in other graph or link homologies (\cite{bn}, \cite{kov}, \cite{kovroz}, \cite{gr}, \cite{gr2}, \cite{grp}, \cite{v}, \cite{lee2}), we also obtain the long exact sequence of homology groups:
$$\cdots \rightarrow H^{i-1}(G/e) \rightarrow H^i(G) \rightarrow H^i(G-e) \rightarrow H^i(G/e) \rightarrow H^{i+1}(G) \rightarrow \cdots $$ 
\indent The analogous long exact sequence was f\mbox{}irst explicitly observed in \cite{v}, but it is also implicit in \cite{kov}, since chain complexes in all link (and graph) homology constructions are in fact mapping cones of certain homomorphisms, and hence we always obtain the long exact sequence in homology.
\label{rem2}
\end{remark}

\section{Categorif\mbox{}ication of the Jones polynomial for alternating links}
\label{l5}

As is well-known, every planar graph is in bijective correspondence with some knot universe (knot shadow), which in turn gives rise to an alternating link diagram, up to mirror image. As it is also known (e.g. \cite{bol}), the Jones polynomial of an alternating link is related (equal up to a multiple) with the value of the Tutte polynomial $T_G(x,1/x)$ of a planar graph corresponding to it. As we saw in the introduction the value of this specialization of the Tutte polynomial is (up to a multiple) equal to $P_G(q,q^2/(q-1))$ (with $q=1-1/x$). This can also be seen directly (see  Appendix \ref{a1} and e.g. \cite{kauf}). Hence, we can obtain a categorif\mbox{}ication of the Jones polynomial of an alternating link, by categorifying the one-variable specialization $J_G(q)=P_G(q,q^2/(q-1))$ of the dichromatic polynomial.
As we saw in the previous Section, we managed to categorify an inf\mbox{}inite set of specializations of the dichromatic polynomial (enough to recover it), and for $n=2$ we obtain that the specialization coincides with $J_G(q)$. However, in this section we will give an explicit categorif\mbox{}ication for this graph (Laurent) polynomial in terms of enhanced states.\\
\indent First of all, notice that we can write the value $v$ (which in the case of $J_G(q)$ is equal to $q^2/(q-1)$) as the following series:
\begin{equation}
v=q/(1-q^{-1})=q\sum_{i\ge 0}q^{-i}=\sum_{i\ge 0}q^{1-i},
\label{v}
\end{equation}
and with that value of $v$, the polynomial $J_G(q)$ is given by 
$$J_G(q)=\sum_{s\subset E(G)}{(-1)^{|s|}q^{|s|}v^{k(s)}}.$$ 
\indent In order to  def\mbox{}ine the graded complex, we will introduce the notion of enhanced states of a graph $G$. An enhanced state $S$ is given by a pair $S=(s,l)$, where $s$ is a state (subset of the set of edges $E(G)$), and $l$ is a labeling of the connected components of $[G:s]$ by nonnegative integers, i.e. a sequence of $k(s)$ nonnegative integers $l_1,\ldots,l_{k(s)}$. Denote by $|l|$, the sum of all $l_i$'s, $i=1,\ldots,k(s)$. To each enhanced state $S$ we will assign two numbers: $i(S)=|s|$, and $j(S)=|s|+k(s)-|l|$. Now, we def\mbox{}ine $C^{i,j}(G)$ (the $j$-th graded component of the $i$-th chain group, $0\le i \le m$, $j\in \Z$) as the span over $\Z$ of all enhanced states $S$ of $G$ such that $i(S)=i$ and $j(S)=j$.\\
\indent We will def\mbox{}ine a (degree preserving) dif\mbox{}ferential $d:C^{i,j}(G)\rightarrow C^{i+1,j}(G)$, by giving its value on each enhanced state $S=(s,l)\in C^{i,j}(G)$. We def\mbox{}ine:
\begin{equation}
d(S)=\sum_{e\in E(G)\setminus s}{(-1)^{n_e}}S_e,
\label{d}
\end{equation}
with $n_e$ being the number of edges in $s$ ordered before $e$, and where$S_e=(s_e,l_e)$ is an enhanced state whose state is $s_e=s \cup \{e\}$ and whose labeling $l_e$ is def\mbox{}ined as follows: denote by $E_1,\ldots,E_k$ the connected components of $[G:s]$.\\
\indent If $e$ connects some $E_i$ to itself then the components of $[G:s_e]$ are $E_1,\ldots,E_i\cup \{e\},\ldots,E_k$, and we def\mbox{}ine the labeling $l_e$ by $l_e(E_i\cup \{e\})=l(E_i)+1$ and $l_e(E_j)=l(E_j)$, for $j \ne i$.\\ 
\indent If $e$ connects some $E_i$ to $E_j$ (say $E_1$ and $E_2$), then the components of $[G:s_e]$ are $E_1\cup E_2 \cup \{e\}, E_3,\ldots, E_k$, and we def\mbox{}ine the labeling $l_e$ by $l_e(E_1 \cup E_2 \cup \{e\})=l(E_1)+l(E_2)$ and $l_e(E_j)=l(E_j)$, for $j>2$.\\

\begin{proposition}
The mapping $d$, def\mbox{}ined by (\ref{d}), satisf\mbox{}ies $d^2=0$.
\end{proposition}
\textbf{Proof:} \\
\indent Let $i$ and $j$ be arbitrary integers, and let $S$ be an enhanced state from $C^{i,j}$. If $e$ and $f$ are two distinct edges of $E(G)\setminus s$, then from the def\mbox{}inition of the states $S_e$ in the image of the dif\mbox{}ferential, we have $(S_e)_f=(S_f)_e$. Suppose that we have introduced an ordering $<$ on the set $E(G)$. Also, denote the number of edges in $s\cup \{f\}$ ordered before $e$ by $n'_e$, and the number of edges in $s\cup \{e\}$ ordered before $f$ by $n'_f$. Now we have:
\begin{eqnarray*}
d(d(S))&=&d(\sum_{e\in E(G)\setminus s}{(-1)^{n_e}}S_e)= \\
       &=&\sum_{e\in E(G)\setminus s}\left(\sum_{f \in E(G)\setminus s_e}{(-1)^{n'_f}}{(-1)^{n_e}}(S_e)_f\right)=\\
       &=&\sum_{{\scriptsize{\begin{array}{c}e,f \in E(G)\setminus s \\ e<f           \end{array}}}}{((-1)^{n'_f}(-1)^{n_e}+(-1)^{n'_e}{(-1)^{n_f}}})(S_e)_f=0,
\end{eqnarray*} 
as required. \kraj
\\
\\
\indent In this way we have obtained a sequence of chain complexes (one for each degree $j$). We could take the direct sum along columns (f\mbox{}ixed $i$), and obtain a (simply) graded chain complex $\C(G)$ with chain groups given by:
$$C^i(G)=\oplus_{j\in \Z}C^{i,j}(G).$$
The graded dimension of $C^i(G)$ is equal to the sum of the graded dimensions corresponding to each state $s$ with $|s|=i$. From the def\mbox{}inition of the gradings of our enhanced states we have that the graded dimension corresponding to each such state $s$ is equal to 
$$q^i{(\sum_{j\ge 0}{q^{1-j}})}^{k(s)}=q^iv^{k(s)},$$
with $v$ given by (\ref{v}).
Furthermore, since we have def\mbox{}ined degree preserving dif\mbox{}ferentials, we have the following:
\begin{theorem}
The graded Euler characteristic of the chain complex $\C(G)$ is equal to $J_G(q)$.
\end{theorem}  

\chapter{Koszul complexes and categorif\mbox{}ication of link and graph invariants}\label{Koszul}

\section{Triply graded link homology}\label{coment}

\subsection{Introduction}

\indent In this section we will introduce the parametrization of the HOMFLYPT polynomial that we will categorify. It is very similar to the one in \cite{KR2}. Throughout the chapter we will consider only braid diagrams $D$ of a link $L$, i.e. regular diagrams which are the closures of upward oriented braids (see Section \ref{braidnot}).\\
\indent  As is well known, every link can be represented by a braid diagram. Also,  the closures of two braid diagrams $D_1$ and $D_2$ are isotopic as oriented links if and only if $D_1$ and $D_2$ are related by a sequence of Markov moves, which are the following (see \cite{mar}):\\
\indent $(i)$ \quad conjugation: $DD'\longleftrightarrow D'D$\\
\indent $(ii)$ \quad transformations in the braid group:
\begin{eqnarray*}
D&\longleftrightarrow& D\sigma_i\sigma_i^{-1}\\
D&\longleftrightarrow& D\sigma_i^{-1}\sigma_i\\
D{\sigma_j}\sigma_i&\longleftrightarrow&D\sigma_i\sigma_j,\quad |i-j|>1\\
D\sigma_i\sigma_{i+1}\sigma_i&\longleftrightarrow&D\sigma_{i+1}\sigma_i\sigma_{i+1}\end{eqnarray*}
\indent $(iii)$\quad transformations $D\longleftrightarrow D\sigma_{n}^{\pm 1}$, for a braid $D$ with $n$ strands.\\
    
\indent In order to def\mbox{}ine the HOMFLYPT polynomial for a link $L$, from its braid diagram representation $D$, we will introduce a function $F$  on braid diagrams with values in the ring of rational functions in $q$ and $t$ def\mbox{}ined uniquely by the following axioms:\\
\indent $\bf{\ast}$  $F(D_1)=F(D_2)$, if $D_1$, $D_2$ are related by Markov move $(i)$\\
\indent $\bf{\ast}$  $F(D_1)=F(D_2)$, if $D_1$, $D_2$ are related by Markov moves $(ii)$\\
\indent $\bf{\ast}$  $F(D\sigma_n)=F(D)$, if a braid $D$ has $n$ strands\\
\indent $\bf{\ast}$  $F(D\sigma_n^{-1})=-t^{-1}q^{-1}F(D)$, if a braid $D$ has $n$ strands\\
\indent $\bf{\ast}$  Skein relation: for every braid diagram $D$ with $n$ strands and $0<i<n$  
$$ q^{-1}F(D\sigma_i)-qF(D\sigma_i^{-1})=(q^{-1}-q)F(D)$$ 
\indent $\bf{\ast}$ If $U$ is the one-strand diagram of the unknot then
$F(U)=1$.\\

\indent In order to obtain a link invariant we need to normalize the function $F$. Let $\alpha=-t^{-1}q^{-1}$ and let 
\begin{equation}
G(D)={\sqrt{\alpha}}^{n_{+}(D)-n_{-}(D)-s(D)+1}F(D),
\label{xxx}
\end{equation}
where $n_+(D)$, $n_-(D)$ and $s(D)$ are the number of positive crossings, negative crossings and the number of strands of $D$, respectively. We denote the expression $n_{+}(D)-n_{-}(D)-s(D)+1$ by $\omega(D)$. Obviously $G(D)$ is invariant under all Markov moves of braids and it satisf\mbox{}ies the HOMFLYPT skein relation
$$ {(q\sqrt{\alpha})^{-1}G(D\sigma_i)-q\sqrt{\alpha}G(D{\sigma_i}^{-1})}=(q^{-1}-q)G(D).$$
Hence, $G(D)$ is equal to the HOMFLYPT polynomial of the link $L$, normalized such that $G(U)=1$.
%\frac{1+t^{-1}q}{1-q^2}$.
In Section \ref{kr2a0}, we will def\mbox{}ine a triply graded chain complex $\C(D)$ whose Euler characteristic is equal to $F(D)$.\\
\indent First of all note  our slightly dif\mbox{}ferent convention compared to \cite{KR2} on the value of unknot. This has the advantage that we can  obtain  the  Alexander polynomial directly by specializing $t$ and $q$ ($t=-q$), and the whole sequence of the $n$-specializations of the (reduced) HOMFLYPT polynomial  (see \cite{kovroz}, \cite{moy}). Specif\mbox{}ically, by taking  $t=-q^{1-2n}$ we obtain polynomials $G_n(D)$ that satisfy the skein relation
$$q^{-n}G_n(D\sigma_i)-q^n G_n(D{\sigma_i}^{-1})=(q^{-1}-q)G_n(D), $$
and whose value on the unknot is $G_n(U)=1$. Hence by suitably collapsing the tri-grading to a bi-grading we get a new categorif\mbox{}ication of the $n$-specializations of the HOMFLYPT polynomial. 

\subsection{Graphs with wide edges}\label{grafwide}  

\indent In order to def\mbox{}ine the function $F(D)$ and hence the HOMFLYPT  polynomial $G(D)$, we introduce the trivalent graphs with wide edges (like in Chapter \ref{not}, Section \ref{homflynot}) as the resolutions of the crossings. Apart from the crossings $\sigma_i$ and $\sigma_i^{-1}$, we introduce wide edges $\bar{E}_i$ placed between the $i$-th and $(i+1)$-th strand of the braid, like in the following picture:

{{
\begin{center}
\setlength{\unitlength}{4mm}
\begin{picture}(30,6) 
\linethickness{0.6pt}
%\line(1.00,2.00,1.00,3.00)
\qbezier(5.00,1.00)(5.00,3.00)(5.00,5.00)
\qbezier(4.80,4.60)(4.90,4.80)(5.00,5.00)
\qbezier(5.20,4.60)(5.10,4.80)(5.00,5.00)

\qbezier(9.00,1.00)(9.00,3.00)(9.00,5.00)
\qbezier(8.80,4.60)(8.90,4.80)(9.00,5.00)
\qbezier(9.20,4.60)(9.10,4.80)(9.00,5.00)

\linethickness{2.0pt}
\qbezier(13.00,2.00)(13.00,3.00)(13.00,4.00)
\linethickness{0.6pt}
\qbezier(12.00,1.00)(12.50,1.50)(13.00,2.00)
\qbezier(14.00,1.00)(13.50,1.50)(13.00,2.00)
\qbezier(12.00,5.00)(12.50,4.50)(13.00,4.00)
\qbezier(14.00,5.00)(13.50,4.50)(13.00,4.00)

\qbezier(14.00,5.00)(13.80,4.90)(13.60,4.80)
\qbezier(14.00,5.00)(13.90,4.80)(13.80,4.60)

\qbezier(12.00,5.00)(12.20,4.90)(12.40,4.80)
\qbezier(12.00,5.00)(12.10,4.80)(12.20,4.60)

\qbezier(12.50,1.50)(12.30,1.40)(12.10,1.30)
\qbezier(12.50,1.50)(12.40,1.30)(12.30,1.10)

\qbezier(13.50,1.50)(13.70,1.40)(13.90,1.30)
\qbezier(13.50,1.50)(13.60,1.30)(13.70,1.10)

\qbezier(17.00,1.00)(17.00,3.00)(17.00,5.00)
\qbezier(16.80,4.60)(16.90,4.80)(17.00,5.00)
\qbezier(17.20,4.60)(17.10,4.80)(17.00,5.00)

\qbezier(21.00,1.00)(21.00,3.00)(21.00,5.00)
\qbezier(20.80,4.60)(20.90,4.80)(21.00,5.00)
\qbezier(21.20,4.60)(21.10,4.80)(21.00,5.00)

\put(6.00,3.00){$\cdot$}
\put(7.00,3.00){$\cdot$}
\put(8.00,3.00){$\cdot$}
\put(18.00,3.00){$\cdot$}
\put(19.00,3.00){$\cdot$}
\put(20.00,3.00){$\cdot$}

\put(5.00,5.50){1}
\put(11.95,5.50){$i$}
\put(13.50,5.50){$i+1$}
\put(21.00,5.50){$p$}
\put(12.50,0.00){$\bar{E}_i$} 
\end{picture}
\end{center}
}}

 Then we can def\mbox{}ine the function $F(D)$ by resolving the crossings by using the following two relations:
\begin{eqnarray}
F(D\sigma_i)&=& F(D\bar{E}_i)-q^2F(D),\label{1}\\
F(D{\sigma_i}^{-1})&=& q^{-2}F(D\bar{E}_i)-q^{-2}F(D).\label{2}
\end{eqnarray}  
Here by $F(D)$ we mean the value of the function $F$ on the diagram that is the closure of the braid diagram $D$, and we have extended the domain of the function $F$ to include trivalent graphs. Then $F$ (restricted to braid diagrams) will satisfy the axioms from the previous subsection, if and only if the values of $F$ on  completely resolved trivalent graphs satisfy
\begin{eqnarray}
&F(U)=1 \label{rel0}\\
&F(D\cup U) = \frac{1+t^{-1}q}{1-q^2} F(D), \quad{\textrm{if }} D {\textrm{ is not an empty diagram}} \label{rel1}\\
&F(D\bar{E}_n)=\frac{1+t^{-1}q^3}{1-q^2} F(D), \quad {\textrm{where }} D {\textrm{  is a diagram with  }} n {\textrm{ strands}}\\
&F(D\bar{E}_i^2)= (1+q^2)F(D\bar{E}_i)\\
&F(D\bar{E}_i\bar{E}_{i+1}\bar{E}_i)+q^2F(D\bar{E}_{i+1})=
F(D\bar{E}_{i+1}\bar{E}_{i}\bar{E}_{i+1})+q^2F(D\bar{E}_{i}).\label{rel4}
\end{eqnarray}

We will use (\ref{1}) and (\ref{2}) in resolving the crossings in the def\mbox{}inition of the triply-graded chain complex that categorif\mbox{}ies  the HOMFLYPT polynomial. For some more details about state-sum models of the HOMFLYPT polynomial, see Appendix \ref{ap2}.

\subsection{Categorif\mbox{}ication of the two-variable HOMFLYPT polynomial}\label{kr2a0}

In this subsection we will give an alternative, simpler and (essentially) equivalent construction to the one in \cite{KR2} (similar simplif\mbox{}ication is also implicit in \cite{kovhoh}). \\
\indent Essentially, we will set the variable $a$ from \cite{KR2} to be 0, but we will keep the double grading of the ring of polynomials $R'=\Q[x_1,\ldots,x_{2n}]$. Specif\mbox{}ically, to every arc (line between two crossings) we will assign a dif\mbox{}ferent variable $x_i$, $i=1,\ldots,2n$, where $n$ is the number of crossings of a given planar projection $D$ of a knot $K$. We def\mbox{}ine the bidegree of every $x_i$ to be (0,2) and we put the f\mbox{}ield of coef\mbox{}f\mbox{}icients $\Q$ in bidegree (0,0). Also, by $\{\cdot,\cdot\}$ we denote a shift in bigrading (like in the case of a simply graded space, see Section \ref{kovnot}).
\begin{remark}
Note that this corresponds to the case $n=-1$ in \cite{kovroz}, but with the introduction of a new grading direction.   
\end{remark}

Let $L$ be a link and let $D$ be its braid diagram presentation. Let $I$ be the ideal of $R'$ generated by the monomials $x_1+x_2-x_3-x_4$ for every crossing of $D$, and let $R=R'/I$ (see the picture below for the notation).\\
\indent   To each crossing we will assign 0- and 1-resolutions according to the following picture (the same as in Section \ref{kovroznot}):
\begin{center}
\setlength{\unitlength}{5mm}
\begin{picture}(25,5) 
\linethickness{0.6pt}
%\line(1.00,2.00,1.00,3.00)
\put(-0.5,0)
{
\put(0,-0.5){
\qbezier(4.00,1.00)(4.00,3.00)(4.00,5.00)
\qbezier(3.80,4.60)(3.90,4.80)(4.00,5.00)
\qbezier(4.20,4.60)(4.10,4.80)(4.00,5.00) \qbezier(2.00,1.00)(2.00,3.00)(2.00,5.00)
\qbezier(1.80,4.60)(1.90,4.80)(2.00,5.00)
\qbezier(2.20,4.60)(2.10,4.80)(2.00,5.00) 
\qbezier(9.5,3)(10.5,3)(11.5,3)
\qbezier(11.5,3)(11.3,3.1)(11.1,3.2)
\qbezier(11.5,3)(11.3,2.9)(11.1,2.8)

\qbezier(4.5,3)(5.5,3)(6.5,3)
\qbezier(4.5,3)(4.7,3.1)(4.9,3.2)
\qbezier(4.5,3)(4.7,2.9)(4.9,2.8)

\put(10,0){
\qbezier(9,1)(8,3)(7,5)
\qbezier(7,1)(7.4,1.8)(7.8,2.6)
\qbezier(9,5)(8.6,4.2)(8.2,3.4)
}

\linethickness{2.5pt}
\qbezier(13.00,2.00)(13.00,3.00)(13.00,4.00)
\linethickness{0.6pt}
\qbezier(12.00,1.00)(12.50,1.50)(13.00,2.00)
\qbezier(14.00,1.00)(13.50,1.50)(13.00,2.00)
\qbezier(12.00,5.00)(12.50,4.50)(13.00,4.00)
\qbezier(14.00,5.00)(13.50,4.50)(13.00,4.00)

\qbezier(14.00,5.00)(13.80,4.90)(13.60,4.80)
\qbezier(14.00,5.00)(13.90,4.80)(13.80,4.60)

\qbezier(12.00,5.00)(12.20,4.90)(12.40,4.80)
\qbezier(12.00,5.00)(12.10,4.80)(12.20,4.60)

\qbezier(12.50,1.50)(12.30,1.40)(12.10,1.30)\qbezier(12.50,1.50)(12.40,1.30)(12.30,1.10)

\qbezier(13.50,1.50)(13.70,1.40)(13.90,1.30)\qbezier(13.50,1.50)(13.60,1.30)(13.70,1.10)
\qbezier(14.5,3)(15.5,3)(16.5,3)
\qbezier(14.5,3)(14.7,3.1)(14.9,3.2)
\qbezier(14.5,3)(14.7,2.9)(14.9,2.8)

\qbezier(19.5,3)(20.5,3)(21.5,3)
\qbezier(21.5,3)(21.3,3.1)(21.1,3.2)
\qbezier(21.5,3)(21.3,2.9)(21.1,2.8)
\put(-10,0){
\qbezier(19,5)(18,3)(17,1)
\qbezier(17,5)(17.4,4.2)(17.8,3.4)
\qbezier(19,1)(18.6,1.8)(18.2,2.6)

\qbezier(19,5)(19,4.75)(19,4.5)
\qbezier(19,5)(18.8,4.8)(18.6,4.6)
\qbezier(17,5)(17,4.75)(17,4.5)
\qbezier(17,5)(17.2,4.8)(17.4,4.6)
}
\put(10,0){
\qbezier(9,5)(9,4.75)(9,4.5)
\qbezier(9,5)(8.8,4.8)(8.6,4.6)
\qbezier(7,5)(7,4.75)(7,4.5)
\qbezier(7,5)(7.2,4.8)(7.4,4.6)
}

\qbezier(22.00,1.00)(22.00,3.00)(22.00,5.00)
\qbezier(21.80,4.60)(21.90,4.80)(22.00,5.00)
\qbezier(22.20,4.60)(22.10,4.80)(22.00,5.00)

\qbezier(24.00,1.00)(24.00,3.00)(24.00,5.00)
\qbezier(23.80,4.60)(23.90,4.80)(24.00,5.00)
\qbezier(24.20,4.60)(24.10,4.80)(24.00,5.00)
\small{
\put(4.70,3.50){0-res.}
\put(9.6,3.50){1-res.}
\put(14.7,3.50){0-res.}
\put(19.6,3.50){1-res.}
}
\footnotesize{
\put(1.2,5.2){$x_1$}
\put(6.2,5.2){$x_1$}
\put(11.2,5.2){$x_1$}
\put(16.2,5.2){$x_1$}
\put(21.2,5.2){$x_1$}

\put(1.2,0.6){$x_4$}
\put(6.2,0.6){$x_4$}
\put(11.2,0.6){$x_4$}
\put(16.2,0.6){$x_4$}
\put(21.2,0.6){$x_4$}

\put(4.2,5.2){$x_2$}
\put(9.2,5.2){$x_2$}
\put(14.2,5.2){$x_2$}
\put(19.2,5.2){$x_2$}
\put(24.2,5.2){$x_2$}

\put(4.2,0.6){$x_3$}
\put(9.2,0.6){$x_3$}
\put(14.2,0.6){$x_3$}
\put(19.2,0.6){$x_3$}
\put(24.2,0.6){$x_3$}

}
}
\put(2.5,-1.5){$\Gamma_0$}
\put(22.5,-1.5){$\Gamma_0$}
\put(12.5,-1.5){$\Gamma_1$}
}
\end{picture}
\end{center}
\bigskip
\bigskip
and we call the resolutions obtained $\Gamma_0$ and $\Gamma_1$, as written in the picture. To the resolution $\Gamma_0$ we assign the following complex: 
$$\C(\Gamma_0):\quad\quad\quad 0\longrightarrow R\{-1,1\}\xrightarrow{x_2-x_3} R\longrightarrow 0$$
and to the resolution $\Gamma_1$ with the wide edge we assign the complex: 
$$\C(\Gamma_1):\quad\quad\quad 0\longrightarrow R\{-1,3\}\xrightarrow{(x_2-x_3)(x_4-x_2)} R\longrightarrow 0.$$

Assume that there are no free circles in the diagram $D$.
If we resolve all the crossings of $D$ we obtain a trivalent graph with wide edges. There are $2^n$ such resolutions $\Gamma$ of $D$ and to each one we assign the tensor product of $\C(\Gamma_0)$ and $\C(\Gamma_1)$, over all crossings $c(D)$, depending
on the type of resolution that appeared. In this way we have obtained a complex $\C(\Gamma)$ and to each resolution $\Gamma$ we will  assign  its cohomology $H(\Gamma)=H(\C(\Gamma))$.  \\
\indent Like in \cite{KR2}, we can obtain that $H(\Gamma)$ categorif\mbox{}ies the relations (\ref{rel1})--(\ref{rel4}). For example, the relation
(\ref{rel1}) becomes 
$$H(\Gamma\cup \textrm{ unknot})\cong
(H(\Gamma)\otimes\Q[x_i])\oplus 
(H(\Gamma)\otimes\Q[x_i]\{-1,1\}),$$
where $x_i$ is the label assigned to the circle (unknot). Note that in all def\mbox{}initions only the dif\mbox{}ferences $x_i-x_j$ appear. Thus we can work with the smaller ring of polynomials $R''=\Q[x_2-x_1,\ldots,x_{2n}-x_1]$ instead of $R'$ (like in \cite{KR2}). \\
\indent If we have free circles in the the diagram $D$, we introduce new variable $y$, with $\deg y=(0,2)$, extend the ring of polynomials to $R[y]$ and replace $R$ by $R[y]$ in the complexes $\C(\Gamma_i)$, $i=0,1$. Finally to every free circle we assign the complex: 
$$0\longrightarrow R'[y]\{-1,1\}\xrightarrow{y} R'[y]\longrightarrow 0,$$
and we tensor these complexes with $\C(\Gamma)$. In this way we obtain good value of the unknot (\ref{rel0}),  i.e. $H(U)\cong \Q$.\\
\indent We again organize the $2^n$ total resolutions $\Gamma$ of the diagram $D$ in the same cubic complex as in the standard categorif\mbox{}ications. To each vertex of the cube (i.e. to each total resolution $\Gamma$) we assign the graded vector space $H(\Gamma)$.

\indent We will introduce the dif\mbox{}ferentials between those cohomology groups as the maps induced by the (grading preserving) homomorphisms between the corresponding complexes $\C(\Gamma)$.
 Since these complexes are built as the tensor products of $\C(\Gamma_0)$ and $\C(\Gamma_1)$ it is enough to specify the homomorphisms between these two complexes. For a positive crossing $c$ we def\mbox{}ine  the following complex of complexes: 
\begin{equation}\C_c:\quad\quad 0\longrightarrow\C(\Gamma_0)\{0,2\}\stackrel{\chi_0}{\longrightarrow}\C(\Gamma_1)\longrightarrow 0,\label{dif1}\end{equation}
where $\C(\Gamma_1)$ is in cohomological degree 0, and the map $\chi_0$
is given by \\

\begin{equation*}
\begin{CD}
0@>>>R\{-1,3\}@>{x_2-x_3}>> R\{0,2\}@>>>0\\
 @.  @V{1}VV     @VV{{x_4-x_2}}V @.   \\
0@>>> R\{-1,3\}@>{(x_2-x_3)(x_4-x_2)}>>R@>>> 0.\end{CD}
\end{equation*}\\

%$$0\longrightarrow R\{-1,3\}\xrightarrow{x_2-x_3} %R\{0,2\}\longrightarrow %0$$
%$$\quad\quad{\scriptstyle{1}} %\downarrow\quad\quad\quad\quad\quad\quad\quad\quad\quad \downarrow %{\scriptstyle {x_4-x_2}}$$
%$$0\longrightarrow R\{-1,3\}\xrightarrow{(x_2-x_3)(x_4-x_2)} %R\longrightarrow 0.$$

\indent For a negative crossing $c$ we def\mbox{}ine  the following complex of complexes: 
\begin{equation}\C_c:\quad\quad 0\longrightarrow\C(\Gamma_1)\{0,-2\}\stackrel{\chi_1}{\longrightarrow}\C(\Gamma_0)\{0,-2\}\longrightarrow 0,\label{dif2}\end{equation}
where $\C(\Gamma_1)$ is in cohomological degree 0, and the map $\chi_1$
is given by \\
\begin{equation*}
\begin{CD}
0@>>>R\{-1,1\}@>{(x_2-x_3)(x_4-x_2)}>> R\{0,-2\}@>>>0\\
 @.  @V{x_4-x_2}VV     @VV{1}V @.   \\
0@>>> R\{-1,-1\}@>{x_2-x_3}>>R\{0,-2\}@>>> 0.\end{CD}
\end{equation*}\\

%$$0\longrightarrow R\{-1,1\}\xrightarrow{(x_2-x_3)(x_4-x_2)} %R\{0,-2\}\longrightarrow 0$$
%$${\scriptstyle{x_4-x_2}} %\downarrow\quad\quad\quad\quad\quad\quad\quad\quad\quad\downarrow %{\scriptstyle {1}}$$
%$$0\longrightarrow %R\{-1,-1\}\xrightarrow{x_2-x_3}R\{0,-2\}\longrightarrow 0.$$

\indent Def\mbox{}ine $\C(D)$ as the tensor product of $\C_c$ over all crossings $c(D)$. It is a complex built out of Koszul complexes $\C(\Gamma)$, over all the total resolutions $\Gamma$ of the diagram $D$, and the dif\mbox{}ferential preserves the bigrading of each term $\C^j(D)$. Every $\C^j(D)$ decomposes as a direct sum of contractible two-term complexes and its cohomology $H(\C^j(D))$, which is denoted by $\C H^j(D)$. The dif\mbox{}ferential induces the grading preserving maps $\delta$ from $\C H^j(D)$ to $\C H^{j+1}(D)$ and we denote the complex obtained in this way by $\C H(D)$. The cohomology $H(D)=H(\C H(D),\delta)$ is triply-graded: 
$$H(D)=\bigoplus_{j,k,l}{H_{k,l}^j(D)}.$$
Here $j$ is the cohomology degree, and $k$ and $l$ come from the internal bigrading of the chain groups.\\
\indent In complete analogy with \cite{KR2} we have that $H(D)$ does not depend on the choice of the braid presentation $D$ of a link $L$, up to  an overall shift in the triple grading. Also, as we saw since $H(\Gamma)$ categorif\mbox{}ies the relations (\ref{rel0})--(\ref{rel4}) and since the dif\mbox{}ferentials are induced by the grading preserving maps (\ref{dif1}) and (\ref{dif2}) which obviously categorify the relations (\ref{1}) and (\ref{2}), we have that the 
 bigraded Euler characteristic of $\C H(D)$ is equal to $F(D)$. Finally, by introducing half-integral shifts as in \cite{wu} (in order to compensate the powers of $\alpha$ from (\ref{xxx})) by:
$$\H(D)=H(D)[\omega(D)/2]\{-\omega(D)/2,-\omega(D)/2\},$$
we obtain a triply-graded homology theory, which does not depend on the choice of the braid presentation $D$ of a link $L$ and whose bi-graded Euler characteristic is equal to the two-variable HOMFLYPT polynomial of a link $L$.

\section{New categorif\mbox{}ications of the chromatic and  dichromatic polynomials for graphs}\label{newgraf}

\subsection{Introduction}

\indent In Chapter \ref{graf} we categorif\mbox{}ied the dichromatic polynomial of graphs by introducing one-variable specializations of it and then we categorif\mbox{}ied each of these one-variable polynomials.  In this section we will def\mbox{}ine a complex of doubly-graded modules whose doubly-graded Euler characteristic is equal to the whole two-variable dicromatic polynomial.  The idea is partially inspired by the categorif\mbox{}ication of HOMFLYPT polynomial described in the previous section. \\
\indent Also, we
give a new categorif\mbox{}ication of the chromatic polynomial for graphs. We do this here dif\mbox{}ferently to \cite{gr}. We will def\mbox{}ine the chain groups (the modules corresponding to the vertices of the cube of resolutions) as the cohomologies of certain chain complexes.\\
\indent We will use the same notation as in Chapter \ref{graf}: namely, a graph $G$ is specif\mbox{}ied by a set of vertices $V(G)$ and a set of
edges $E(G)$. If $e$ is an arbitrary edge of the graph $G$,
then by $G-e$ we denote the graph $G$ with the edge $e$ deleted, and
by $G/e$ the graph obtained by contracting edge $e$ (i.e. by
identifying the vertices incident to $e$ and deleting $e$). 

\subsection{The chromatic polynomial}

 If $q$ is a positive integer, the chromatic polynomial $P_G(q)$ is def\mbox{}ined as the number of ways to color the vertices of $G$ by using at most $q$ colors, such that every two vertices which are connected by an edge receive a dif\mbox{}ferent color.
It is well-known that the chromatic polynomial can be def\mbox{}ined equivalently by the following two  axioms:
\begin{eqnarray*}
(C1)&\quad P_G=P_{G-e}-P_{G/e},\\
(C2)&\quad P_{N_k}=q^k,
\end{eqnarray*}
where $N_k$ is the graph with $k$ vertices and no edges. By using these axioms we can obviously extend the domain of the polynomial to the set of complex numbers, and, furthermore, instead of $q$ in the axiom $(C2)$   we will put $1/(1-q)$, with $|q|<1$.\\
\indent By repeated use of $(C1)$ (which is the famous deletion-contraction rule) we will obtain the value of the chromatic polynomial as a sum of  contributions from all spanning subgraphs of $G$ (subgraphs that contain all vertices of $G$), which we will call states.
Furthermore, if for each subset $s\subset E(G)$ we denote by $[G:s]$ the graph whose set of vertices is $V(G)$ and set of edges is $s$, then the contribution of the graph $[G:s]$ is $(-1)^{|s|}(1-q)^{-k(s)}$, where $|s|$ is the number of elements of $s$ and $k(s)$ is the number of connected components of $[G:s]$. Hence, we obtain the expression:
$$P_G(q)=\sum_{s\subset E(G)}{(-1)^{|s|}(1-q)^{-k(s)}}=\sum_{i\ge 0}{(-1)^i \sum_{s\subset E(G),|s|=i}{(1-q)^{-k(s)}}},$$
which is called the state-sum expansion of the polynomial $P_G(q)$.\\
 \indent In  subsection \ref{iiglavni} we will def\mbox{}ine a graded chain complex of modules $\C(G)$ whose graded Euler characteristic is equal to $P_G(q)$.\\

\subsection{The dichromatic polynomial}

The dichromatic polynomial $P_G(q,v)$ of the graph $G$ is a two-variable generalization of the chromatic polynomial given by the following two axioms:
\begin{eqnarray*}
(D1)&\quad P_G=P_{G-e}-qP_{G/e},\\
(D2)&\quad P_{N_k}=v^k,
\end{eqnarray*}
where $N_k$ is the graph with $k$ vertices and no edges.\\
\indent From $(D1)$ we have a recursive expression for the dichromatic polynomial in terms of the value of the polynomial on graphs with a smaller number of edges. Indeed, as in the case of the chromatic polynomial we obtain that the contribution of the state $[G:s]$ is $(-1)^{|s|}q^{|s|}v^{k(s)}$, where $|s|$ is the number of elements of $s$ and $k(s)$ is the number of connected components of $[G:s]$. Hence, we obtain the expression:
$$P_G(q,v)=\sum_{s\subset E(G)}{(-1)^{|s|}q^{|s|}v^{k(s)}}=\sum_{i\ge 0}{(-1)^i q^i\sum_{s\subset E(G),|s|=i}{v^{k(s)}}},$$
which is called the state-sum expansion of the polynomial $P_G(q,v)$. However, we will use a slightly dif\mbox{}ferent parametrization of the dichromatic polynomial, given by:
$$D_G(t,q)=(1+t^{-1}q)^{m}P_G(q,\frac{1+t^{-1}q}{1-q}),$$
where 
%$n$ is the number of vertices of the graph $G$ and 
$m$ is the number of edges of the graph $G$.\\
\indent In subsection \ref{iiglavni2} we will def\mbox{}ine a chain complex $\D(G)$ of doubly graded modules whose doubly graded Euler characteristic is equal to $D_G(t,q)$.

\subsection{The categorif\mbox{}ication of the chromatic polynomial}\label{iiglavni}

Let $n$ denote the number of vertices of the graph $G$. Let $R$ be the ring of polynomials in $n$ variables over $\Q$, i.e. $R=\Q[x_1,\ldots,x_n]$. We introduce a grading in $R$, by giving the degree 1 to every $x_i$. Order the set of vertices of $G$ and to the $i$-th  vertex assign the variable $x_i$. Finally, to every edge $e\in E(G)$, whose endpoints are the vertices $i_e$ and $j_e$, assign the monomial $m_e=x_{i_e}-x_{j_e}$ (the ambiguity of the sign does not af\mbox{}fect the later construction).

\subsubsection{The cubic complex construction}
\indent Let $s\subset E(G)$ be a subset of the set of edges of $G$, and let $[G:s]$ be the corresponding state in the resolution of a graph $G$. Def\mbox{}ine the ideal $I_s$ as the ideal generated by the monomials $m_e$, for all edges $e\in s$. Finally, to the state $[G:s]$ assign the module $R_s=R/{I_s}$.\\
\begin{proposition}{\label{prop1}} The quantum graded dimension of $R_s$ is equal to
$(1-q)^{-k(s)}$, where $k(s)$ denotes the number of connected components of $[G:s]$.
\end{proposition}

\textbf{Proof:}\\
\indent Let $i$ and $j$ be two arbitrary vertices of $G$. They obviously belong to the same connected component of $[G:s]$ if and only if there exist a sequence of edges belonging to $s$ which connects $i$ and $j$, which obviously happens if and only if $x_i-x_j$ belongs to $I_s$. Hence, all the variables corresponding to the vertices from the same component, must be the same in $R_s$. In other words, $R_s$ is isomorphic to the ring of polynomials (over $\Q$) in $k(s)$ variables, and hence:
$$q\dim R_s={\left(\sum_{i\ge 0}{q^i}\right)}^{k(s)}={(1-q)}^{-k(s)}.$$ \kraj
\\
\\
\indent Denote by $m$ the number of edges of $G$, and f\mbox{}ix an ordering on the set $E(G)$, denoted by $(e_1,\ldots,e_m)$. Now we will def\mbox{}ine 
the chain complex $\C$ in a standard way, by ''summing over columns" of our cubic complex:
for each $i$, with $0\le i \le m$, we will def\mbox{}ine the $i$-th chain group, ${\C^i(G)}$ as the direct sum of $R_s$, over all $s\subset E(G)$, such that $|s|=i$. \\
\indent Now, let us turn to the dif\mbox{}ferential. We will def\mbox{}ine the map $d^i$ from ${\C^i(G)}$ to ${\C^{i+1}(G)}$ as a sum of maps between the direct summands of the chain groups. The only nonzero maps are the maps from $R_s$ to $R_{s\cup \{e\}}$, with $e \notin s$ (which are exactly the ones that correspond to the edges of the cube), and we set them (up to a sign) to be the identity (i.e. the map that sends $f+I_s$ to $f+I_{s\cup \{e\}}$ for every $f\in R$). \\
\indent We now introduce signs in a standard way in order to make the cube anticommutative, and hence to make the square of the dif\mbox{}ferential equal to 0. Namely, we put minus signs exactly for those maps $R_s\to R_{s\cup \{e\}}$, with an odd number of edges in $s$ which are ordered before $e$.  \\
\indent In this way we have obtained a chain complex, $\C(G)$, of  graded $R$-modules with grading preserving dif\mbox{}ferential. Its homology groups obviously don't depend of the ordering of the vertices, and also don't depend of the ordering of the edges of $G$ (like in \cite{gr}, section 2.2.3), and hence we obtain
\begin{theorem}  
The homology groups of the chain complex $\C(G)$ are invariants of the graph $G$, and the graded Euler characteristic of $\C(G)$ is equal to the chromatic polynomial $P_G(q)$.
\end{theorem}

\subsubsection{Alternative description}

Now we will give an equivalent def\mbox{}inition of the chain complex $\C(G)$ in terms of Koszul complexes.\\
\indent To each edge $e\in E(G)$ we assign two complexes, $\C_{e-}$ and $\C_{e+}$ def\mbox{}ined in the following way:
$$\C_{e-}: \quad 0 \longrightarrow R \buildrel{0}\over\longrightarrow R \longrightarrow 0,$$
$$\C_{e+}: \quad 0 \longrightarrow R \buildrel{x_i-x_j}\over\longrightarrow R \longrightarrow 0,$$
where $i$ and $j$ are the vertices of $G$ which are the endpoints of the edge $e$.
Now, to every subset $s\subset E(G)$ we assign a complex $\C_s$ which is the tensor product of $\C_{e\pm}$, where we take $+$ if $e\in s$ and $-$ if $e \notin s$. Finally, to the state $[G:s]$ we assign the cohomology of $\C_s$ at the rightmost position.\\
\indent To build the dif\mbox{}ferentials, we introduce the (grading preserving) maps $d_e$, as the maps induced on the cohomology by the following homomorphism from $\C_{e-}$ to $\C_{e+}$:\\
\begin{equation}
\begin{CD}
0@>>>R@>{0}>> R@>>>0\\
 @.  @V{0}VV     @VV{1}V @.  \\
0@>>> R@>{x_i-x_j}>>R@>>> 0.\end{CD}
\label{dif}\end{equation}\\

Here we put the upper row in cohomological degree 0, and the lower one in cohomological degree 1.\\
\indent Now, in order to def\mbox{}ine the dif\mbox{}ferentials, just tensor all the
chain complexes and maps between them from (\ref{dif}) over all edges $e$ of $E(G)$. If we take the cohomology only at the rightmost position in  each ``horizontal" complex (the ones in the same cohomological degree with respect to the def\mbox{}inition after (\ref{dif})), and as the dif\mbox{}ferentials are the induced maps between them, we obtain a complex  $\C'(G)$ which is isomorphic to the complex $\C(G)$ from the previous subsection.

\subsection{The categorif\mbox{}ication of the dichromatic polynomial}\label{iiglavni2}

In order to categorify the dichromatic polynomial we will have to introduce a new grading direction, and we will use the whole Koszul complex (actually a slightly modif\mbox{}ied one) that we have used in the previous subsection.\\
\indent We order the vertices of $G$, and to the $i$-th one ($1\le i \le n=\sharp V(G))$, we assign the variable $x_i$. We def\mbox{}ine the bidegree of all $x_i$ as (0,1). Def\mbox{}ine the bigraded ring $R$ by $R=\Q[x_1,\ldots,x_n]$, where we put the f\mbox{}ield $\Q$ in bidegree (0,0).\\
\indent To each edge $e$, such that its endpoints are the $i$-th and $j$-th vertex, we can associate two resolutions of the graph $G$: the f\mbox{}irst one with the edge $e$ contracted (i.e. when we identify the vertices $i$ and $j$), and the second one with the edge $e$ deleted. To the f\mbox{}irst resolution we assign the following complex (denoted by $\D(e+)$):
$$\D(e+):\qquad 0 \longrightarrow R\{-1,1\} \stackrel{x_i-x_j}{\longrightarrow} R \longrightarrow 0,$$
and to the second one we assign the complex $\D(e-)$ given by:
$$\D(e-):\qquad 0 \longrightarrow R\{-1,1\} \stackrel{0}{\longrightarrow} R \longrightarrow 0.$$
Let $s$ be an arbitrary subset of $E(G)$ and let $[G:s]$ be the corresponding state of $G$. Then, to that state we assign the (Koszul) complex $\D'(s)$ of bigraded $R$-modules obtained by tensoring the complexes $\D(e+)$, where $e$ runs over all edges in $s$, and $\D(f-)$, where $f$ runs over all edges in $E(G)\setminus s$. We denote its (bigraded) cohomology by $H'(s)$ (the direct sum of the cohomology groups of $\D'(s)$).
\begin{proposition}
The quantum bigraded dimension of $H'(s)$ is equal to:
$${(1+t^{-1}q)}^{m-n}{\left(\frac{1+t^{-1}q}{1-q}\right)}^{k(s)},$$
where $k(s)$ is the number of connected components of $[G:s]$, and $n$ and $m$ are the number of vertices and edges of $G$, respectively. 
\end{proposition}
\textbf{Proof:}\\
\indent Like in the proof of Proposition \ref{prop1} we obtain that the cohomology at the rightmost position of $\D'(s)$ is isomorphic to the ring of polynomials in $k(s)$ variables. However, here we will also have the cohomology at the leftmost position in each of the $\D(e\pm)$, which is isomorphic to the same ring of polynomials in $k(s)$ variables, but shifted by the bidegree $\{-1,1\}$, for all $\D(e-)$ and for a certain number of the $D(e+)$. We will show by induction on $|s|$ that the total number of such $e$'s, denoted by $c(s)$, is equal to $k(s)-n+m$.\\
\indent If $|s|=0$ then we have the tensor product of $m$ complexes with all the mappings equal to zero, and hence we have that the number of edges which contribute with nontrivial cohomology at the leftmost position is equal to $m=k(s)-n+m$ (note that in this case $k(s)=n$). Now suppose that the formula is true for some subset $s$ and consider the state $[G:{(s\cup e)}]$ with $e \in E(G) \setminus s$. Denote the endpoints of $e$ by $i$ and $j$, and denote $s'=s\cup e$. This means that  $\D'(s')$ is formed by the tensor product of the same complexes as $\D'(s)$ with $\D(e+)$ instead of $\D(e-)$. Now, $\D(e+)$ will have nontrivial cohomology at the leftmost position if and only if $x_i-x_j$ belongs to the ideal generated by the monomials def\mbox{}ined by the edges of $s$, i.e. if and only if the vertices $i$ and $j$ belong to the same connected component of $[G:s]$. In other words, we have $c(s')=c(s)$ if $k(s')=k(s)$ and $c(s')=c(s)-1$ if $k(s')=k(s)-1$. So we have $c(s)=k(s)-n+m$ as we wanted to prove.\\
\indent Hence the total bigraded dimension of $H'(s)$ is equal to:
$${(1+t^{-1}q)}^{k(s)-n+m}{(1-q)}^{-k(s)}.$$ \kraj

Furthermore, to every vertex $v$ of the graph $G$, we assign the same complex as $\D(e-)$. Now, if we tensor these complexes over all vertices of $G$ and tensor the complex obtained with $\D'(s)$, we obtain the complex $\D(s)$. We denote the cohomology of $\D(s)$ by $H(s)$, and that is the space that we will assign to the state $[G:s]$. Obviously, we have:
$$q\dim H(s)={(1+t^{-1}q)}^m{\left(\frac{1+t^{-1}q}{1-q}\right)}^{k(s)}. $$

\indent In order to introduce the dif\mbox{}ferentials between the cohomologies $H(s)$, we will def\mbox{}ine the (grading preserving) homomorphism $d(e)$ from $\D(e+)\{0,1\}$ to $\D(e-)$, and then for the dif\mbox{}ferentials we take the induced mappings on cohomology. We def\mbox{}ine $d(e)$ by:\\
\begin{equation*}
\begin{CD}
0@>>>R\{-1,2\}@>{x_i-x_j}>> R\{0,1\}@>>>0\\
 @.  @V{x_i-x_j}VV     @VV{0}V @.   \\
0@>>> R\{-1,1\}@>{0}>>R@>>> 0.\end{CD}
\end{equation*}\\

%$$0 \longrightarrow R\{-1,2\} %\stackrel{\scriptstyle{x_i-x_j}}{\longrightarrow} R\{0,1\} %\longrightarrow 0$$
%$${\scriptstyle{x_i-x_j}}\downarrow \quad \quad \quad \qquad\downarrow %0$$
%$$0 \longrightarrow R\{-1,1\} \stackrel{0}{\longrightarrow} \quad R %\quad\longrightarrow 0$$ 
We put the upper row in cohomological degree -1, and the lower row in cohomological degree 0. We will denote this complex of complexes 
by $\D_e$. \\
\indent For the graph $G$ def\mbox{}ine the complex of (Koszul) complexes by tensoring $\D_e$ over all edges $e$ of $G$. By taking the cohomology $H^j(G)$, $-m\le j\le 0$,  in each ``horizontal" complex and def\mbox{}ining the dif\mbox{}ferentials between them to be the ones induced by the tensor product of $d(e)$'s, we obtain a triply graded  complex $\D(G)$.\\
\indent From the def\mbox{}inition we have

\begin{theorem}
The homotopy class of the complex $\D(G)$ is an invariant of the graph $G$ whose bigraded Euler characteristic is equal to the  dichromatic polynomial $D_G(t,q)$ of the graph $G$.
\end{theorem}

\chapter{Properties of link homology for positive braid knots}
\label{braid}

\section{Introduction}

\indent In this chapter we prove the conjecture from \cite{pat} that the f\mbox{}irst Khovanov homology group of the positive braid knot is trivial (see Theorem \ref{braidgl}). In the proof we use the basic ingredients of the construction of the $sl(2)$ link homology (cubic complex, independence of the planar projection chosen), and so the major part can be directly applied to other link homology theories. Especially, in the case of $sl(n)$-link homology (see Section \ref{kovroznot} of Chapter \ref{not}), whose particular def\mbox{}initions of chain groups and dif\mbox{}ferentials make it practically incalculable, we modify slightly our proof to show that the f\mbox{}irst $sl(n)$ homology group of the positive braid knot is trivial for every positive integer $n$, as well.

\section{Positive braid knots}\label{posbraid}

The positive braid knots are the knots (or links) that are the closures of positive braids. Let $K$ be an arbitrary positive braid knot and let $D$ be its planar projection which is the closure of a positive braid. Denote the number of strands of that braid by $p$. We say that  the crossing $c$ of $D$ is of the type $\sigma_i$, $i<p$, if it corresponds to the generator $\sigma_i$ in the braid word of which $D$ is the closure.

{{
\begin{center}
\setlength{\unitlength}{4mm}
\begin{picture}(30,6) 
\linethickness{0.8pt}
%\line(1.00,2.00,1.00,3.00)
\qbezier(5.00,1.00)(5.00,3.00)(5.00,5.00)
\qbezier(9.00,1.00)(9.00,3.00)(9.00,5.00)

%\put(15,5){\line(-1,-1){4}}
\qbezier(15.00,5.00)(13.00,3.00)(11.00,1.00)
\qbezier(15.00,1.00)(14.30,1.70)(13.60,2.40)
\qbezier(11.00,5.00)(11.70,4.30)(12.40,3.60)

\qbezier(17.00,1.00)(17.00,3.00)(17.00,5.00)
\qbezier(21.00,1.00)(21.00,3.00)(21.00,5.00)

\put(6.00,3.00){$\cdot$}
\put(7.00,3.00){$\cdot$}
\put(8.00,3.00){$\cdot$}
\put(18.00,3.00){$\cdot$}
\put(19.00,3.00){$\cdot$}
\put(20.00,3.00){$\cdot$}

\put(4.90,5.50){1}
\put(10.90,5.50){$i$}
\put(14.40,5.50){$i+1$}
\put(20.90,5.50){$p$}
\put(12.70,0.00){{\Large $\sigma_i$}}

\end{picture}
\end{center}
}}

 Denote the number of crossings of the type $\sigma_i$ by $l_i$, $i=1,\ldots,p-1$ and order them from top to bottom. Then  we can write each crossing $c$ of $D$ as a pair $(i,\alpha)$ (we will also write $(i\alpha)$ if there is no possibility of confusion), $i=1,\ldots,p-1$ and $\alpha=1,\ldots,l_i$, if $c$ is of the type $\sigma_i$ and it is ordered as $\alpha$-th among the crossings of the type $\sigma_i$. Finally, we order the crossings of $D$ by the following ordering: $c=(i\alpha)<d=(j\beta)$ if and only if $i<j$, or $i=j$ and $\alpha<\beta$.\\
\indent Since a positive braid knot is a positive knot, we have that $\H^i(K)=0$ for $i<0$. Also if a positive braid knot $K$ has a regular diagram $D$, which is the closure of the $p$-strand braid $\sigma$ such that all letters $\sigma_i$, $i=1,\ldots,p-1$, appear in the braid word of $\sigma$ (i.e. $D$ is not the disjoint union of two nonempty diagrams), we know that its zeroth homology group (see e.g. \cite{pat}) is two-dimensional (without torsion) and that the $q$-gradings of the two generators are $1-p+n(D)\pm 1$, where $n(D)$ is the number of crossings of $D$.

\section{The main result}\label{iiiglavni} 
  
 In the following theorem we prove that the f\mbox{}irst homology group of a positive braid knot is trivial.
\begin{theorem}\label{braidgl}
If $K$ is a positive braid knot, then $\H^1(K)$ is trivial.
\end{theorem}
\textbf{Proof:}\\
\indent First of all, if $D$  is a regular diagram of $K$, which is the closure of a positive braid, then we have that $\H^{i,j}(K)=H^{i,j-n(D)}(D)$, where $n(D)$ is the number of crossings of $D$, and so we have that $\H^1(K)=0$ if and only if $H^1(D)=0$.  So we are going to show that the (unnormalized) f\mbox{}irst homology group of $D$ is trivial. \\
\indent In order to prove that $H^1(D)=0$, we will use the def\mbox{}inition of Khovanov homology (see Section \ref{kovnot} of Chapter \ref{not}), i.e. we will prove that for any element $t'$ of the  chain group $C^1(D)$ such that $d^1(t')=0$, there exists an element $t\in C^0(D)$ such that $t'=d^0(t)$. For this we f\mbox{}irst need to understand the chain groups $C^0(D)$, $C^1(D)$ and $C^2(D)$ and the dif\mbox{}ferentials  $d^0$ and $d^1$.\\
\indent $C^0(D)$ ``comes" from all resolutions $D_{\epsilon}$ with $|{\epsilon}|=0$. However, since $D$ is the closure of a positive braid, we have only one such resolution ${\epsilon}_0$ (all crossings are resolved into 0-resolutions) and it is  an unlink  that consists  of $p$ unknots (i.e. the closure of the trivial braid with $p$ strands). Hence, we have that $C^0(D)=V^{\otimes p}$ where we assigned the  $i$-th  copy of $V$ (denoted by $V^i$) to the circle  which is the closure of the $i$-th strand of $D_{{\epsilon}_0}$. \\
\indent Now, we pass to $C^1(D)$. It ``comes"  from all  resolutions $D_{\epsilon}$ with $|{\epsilon}|=1$, i.e. all resolutions where we resolve all except one crossing of $D$ in a 0-resolution   
 and the remaining one in a 1-resolution. 
In such a way, if the crossing $c$ that is resolved into a 1-resolution is of the type $\sigma_i$ (i.e. if $c=(i\alpha)$, for some $\alpha=1,\ldots,l_i$), then the corresponding resolution, $D_c$, is the closure of the plat diagram $E_i$ (see picture).

{{
\begin{center}
\setlength{\unitlength}{4mm}
\begin{picture}(30,6)
\linethickness{0.8pt}
%\line(1.00,2.00,1.00,3.00)
\qbezier(5.00,1.00)(5.00,3.00)(5.00,5.00)
\qbezier(9.00,1.00)(9.00,3.00)(9.00,5.00)

\qbezier(15.00,5.00)(13.00,2.00)(11.00,5.00)
\qbezier(11.00,1.00)(13.00,4.00)(15.00,1.00)

\qbezier(17.00,1.00)(17.00,3.00)(17.00,5.00)
\qbezier(21.00,1.00)(21.00,3.00)(21.00,5.00)

\put(6.00,3.00){$\cdot$}
\put(7.00,3.00){$\cdot$}
\put(8.00,3.00){$\cdot$}
\put(18.00,3.00){$\cdot$}
\put(19.00,3.00){$\cdot$}
\put(20.00,3.00){$\cdot$}

\put(5.00,5.50){1}
\put(11.00,5.50){$i$}
\put(14.50,5.50){$i+1$}
\put(21.00,5.50){$p$}
\put(12.70,0.00){$E_i$}

\end{picture}
\end{center}
}}

To that resolution we assign the vector space $V_c=V^{\otimes(p-1)}$, where we have assigned the f\mbox{}irst $i-1$ and the last $p-1-i$ copies of $V$ to the circles that are the closures of the f\mbox{}irst $i-1$ and the last $p-i-1$ strands of the resolution $D_c$, respectively, and the $i$-th copy of $V$ corresponds to the remaining circle (closure of $E_i$). So, we have that $C^1(D)=\bigoplus_{c\in c(D)}{V_c \{1\}}$. Further on, we denote the $k$-th copy of $V$ in $V_c$ by $V^k_c$.\\ 

\indent The dif\mbox{}ferential $d^0: C^0(D)\to C^1(D)$ is given by the maps $f_{c}: V^{\otimes p}\to V_{c}$ where if $c$ is of the type $\sigma_i$ than $f_c$ acts as the identity on the f\mbox{}irst $i-1$ and the last $p-i-1$ copies of $V$ and as the multiplication $m$ on the remaining two copies of $V$. In other words, $f_c$ maps the copies $V^j$ as the identity onto $V^j_c$, for $j<i$, maps the copies $V^{j+1}$ as the identity onto $V^{j}_c$, for $i<j<p$, and on the remaining two factors $f_c$ acts as the  multiplication $m:V^i\otimes V^{i+1} \to V_c^i$.\\

\indent Furthermore, $C^2(D)$ comes from the resolutions where exactly two of the crossings are resolved into 1-resolutions and the remaining ones are resolved into 0-resolutions. Denote the two  crossings that are resolved into   1-resolutions  by $c$ and $d$ (where $c$ is ordered before $d$), and let $c$ be of the type $\sigma_i$ and $d$ of the type $\sigma_j$. Then we have that $i \le j$.\\
\indent If $i=j$ then the corresponding resolution $D_{c,d}$ is the  closure of the plat diagram $E_i^2$, has $p$ circles and hence the corresponding summand  $V_{c,d}$ of  $C^2(D)$ is isomorphic to $V^{\otimes p}$. Here we have assigned the f\mbox{}irst $i-1$ and the last $p-i-1$ copies of $V$ to the circles that are the closures of the f\mbox{}irst $i-1$ and the last $p-i-1$ strands of the resolution $D_{c,d}$, respectively, while the  $i$-th and $(i+1)$-th copy of $V$ correspond to the remaining two circles that are formed from the  $i$-th and $(i+1)$-th strand (the $i$-th copy of $V$ corresponds to the outer, and the $(i+1)$-th copy of $V$ to the inner circle). \\
\indent If $i<j$ then the corresponding resolution is the closure of the plat diagram $E_i E_j$ (or $E_j E_i$), has $p-2$ circles and hence the corresponding summand $V_{c,d}$ of $C^2(D)$ is isomorphic to $V^{\otimes (p-2)}$. \\
\indent If $i+1<j$, we assign the f\mbox{}irst $i-1$ copies of $V$ to the closures of the f\mbox{}irst $i-1$ strands. The $i$-th copy of $V$ is assigned to the circle that is obtained by joining the $i$-th and $(i+1)$-th strand (the 1-resolution of $c$). The following $j-i-2$ copies of $V$ are assigned to the closures of the strands from the $(i+2)$-th to the $(j-1)$-th of $D_{c,d}$, respectively.
The $(j-1)$-th copy of $V$ is assigned to the circle that is obtained by joining the $j$-th and $(j+1)$-th strand (the 1-resolution of $d$). 
The remaining $p-j-1$ copies of $V$ are assigned to the closures of the last $p-j-1$ strands of $D_{c,d}$. \\
\indent If $i+1=j$, then we assign the f\mbox{}irst $i-1$ copies of $V$ to the closures of the f\mbox{}irst $i-1$ strands. The $i$-th copy of $V$ is assigned to the circle that is obtained by joining the $i$-th, $(i+1)$-th and $(i+2)$-th strand (the 1-resolutions of $c$ and $d$). The remaining $p-i-2$ copies of $V$ are assigned to the closures of the strands from the $(i+3)$-th to the $p$-th of $D_{c,d}$, respectively.\\
\indent In all previous cases, we denote the $k$-th copy of $V$ in $V_{c,d}$ by $V^k_{c,d}$.\\
 
\indent Finally, the second chain group is $$C^2(D)=\bigoplus_{c,d\in c(D),\,c<d} V_{c,d}\{2\}.$$

\indent Now, we can describe the dif\mbox{}ferential $d^1: C^1(D)\to C^2(D)$. 
It is given as the sum of the maps of the form $f_{ecd}: V_e\to V_{c,d}$ (with $c,d,e\in c(D),\,c<d$), where $f_{ecd}$ is zero unless $e=c$ or $e=d$. Let  $c=(i\alpha)$ and $d=(j\beta)$, for some $1\le i \le j \le p-1$, $\alpha=1,\ldots,l_i$ and $\beta=1,\ldots,l_j$. The maps $f_{ccd}: V_c\to V_{c,d}$ are as follows:  \\
\indent If $i=j$ then $f_{ccd}$ is given by the identity maps $id:V^l_c\to V^l_{c,d}$, for $l<i$, and $id:V^{l}_c\to V^{l+1}_{c,d}$, for $i<l<p$, and by the comultiplication map $\Delta:V^i_c\to V^i_{c,d}\otimes V^{i+1}_{c,d}$.\\
\indent If $i<j$, then $f_{ccd}$ is given by the identity maps
$id:V^l_c\to V^l_{c,d}$, for $l<j-1$, and $id:V^{l}_c\to V^{l-1}_{c,d}$, for $j<l\le p-1$, and by the multiplication map $m:V^{j-1}_c\otimes V^j_c \to V^{j-1}_{c,d}$.\\
\indent The other class of nonzero maps $f_{dcd}: V_d\to V_{c,d}$ is in the case $i=j$ given by $-f_{ccd}$ (because of the signs $(-1)^{f_{\nu}}$ sprinkled around the cubic complex and the fact that in this case $c<d$ and $V_c=V_d$). In the case $i<j$, $f_{dcd}$ is given as $-g_{cd}$ (the minus sign is because of the same reasons as above), where $g_{cd}$ is given by the identity maps:
$id:V^l_d\to V^l_{c,d}$, for $l<i$, and $id:V^{l+1}_d\to V^{l}_{c,d}$, for $i<l<p-1$, and by the multiplication map $m:V^{i}_d\otimes V^{i+1}_d \to V^{i}_{c,d}$.\\ 

\indent Now, we can go back to the proof. Let $$t'=(t_{1,1},\ldots,t_{1,l_1},t_{2,1},\ldots,t_{2,l_2},\ldots,t_{p-1,1},\ldots,t_{p-1,l_{p-1}})\in C^1(D)$$ be such that $d^1(t')=0$. Here we have that  $t_{i,\alpha}\in V_{(i\alpha)}$ for $i=1,\ldots,p-1$, $\alpha=1,\ldots,l_i$. We shall also write $t_{i\alpha}$ for $t_{i,\alpha}$ if there is no possibility of confusion. Our aim is to f\mbox{}ind an element $t\in C^0=V^{\otimes p}$ such that $d^0(t)=t'$. \\
\indent Since $d^1(t')=0$ we have that its projection, denoted by $d^1_{i\alpha\beta}(t')$, to the space $V_{(i\alpha),(i\beta)}$  is equal to zero, for every $i=1,\ldots,p-1$, and $\alpha,\beta=1,\ldots, l_i$. 
However, this implies that $t_{i\alpha}=t_{i\beta}$ for every
$i=1,\ldots,p-1$, $\alpha,\beta=1,\ldots, l_i$, since only the maps from $V_{(i\alpha)}$ and $V_{(i\beta)}$ to $V_{(i\alpha),(i\beta)}$ are nonzero, and the only nonidentity part of the mappings is the comultiplication $\Delta$ on the same ($i$-th) copy of $V$ in both $V_{(i\alpha)}$ and $V_{(i\beta)}$. \\
\indent 
Hence, we have obtained that $t'\in\ker d^1$ if and only if     
$t_{i\alpha}=t_{i\beta}$ for every
$i=1,\ldots,p-1$, $\alpha,\beta=1,\ldots, l_i$ and $\bar{t}=(t_{1,1},t_{2,1},\ldots,t_{p-1,1})\in\ker\bar{d}^1$, where $\bar{d}^1$ is the restriction of $d^1$ to $W_1=V_{(1,1)}\oplus V_{(2,1)}\oplus\cdots V_{(p-1,1)}$ (we omit $V_{(i,1)}$ if there are no crossings of the type $\sigma_i$ in $D$), composed with the projection onto $W_2=\bigoplus_{1\le i<j\le p-1}{V_{(i1),(j1)}}$ (omitting the indices $i,j$ of crossings which do not appear in $D$).\\
\indent Furthermore, note that for every $i=1,\ldots,p-1$, and $\alpha,\beta=1,\ldots, l_i$, the projection of $d^0(t)$, for any $t\in C^0(D)$ to $V_{(i\alpha)}$ and $V_{(i\beta)}$ is equal, i.e. the dif\mbox{}ferential $d^0$ is completely determined by the map $\bar{d}^0$ which is the projection of $d^0$ on $W_1$. Finally, we have that 
if there exists $t\in C^0$ such that 
$\bar{d}^0(t)=\bar{t}$, 
then $d^0(t)=t'$. Hence, to f\mbox{}inish the proof, we are left with proving that for every $y\in \ker \bar{d}^1 \subset W_1$, there exists $x \in C^0(D)$, such that $\bar{d}^0(x)=y$.\\
  
\indent 
Now, observe the positive braid  knot $K'$ which has the regular diagram $D'$ which is the closure of the braid $\sigma_1\sigma_2\cdots\sigma_{p-1}$ (where we omit the symbols $\sigma_i$ which are not contained in the braid word of which $D$ is the closure).  Its zeroth chain group $C^0(D')$ is obviously equal to $C^0(D)=V^{\otimes p}$, its f\mbox{}irst chain group $C^1(D')$ is equal to $W_1$ and its second chain group $C^2(D')$ is equal to $W_2$. Its zeroth dif\mbox{}ferential is equal to the previously def\mbox{}ined $\bar{d}^0$, while its f\mbox{}irst dif\mbox{}ferential is equal to $\bar{d}^1$. Since $K'$ is isotopic to the unknot (or to an  unlink consisting of unknots), its f\mbox{}irst homology group is trivial and hence for every $y\in\ker \bar{d}^1$ there exists $x\in C^0(D')=C^0(D)$ such that $\bar{d}^0(x)=y$. This concludes the proof. \kraj \\

\section{The $sl(n)$ case}\label{sln}

In this section we will adapt our proof from the previous section to show that the f\mbox{}irst $sl(n)$ homology group of a positive braid knot is trivial. As we saw in Section \ref{kovroznot}, the basic principle of the construction, namely the cubic complex, 0- and 1-resolutions of the crossings, per-edge maps and summing of the columns are the same as in the $sl(2)$ case (Section \ref{kovnot} and the previous section).
The slight dif\mbox{}ferences in the $sl(n)$ case for positive braid knots are: f\mbox{}irst of all, the 1-resolution of the crossing $c=(i\alpha)$ is given by the thick edge resolution $\bar{E}_i$:
\begin{center}
\setlength{\unitlength}{4mm}
\begin{picture}(30,6) 
\linethickness{0.6pt}
%\line(1.00,2.00,1.00,3.00)
\qbezier(5.00,1.00)(5.00,3.00)(5.00,5.00)
\qbezier(4.80,4.60)(4.90,4.80)(5.00,5.00)
\qbezier(5.20,4.60)(5.10,4.80)(5.00,5.00)

\qbezier(9.00,1.00)(9.00,3.00)(9.00,5.00)
\qbezier(8.80,4.60)(8.90,4.80)(9.00,5.00)
\qbezier(9.20,4.60)(9.10,4.80)(9.00,5.00)

\linethickness{2.0pt}
\qbezier(13.00,2.00)(13.00,3.00)(13.00,4.00)
\linethickness{0.6pt}
\qbezier(12.00,1.00)(12.50,1.50)(13.00,2.00)
\qbezier(14.00,1.00)(13.50,1.50)(13.00,2.00)
\qbezier(12.00,5.00)(12.50,4.50)(13.00,4.00)
\qbezier(14.00,5.00)(13.50,4.50)(13.00,4.00)

\qbezier(14.00,5.00)(13.80,4.90)(13.60,4.80)
\qbezier(14.00,5.00)(13.90,4.80)(13.80,4.60)

\qbezier(12.00,5.00)(12.20,4.90)(12.40,4.80)
\qbezier(12.00,5.00)(12.10,4.80)(12.20,4.60)

\qbezier(12.50,1.50)(12.30,1.40)(12.10,1.30)
\qbezier(12.50,1.50)(12.40,1.30)(12.30,1.10)

\qbezier(13.50,1.50)(13.70,1.40)(13.90,1.30)
\qbezier(13.50,1.50)(13.60,1.30)(13.70,1.10)

\qbezier(17.00,1.00)(17.00,3.00)(17.00,5.00)
\qbezier(16.80,4.60)(16.90,4.80)(17.00,5.00)
\qbezier(17.20,4.60)(17.10,4.80)(17.00,5.00)

\qbezier(21.00,1.00)(21.00,3.00)(21.00,5.00)
\qbezier(20.80,4.60)(20.90,4.80)(21.00,5.00)
\qbezier(21.20,4.60)(21.10,4.80)(21.00,5.00)

\put(6.00,3.00){$\cdot$}
\put(7.00,3.00){$\cdot$}
\put(8.00,3.00){$\cdot$}
\put(18.00,3.00){$\cdot$}
\put(19.00,3.00){$\cdot$}
\put(20.00,3.00){$\cdot$}

\put(4.90,5.50){1}
\put(11.85,5.50){$i$}
\put(13.40,5.50){$i+1$}
\put(20.90,5.50){$p$}
\put(12.50,-0.3){$\bar{E}_i$} 
\end{picture}
\end{center}
instead of the plat $E_i$. Hence the total resolutions, $D_{\epsilon}$, are in this case trivalent graphs with thick edges. 
We denote the total resolution when all crossings of $D$ are resolved into 0-resolution by $D_0$,  the specif\mbox{}ic total resolution with the crossing $c$ resolved into a 1-resolution and all other crossings into 0-resolutions by $D_c$, and the specif\mbox{}ic total resolution with the crossings $c$ and $d$ resolved into 1-resolutions and all other crossings into 0-resolutions  by $D_{c,d}$ (like in the $sl(2)$ case). Furthermore, we assign to the diagram $D_{\epsilon}$ the $\Q$-vector space $\bar{V}_{\epsilon}$ (we will also write $\bar{V}(D_{\epsilon})$), and we def\mbox{}ine the $i$-th chain group by
$$C_n^i(D)=\oplus_{|\epsilon|=i}\bar{V}_{\epsilon}\{-i\}.$$
Also, we def\mbox{}ine the dif\mbox{}ferentials as the sum (we have already included the same signs as in the $sl(2)$ case) of the grading preserving per-edge maps $\bar{d}_{\nu}:\bar{V}_{\epsilon} \to \bar{V}_{\epsilon'}\{-1\}$. 
Again, we denote $\bar{d}_{\nu}:\bar{V}_0 \to \bar{V}_c\{-1\}$ by $f_c$, and $\bar{d}_{\nu}:\bar{V}_e \to \bar{V}_{c,d}\{-1\}$ by $f_{ecd}$. Since we sum only per-edge maps, $f_{ecd}$ is nonzero only if $e=c$ or $e=d$.\\
%Here, as 
%before, $\epsilon'$ has the same $n-1$ entries as $\epsilon$ and the %remaining one, say $j$-th one, is equal to 0 in $\epsilon$ while %it is 
%equal to 1 in $\epsilon'$, and $\nu$ is the $n$-tuple that 
%coincide with 
%$\epsilon$ at all entries except $j$-th one, where it is %equal to 
%$\ast$.\\
\indent The main dif\mbox{}ference in the def\mbox{}inition of $sl(n)$ homology compared to the $sl(2)$ case, is the very complicated assignment of the vector spaces $\bar{V}_{\epsilon}$ as well as the maps $\bar{d}_{\nu}$, which make the theory practically incalculable, and very dif\mbox{}f\mbox{}icult to deal with. However, in the case that we are interested in, we only need the chain groups $C_n^0(D)$, $C_n^1(D)$ and $C_n^2(D)$, and the dif\mbox{}ferentials $d_n^0$ and $d_n^1$, and their def\mbox{}initions can be extracted explicitly from the def\mbox{}inition (see Section \ref{kovroznot2} and also \cite{kovroz}).\\
\indent First of all, if we have two disjoint trivalent graphs with wide edges, $\Gamma_1$ and $\Gamma_2$, then the vector space assigned to their disjoint union, $\bar{V}(\Gamma_1 \coprod \Gamma_2)$, is isomorphic to the tensor product of the spaces $\bar{V}(\Gamma_1)$ and $\bar{V}(\Gamma_2)$. Next, if we have a braid diagram with $p$ strands, we assign the variable $x_i$ to the $i$-th strand. We introduce a grading on the ring of polynomials $R=\Q[x_1,\ldots,x_p]$, by def\mbox{}ining $\deg x_i=2$. To the circle (unknot) labeled by the variable $x_i$, we assign the vector space $A_i$ given by:
$$A_i=\Q[x_i]/(x_i^n)\{1-n\}.$$
Note that the quantum dimension of $A_i$ is equal to $[n]$. Another type of diagrams that we will encounter is:
%will be the following one:
\begin{center}
\setlength{\unitlength}{4mm}
\begin{picture}(20,3.5) 
\linethickness{0.6pt}
\put(9,2){\oval(2,4)}
\put(11,2){\oval(2,4)}
\put(8,2){\qbezier(0,0)(-0.13,0.13)(-0.26,0.26)\qbezier(0,0)(0.13,0.13)(0.26,0.26)}
\put(7,3){$x_i$}
\put(12,2){\qbezier(0,0)(-0.13,0.13)(-0.26,0.26)\qbezier(0,0)(0.13,0.13)(0.26,0.26)}
\put(12.5,3){$x_{i+1}$}
\linethickness{2.5pt}
\put(10,0.85){\qbezier(0,0)(0,1.2)(0,2.4)}
\end{picture}
\end{center}
i.e. the closure of $\bar{E}_i$ (below, by $\bar{E}_i$ we mean only the part of the whole braid formed by the strands that are connected to a wide edge). To this graph we assign the vector space $B_{i}$ given by:
$$B_{i}=\Q[x_i,x_{i+1}]/(\pi_{i,i+1}^n,\pi_{i,i+1}^{n-1})\{3-2n\}.$$
Here, by $\pi_{i,j}^k$ we denote the following polynomial:
$$\pi_{i,j}^k=\frac{x_i^{k+1}-x_j^{k+1}}{x_i-x_j}
=\sum_{l=0}^k {x_i^l x_j^{k-l}}.$$
Next, to the closure of the diagrams $\bar{E}_i^2$

\begin{center}
\setlength{\unitlength}{4mm}
\begin{picture}(20,4) 
\linethickness{0.6pt}

\put(10,0){
\linethickness{2.5pt}
\qbezier(0,0.5)(0,0.75)(0,1)
\qbezier(0,3)(0,3.25)(0,3.5)
\linethickness{0.6pt}
\put(0,0.1){\qbezier(-0.5,2)(-0.6,1.85)(-0.7,1.7)
\qbezier(-0.5,2)(-0.4,1.85)(-0.3,1.7)}
\put(0,0.1){\qbezier(0.5,2)(0.6,1.85)(0.7,1.7)
\qbezier(0.5,2)(0.4,1.85)(0.3,1.7)}

\qbezier(0,1)(-0.5,1.6)(-0.5,2)
\qbezier(0,3)(-0.5,2.4)(-0.5,2)
\qbezier(0,1)(0.5,1.6)(0.5,2)
\qbezier(0,3)(0.5,2.4)(0.5,2)

\qbezier(0,0.5)(-0.5,-0.1)(-0.5,-0.5)
\qbezier(0,0.5)(0.5,-0.1)(0.5,-0.5)
\qbezier(0,3.5)(-0.5,4.1)(-0.5,4.5)
\qbezier(0,3.5)(0.5,4.1)(0.5,4.5)

\qbezier(-2,-0.5)(-2,2)(-2,4.5)
\qbezier(2,-0.5)(2,2)(2,4.5)

\qbezier(-2,4.5)(-1.25,5.5)(-0.5,4.5)
\qbezier(-2,-0.5)(-1.25,-1.5)(-0.5,-0.5)
\qbezier(2,4.5)(1.25,5.5)(0.5,4.5)
\qbezier(2,-0.5)(1.25,-1.5)(0.5,-0.5)

\put(-2,2){\qbezier(0,0)(-0.1,0.15)(-0.2,0.3)
\qbezier(0,0)(0.1,0.15)(0.2,0.3)
}
\put(2,2){\qbezier(0,0)(-0.1,0.15)(-0.2,0.3)
\qbezier(0,0)(0.1,0.15)(0.2,0.3)
}
\put(-3,3.5){$x_i$}
\put(2.2,3.5){$x_{i+1}$}
\put(-1.1,2.2){$y$}
\put(0.8,2.2){$z$}
}
\end{picture}
\end{center}
%i.e. the %graphs:
%
%\begin{center}
%\setlength{\unitlength}{4mm}
%\begin{picture}(20,3) 
%\linethickness{0.6pt}
%\put(9,2){\oval(2,4)}
%\put(11,2){\oval(2,4)}
%\put(8,2){\qbezier(0,0)(-0.13,0.13)(-0.26,0.26)\qbezier(0,0)(0.13,0.13)(%0%.26,0.26)}
%\put(7,3){$x_i$}
%\put(12,2){\qbezier(0,0)(-0.13,0.13)(-0.26,0.26)\qbezier(0,0)(0.13,0.13)%(%0.26,0.26)}
%\put(12.5,3){$x_j$}
%\linethickness{2.5pt}
%\put(10,0.85){\qbezier(0,0)(0,1.2)(0,2.4)}
%\end{picture}
%\end{center}
%
%To this graph 
we assign the graded vector space $C_{i}$, def\mbox{}ined by:
$$C_{i}=\Q[x_i,x_{i+1},y]/(\pi_{i,i+1}^n,\pi_{i,i+1}^{n-1},(y-x_i)(y-x_{i+1}))\{2-2n\},$$
where $y$ is a new variable with $\deg y=2$.\\
\indent To the graph which is the closure of the diagram $\bar{E}_i\bar{E}_j$, for $|i-j|>1$ we assign the vector space $F_{i,j}$ def\mbox{}ined by 
$$F_{i,j}=\Q[x_i,x_{i+1},x_j,x_{j+1}]/(\pi_{i,i+1}^n,\pi_{i,i+1}^{n-1},\pi_{j,j+1}^n,\pi_{j,j+1}^{n-1})\{6-4n\}.$$
Finally, to the closure of $\bar{E}_i\bar{E}_{i+1}$ we assign the vector space $F_{i,i+1}$ def\mbox{}ined by
$$F_{i,i+1}=\Q[x_i,x_{i+1},x_{i+2}]/(\pi_{i,i+1}^n,\pi_{i,i+1}^{n-1},\pi_{i+1,i+2}^n,\pi_{i+1,i+2}^{n-1})\{5-3n\}.$$
\indent Now we can give the main theorem:

\begin{theorem}\label{slngl}
For every positive braid knot $K$, we have that the homology group $\H^1_n(K)$ is trivial.
\end{theorem}
\textbf{Proof:}\\
\indent Mainly we will adapt our proof from the previous Section for this case. First of all, we again take a diagram $D$ of $K$, which is the closure of the positive braid $\sigma$. Denote the number of strands  of $\sigma$ by $p$, and assign the variable $x_i$, with $\deg x_i=2$, to the $i$-th strand. Since $D$ has only positive crossings we have that $\H_n^1(K)$ is trivial if and only if $H_n^1(D)$ is trivial. So, we have to prove that the latter group is trivial. \\
\indent Again, we use the direct sum def\mbox{}inition of  the chain groups $C_n^0(D)$, $C_n^1(D)$, $C_n^2(D)$, and of the dif\mbox{}ferentials $d_n^0$ and $d_n^1$.  We will prove that for each $t' \in C_n^1(D)=\bigoplus_{c}{\bar{V}_c}\{-1\}$ 
such that $d_n^1(t')=0$, there exists $t \in C_n^0(D)=\bar{V}_0=\bigotimes_{i=1}^p {A_i}$ such that $d_n^0(t)=t'$.\\
\indent Let $i\in \{1,\ldots,p-1\}$ and $1\le \alpha < \beta \le l_i$ be arbitrary. Denote $c=(i\alpha)$ and $d=(i\beta)$. We denote the restriction of $t'$ to the space $\bar{V}_c=\bar{V}_{(i\alpha)}$ by $t_{i\alpha}$. Again the dif\mbox{}ferential ${d}_n^1$ maps only the spaces $\bar{V}_{(i\alpha)}$ and $\bar{V}_{(i\beta)}$ into  $\bar{V}_{(i\alpha), (i\beta)}\{-1\}$. Let us see what these spaces and maps are.\\
\indent The space $\bar{V}_{c}=\bar{V}_{(i\alpha)}$ is assigned to the diagram $D_{c}$ which is the disjoint union of $p-2$ circles and the closure of the diagram $\bar{E}_i$. Hence, we have:
$$\bar{V}_{c}=A_1\otimes \cdots \otimes A_{i-1} \otimes B_i \otimes A_{i+2} \otimes \cdots \otimes A_p,$$
and we have completely the same expression for $\bar{V}_{d}=\bar{V}_{(i\beta)}$. The space $\bar{V}_{c,d}=\bar{V}_{(i\alpha),(i\beta)}$ is assigned to the diagram $D_{c,d}$ which is the disjoint union of $p-2$ circles and the closure of the diagram $\bar{E}_i^2$, and so:
$$\bar{V}_{(i\alpha),(i\beta)}=A_1\otimes \cdots \otimes A_{i-1} \otimes C_i \otimes A_{i+2} \otimes \cdots \otimes A_p.$$
Finally, the map $f_{ccd}:\bar{V}_c \to \bar{V}_{c,d}$ is the identity on the factors $A_j$, for $1\le j <i$ and $i+1<j \le p$, and on the remaining part is given as the multiplication by $y-x_{i+1}$.
%$$m_{y-x_{i+1}}:B_{i,i+1} \to C_{i,i+1}\{-1\},$$
The map $f_{dcd}:\bar{V}_d \to \bar{V}_{c,d}$ is just the negative of $f_{ccd}$ (i.e. $f_{ccd}=-f_{dcd}$), because of the signs sprinkled around the cube of resolutions that make it anticommutative.\\
\indent Since $d_n^1(t')=0$ we have that the projection of $d_n^1(t')$ to $\bar{V}_{c,d}$ is zero. However, the only nonzero maps that have $\bar{V}_{c,d}$ as the codomain are $f_{ccd}$ and $f_{dcd}$, and so we have that 
$t_{i\alpha}=t_{i\beta}$, like in the $sl(2)$ case.\\
\indent Furthermore, if $c=(i\alpha)$ and $d=(j\beta)$, for $i<j$, then the vector spaces $\bar{V}_{c,d}$ do not depend on $\alpha$ and $\beta$ (they are the tensor product of $F_{i,j}$ and $A_k$'s), and the maps $f_{ccd}$ and $f_{dcd}$ are, up to a sign, identities (multiplications by 1). Also, for all crossings $c=(i\alpha)$, the spaces $\bar{V}_c$ do not depend on $\alpha$ (they are the tensor product of $B_i$ and $A_k$'s), and the maps $f_c$ are identities.\\
\indent Hence, like in the $sl(2)$ case, we have that $d_n^1(t')=0$ if and only if $t_{i\alpha}=t_{i\beta}$ for every 
$i=1,\ldots,p-1$, $\alpha,\beta=1,\ldots, l_i$ and $\bar{t}=(t_{1,1},t_{2,1},\ldots,t_{p-1,1})\in\ker\bar{d}^1$, where $\bar{d}^1$ is the restriction of $d^1$ to $\bar{W}_1=\bar{V}_{(1,1)}\oplus \bar{V}_{(2,1)}\oplus\cdots \bar{V}_{(p-1,1)}$ (here we have omitted $\bar{V}_{(i,1)}$ if there are no crossings of type $\sigma_i$ in $D$) composed with the projection onto $\bar{W}_2=\bigoplus_{1\le i<j\le p-1}{\bar{V}_{(i1),(j1)}}$ (omitting the indices $i,j$ of crossings which do not appear in $D$). Also,  for every $i=1,\ldots,p-1$, $\alpha,\beta=1,\ldots, l_i$ and $t\in C_n^0(D)$, the projection of $d_n^0(t)$  to $\bar{V}_{(i\alpha)}$ and $\bar{V}_{(i\beta)}$ is equal, i.e. the dif\mbox{}ferential $d_n^0$ is completely determined by the map $\bar{d}_n^0$ which is the projection of $d_n^0$ onto $\bar{W}_1$. So, we have that 
if there exists $t\in C_n^0(D)$ such that
$\bar{d}_n^0(t)=\bar{t}$, 
then $d_n^0(t)=t'$. Hence, to f\mbox{}inish the proof, we are left with proving that for every $y\in \ker \bar{d}_n^1 \subset \bar{W}_1$, there exists $x \in C_n^0(D)$, such that $\bar{d}_n^0(x)=y$.\\
\indent However, the last statement again follows from observing the 
positive braid  knot $K'$ which has the regular diagram $D'$ which is the closure of the braid $\sigma_1\sigma_2\cdots\sigma_{p-1}$ (where we omit the symbols $\sigma_i$ which are not contained in the braid word of which $D$ is the closure).  Its zeroth chain group $C_n^0(D')$ is obviously equal to $C_n^0(D)$, its f\mbox{}irst chain group $C_n^1(D')$ is equal to $\bar{W}_1$ and its second chain group $C^2(D')$ is equal to $\bar{W}_2$. Its zeroth dif\mbox{}ferential is equal to the previously def\mbox{}ined $\bar{d}_n^0$, while its f\mbox{}irst dif\mbox{}ferential is equal to $\bar{d}_n^1$. Since $K'$ is isotopic to the unknot (or to an  unlink consisting of unknots), its f\mbox{}irst homology group is trivial (from Theorem \ref{krgl}) and hence for every $y\in\ker \bar{d}_n^1$ there exists $x\in C_n^0(D')=C_n^0(D)$ such that $\bar{d}_n^0(x)=y$. This concludes the proof. \kraj \\

\chapter{Thickness and stability of link homology for torus knots}
\label{torus}

\section{Introduction}

In this chapter we f\mbox{}irst show that the torus knots $T_{p,q}$ for $3\le p\le q$ (non-alternating torus knots) are homologically thick, i.e. that their Khovanov homology occupies at least three diagonals. Furthermore, in the course of the proof we obtain even stronger results that relate the homology of the torus knots $T_{p,q}$ and $T_{p,q+1}$. Namely, we prove that, up to a certain homological degree, their (unnormalized) homologies coincide. \\
\indent As the f\mbox{}irst application of this result we calculate the homology of torus knots for low homological degrees. We also obtain the proof of the existence of stable Khovanov homology of torus knots, conjectured in \cite{dgr}.\\
\indent Furthermore, we conjecture that the homological width of the torus knot $T_{p,q}$ is at least $p$, and we reduce this problem to determining the nontriviality of certain homological groups.\\
\indent Since in the proofs we mainly use the long exact sequence of Khovanov homology (\ref{eq2}) -- we do not rely heavily on the explicit def\mbox{}inition of $sl(2)$ homology -- we also obtain most of the analogous results for the stability of $sl(n)$ homology  of torus knots.

\section{Torus knots}

\indent A knot or a link is a torus knot if it is isotopic to  a knot or a link that can be drawn without any points of intersection on the trivial torus. Every torus link is, up to a mirror image,  determined by two nonnegative integers $p$ and $q$, i.e. it is isotopic to a unique torus knot $T_{p,q}$ which has the diagram $D_{p,q}$ - the closure of the braid $(\sigma_1\sigma_2\ldots\sigma_{p-1})^q$ - as a planar projection. In other words, $D_{p,q}$ is the closure of the $p$-strand braid with $q$ full twists. Since $T_{p,q}$ is isotopic to $T_{q,p}$ we can assume that $p\le q$.\\
\indent The torus knots are obviously positive braid knots (since the diagram $D_{p,q}$ is the closure of a positive braid). Hence, their homology satisf\mbox{}ies the properties from the previous chapter. Also, we will use the same notation and conventions from  Section \ref{posbraid}. 
  
\section{Thickness of torus knots}\label{torglavni}

  If $p=1$ then the torus knot $T_{p,q}$ is trivial and for $p=2$ the torus knot $T_{2,q}$ is alternating, hence its homology occupies exactly two diagonals. However, if $p\ge 3$, the torus knot $T_{p,q}$ is non-alternating and we will prove that its homology occupies at least three diagonals. Namely, we prove the following theorem:
\begin{theorem}\label{T1}
Let $K=T_{p,q}$, $3\le p\le q$ be a torus knot. Then 
$$\rank{\H^{4,(p-1)(q-1)+5}(K)}> 0.$$
\end{theorem}

From this theorem we obtain 
\begin{corollary}\label{C1}
Every torus knot $T_{p,q}$, $p,q\ge 3$ is H-thick, i.e. its Khovanov homology occupies at least three diagonals.
\end{corollary}
\textbf{Proof} (of Corollary \ref{C1}):\\

\indent Since $T_{p,q}$ is a positive knot, its zeroth homology group is two dimensional and the $q$-gradings (and consequently the $\delta$-gradings) of its two generators are $(p-1)(q-1)-1$ and $(p-1)(q-1)+1$, respectively (see e.g. Section \ref{posbraid}). However, from Theorem \ref{T1} we have that there exists a generator with $t$-grading equal to 4 and $q$-grading equal to $(p-1)(q-1)+5$, and so its $\delta$-grading is equal to $(p-1)(q-1)+5-2\cdot 4=(p-1)(q-1)-3$. Thus, we have obtained three generators    
of the homology of the torus knot $T_{p,q}$ which have three dif\mbox{}ferent values of the $\delta$-grading and hence its Khovanov homology occupies at least three diagonals. \kraj\\

\indent Now we give a proof of Theorem \ref{T1}.\\
\\
\textbf{Proof}:\\
\indent 
First of all, since $K=T_{p,q}$ is a positive braid knot whose regular diagram $D_{p,q}$ is the closure of the braid $(\sigma_1\sigma_2\ldots\sigma_{p-1})^q$
with $(p-1)q$ crossings, we have that $\H^{4,(p-1)(q-1)+5}(K)=H^{4,6-p}(D_{p,q})$. So, we will ``concentrate" on calculating the latter homology group, i.e. showing that its rank is nonzero. In order to do this we will use the long exact sequence (\ref{eq2}) and we will relate the unnormalized fourth homology groups of the standard regular diagrams of the torus knots $D_{p,q}$ and $D_{p,q-1}$ for $p<q$.\\

\indent Let $3\le p <q$. Let $c_{p-1}$ be the crossing $(p-1,1)$ of the diagram $D_{p,q}$. Now denote by $E_{p,q}^1$ and $D_{p,q}^1$ the 1- and 0-resolutions, respectively, of the diagram $D_{p,q}$ at the crossing $c_{p-1}$. Then from (\ref{eq2}) we obtain
{\small{
\begin{equation*}\cdots\rightarrow H^{3,j-1}(E^1_{p,q})\to H^{4,j}(D_{p,q})\to H^{4,j}(D_{p,q}^1)\to H^{4,j-1}(E^1_{p,q})\to H^{5,j}(D_{p,q})\to \cdots\end{equation*}
}}
\indent Now, we can continue the process, and resolve the crossing $c_{p-2}=(p-2,1)$ of $D_{p,q}^1$ in two possible ways. Denote the diagram obtained by the 1-resolution by $E^2_{p,q}$, and the diagram obtained by the 0-resolution by $D^2_{p,q}$. Then from the long exact sequence (\ref{eq2})
we have:
{\small{
\begin{equation*}\cdots\rightarrow H^{3,j-1}(E^2_{p,q})\to H^{4,j}(D^1_{p,q})\to H^{4,j}(D_{p,q}^2)\to H^{4,j-1}(E^2_{p,q})\to H^{5,j}(D^1_{p,q})\to \cdots\end{equation*}
}}
\indent After repeating this process $p-1$ times (resolving the crossing $c_{p-k}=(p-k,1),\,k=1,\ldots,p-1$, of $D^{k-1}_{p,q}$, obtaining the 1-resolution $E^k_{p,q}$ and 0-resolution $D^k_{p,q}$ and applying the same long exact sequence in homology), we obtain that for every $i=1,\ldots,p-1$:
{\small{
\begin{equation}\cdots\rightarrow H^{3,j-1}(E^i_{p,q})\to H^{4,j}(D^{i-1}_{p,q})\to H^{4,j}(D_{p,q}^i)\to H^{4,j-1}(E^i_{p,q})\to H^{5,j}(D^{i-1}_{p,q})\to \cdots\label{eq3}\end{equation}
}}
Here  $D^0_{p,q}$ denotes $D_{p,q}$, and we obviously have that $D^{p-1}_{p,q}=D_{p,q-1}$.\\
\indent Our goal is to show $H^3(E^i_{p,q})$ and $H^4(E^i_{p,q})$ are trivial for every $0<i<p$. This is done in the following lemma.

\begin{lemma} \label{lem1}
For every three positive integers $p$, $q$ and $i$, such that $3 \le p < q$ and $i<p$, the knot with the diagram $E^i_{p,q}$ is positive, and the diagram $E^i_{p,q}$ has at least $p+q-3$ negative crossings.
\end{lemma}
\textbf{Proof:}\\

\indent Since for every $0<i<p$, $E^{i}_{p,q}$ is obtained by the 1-resolution of the crossing $c_{p-i}=(p-i,1)$ of the (positive braid knot) diagram $D^{i-1}_{p,q}$, it is the closure of the plat braid diagram with only one plat $E_{p-i}$:

{\begin{center}\setlength{\unitlength}{4mm}
\begin{picture}(30,6)
\linethickness{0.8pt}
\qbezier(5.00,1.00)(5.00,3.00)(5.00,5.00)
\qbezier(9.00,1.00)(9.00,3.00)(9.00,5.00)
\qbezier(15.00,5.00)(13.00,2.00)(11.00,5.00)
\qbezier(11.00,1.00)(13.00,4.00)(15.00,1.00)
\qbezier(17.00,1.00)(17.00,3.00)(17.00,5.00)
\qbezier(21.00,1.00)(21.00,3.00)(21.00,5.00)
\put(6.00,3.00){$\cdot$}
\put(7.00,3.00){$\cdot$}
\put(8.00,3.00){$\cdot$}
\put(18.00,3.00){$\cdot$}
\put(19.00,3.00){$\cdot$}
\put(20.00,3.00){$\cdot$}
\put(4.90,5.50){1}
\put(10.50,5.50){$p-i$}
\put(14.00,5.50){$p-i+1$}
\put(20.90,5.50){$p$}
\put(12.70,0.00){$E_{p-i}$}
\end{picture}
\end{center}
}

 Now, note that two lower strands of $E_{p-i}$ are always ``neighbour" strands, i.e. they form a ribbon, through the diagram, until they reach the upper part of $E_{p-i}$. So, we can``slide" the lower part of the plat $E_{p-i}$ through the diagram (by using the second Reidemeister move -- R2 -- and also the f\mbox{}irst Reidemeister move -- R1 -- where the two strands intersect each other) until it reaches the left or the right hand side of the upper part of $E_{p-i}$. If it f\mbox{}irst reaches the right hand side, it automatically (or after a R1 move) becomes the closure of a positive braid diagram. If it f\mbox{}irst reaches the left hand side then after a R1 move and a ``slide" (sequence of R2 moves) we obviously obtain a positive braid diagram.\\
\indent Concerning the number of negative crossings of $E^i_{p,q}$, 
note that by performing the f\mbox{}irst sequence of R2 moves (sliding the lower part of $E_{p-i}$ through the diagram, from the top to the bottom) in each move we have ``canceled" one positive and one negative crossing. Furthermore, since $p<q$, the two strands of the lower part of the plat $E_{p-i}$ will make a full twist at least once, and so they will have two crossings with each other, which are both obviously negative crossings. So, we have that on each of the last (lower) $q-1$ blocks ($\sigma_1\ldots\sigma_{p-1}$) of $E^i_{p,q}$ we have at least one negative crossing. Furthermore, on the part where both lower strands of the plat $E_{p-i}$ make full twists, we have applied $p-2$ R2 moves, and hence we have in addition, at least, $p-2$ negative crossings.\\
\indent  Altogether, this gives at least $q-1+p-2=p+q-3$ negative crossings of $E^i_{p,q}$, as required. \kraj\\

\indent Now, we can go back to the proof. From the previous Lemma, we conclude that for every $0<k<p$, $\H^{i,j}(E^k_{p,q})$ is trivial for $i<0$, and that \begin{equation}
H^{i,j}(E^k_{p,q}) \textrm{ is trivial for } i<p+q-3,
\label{fspust}
\end{equation}
Hence, if $p+q-3>4$, from 
(\ref{eq3}) and (\ref{fspust}), we obtain that
$$H^{4,j}(D^{i-1}_{p,q})=H^{4,j}(D^i_{p,q}),\quad i=1,\ldots,p-1,$$
and thus we have that 
\begin{equation}
H^{4,j}(D_{p,q})=H^{4,j}(D_{p,q-1}),\quad\textrm {for } p+q>7, \,\, p<q.
\label{spust}
\end{equation}
So, we can decrease the number of full twists, $q$, without changing the fourth homology group.\\
\indent If $p=q=3$, then e.g. by using programs for computing Khovanov homology (\cite{bnprog}, \cite{shum}), we obtain that $\rank H^{4,3}(D_{3,3})=\rank \H^{4,9}(T_{3,3})=1$, as wanted. \\
\indent If $3=p<q$, then from (\ref{spust}), we have that $H^{4,3}(D_{3,q})=H^{4,3}(D_{3,4})$. However, by using programs for computing Khovanov homology, we obtain that the rank of the latter group (which is equal to $\H^{4,11}(T_{3,4})$) is equal to 1. 
\begin{remark}
In Bar-Natan's tables of knots in \cite{bn}, the link $T_{3,3}$ is denoted by $6^3_3$, and the torus knot $T_{3,4}$ is isotopic to the knot $8_{19}$. For the general notation of knots and links see \cite{rol} and \cite{bnprog}.
\end{remark}
\indent Now, let us move to the general case $4\le p\le q$. Then if $p<q$ we can apply (\ref{spust}) and obtain $H^{4,6-p}(D_{p,q})=H^{4,6-p}(D_{p,p})$. Thus, we are left with proving that the latter homology group is of nonzero rank.\\
\indent Now, apply the set of long exact sequences ($\ref{eq3}$) for the case $p=q$. In this case, like in Lemma \ref{lem1}, we obtain that $E^k_{p,p}$ is the diagram of a positive knot, for every $k=1,\ldots,p-1$. Furthermore, every diagram $E^k_{p,p}$ has exactly $2p-3$ negative crossings, and so we have:
\begin{equation}
H^i(E^k_{p,p}) \textrm{ is trivial for every }i<2p-3.
\label{novspust}
\end{equation} Since $p>3$, we have that $H^3(E^k_{p,p})$ and 
$H^4(E^k_{p,p})$ are trivial for all $k$, and so we have that
\begin{equation}
H^4(D_{p,p})=H^4(D_{p,p-1}).
\label{sp1}
\end{equation}
\indent On the other hand, we have that 
$$H^{4,6-p}(D_{p,p-1})=\H^{4,p^2-3p+7}(T_{p,p-1})=\H^{4,p^2-3p+7}(T_{p-1,p})=H^{4,7-p}(D_{p-1,p}).$$
\indent If $p>4$, then we have from (\ref{spust}) that $H^{4,7-p}(D_{p-1,p})=H^{4,7-p}(D_{p-1,p-1})$. 
By repeating this process, we can decrease the number of strands $p$, and obtain that:
$$H^{4,6-p}(D_{p,p})=H^{4,2}(D_{4,4}).$$
Finally, from (\ref{sp1}) we have $$H^{4,2}(D_{4,4})=H^{4,2}(D_{4,3})=\H^{4,11}(T_{3,4}),$$
and the last homology group, as we saw previously, is of rank 1. This concludes our proof.   \kraj 

\section{Stability of Khovanov homology for torus knots}\label{sledeci}

In the course of proving Theorem \ref{T1}, apart from showing that the torus knots are H-thick, we have obtained some other properties of the homology of the torus knots. Namely, we proved that we can reduce the number of full twists, $q$, of the standard diagram $D_{p,q}$ of the torus knot $T_{p,q}$ without changing the f\mbox{}irst $p+q-3$ homology groups. In other words, we have obtained the existence of stable homology of torus knots (see Section \ref{stable} and \cite{dgr}).\\
\indent In the following theorem we summarize the stability properties obtained in the previous section.
\begin{theorem} \label{T3}
Let $p$, $q$ and $i$ be integers such that $2\le p<q$ and $i<p+q-3$. Then for every $j\in\Z$
\begin{equation}
H^{i,j}(D_{p,q})=H^{i,j}(D_{p,q-1}).
\label{f1}
\end{equation}
Furthermore, for every $2\le p<q$ and $i<2p-1$ and $j\in\Z$ we have 
\begin{equation}
H^{i,j}(D_{p,p+1})=H^{i,j}(D_{p,p+2})=\cdots=H^{i,j}(D_{p,q}).
\label{f2}
\end{equation}
Also, for every $p\ge 2$, $i<2p-3$ and $j\in\Z$, we have
\begin{equation}
H^{i,j}(D_{p,p})=H^{i,j+1}(D_{p-1,p}).
\label{f3}
\end{equation}
\end{theorem} 
\textbf{Proof:}\\
\indent The equations (\ref{f1}) and (\ref{f3}) we have already obtained in the course of proving Theorem \ref{T1} (long exact sequences and formulas (\ref{fspust}) and (\ref{novspust})). Formula (\ref{f2}) obviously follows from (\ref{f1}) since $i<2p-1=p+(p+2)-3$. \kraj\\

The torus knots $T_{2,q}$ are alternating and their homology is well-known (see e.g. \cite{kov}). However, as the f\mbox{}irst corollary of the previous theorem we obtain the homology groups of $T_{p,q}$, for $3\le p \le q$, with low homological degree. 

\begin{theorem}
Let $3\le p\le q$ with $p$ and $q$ not both equal to 3. Then we have 
\begin{eqnarray*}
\H^{0,(p-1)(q-1)\pm 1}(T_{p,q})&=&\Z\\
\H^{2,(p-1)(q-1)+3}(T_{p,q})&=&\Z\\
\H^{3,(p-1)(q-1)+7}(T_{p,q})&=&\Z\\
\H^{3,(p-1)(q-1)+5}(T_{p,q})&=&\Z_2\\
\H^{4,(p-1)(q-1)+6 \pm 1}(T_{p,q})&=&\Z.
\end{eqnarray*}
All other $\H^{i,j}(T_{p,q})$ for $i=0,\ldots,4$, are trivial. 
\end{theorem}
\textbf{Proof:}\\
\indent  Suppose that $p=3$. Then by applying  (\ref{f2}), we obtain that $H^{i,j}(D_{3,q})=H^{i,j}(D_{3,4})$ for $i=0,\ldots,4$.
If $p>3$, then by applying  (\ref{f1}) repeatedly, we obtain
$H^{i,j}(D_{p,q})=H^{i,j}(D_{p,p})$.
Furthermore, by applying (\ref{f3}) (and then (\ref{f1})) repeatedly we obtain $H^{i,j}(D_{p,p})=H^{i,j+p-3}(D_{3,4})$ for $i=0,1,2,3,4$.
Finally, the homology of the last torus knot is well-known, see e.g. \cite{shtor}-knot $8_{19}$, and thus we obtain the required result. \kraj\\

\section{Further thickness results}\label{zadnji}

Even though we have shown that the torus knots $T_{p,q}$, $p\ge 3$ are H-thick, from the existing experimental results one can see that the homology of torus $T_{p,q}$ knots occupies at least $p$ diagonals (i.e. that its homological width is at least $p$). In fact, one can see that in all  examples we have that $H^{2p-2,p}(D_{p,q})$ is of nonzero rank. 
\begin{proposition}
If $\rank{H^{2p-2,p}(D_{p,q})}>0$ then the homological width of the torus knot $T_{p,q}$ is at least $p$.
\end{proposition}
\textbf{Proof:}\\
\indent As we know, see e.g. the proof of Corollary \ref{C1}, there exists a generator of the homology group $\H^{0,(p-1)(q-1)+1}$ and its $\delta$-grading is equal to $(p-1)(q-1)+1$. Since we assumed that   
$$\rank{H^{2p-2,p}(D_{p,q})}=\rank{\H^{2p-2,p+(p-1)q}(T_{p,q})}>0$$
we have that there exists a generator of this homology group whose $\delta$-grading is equal to $p+(p-1)q-2(2p-2)=(p-1)(q-1)+3-2p$. So, we have two generators whose $\delta$-gradings dif\mbox{}fer by $2p-2$, and hence they lie on two dif\mbox{}ferent diagonals between which there are $p-2$ diagonals. Hence the homological width of the torus knot $T_{p,q}$ is at least $p$. \kraj\\

Thus we are left with proving that $\rank{H^{2p-2,p}(D_{p,q})}>0$ . From (\ref{f2}) we have that ${H^{2p-2,p}(D_{p,q})}={H^{2p-2,p}(D_{p,p+1})}$. Furthermore we have
\begin{lemma} 
${H^{2p-2,p}(D_{p,p})}={H^{2p-2,p}(D_{p,p+1})}.$
\end{lemma} 
\textbf{Proof:}\\

\indent In order to prove this, we will start from the diagram $D_{p,p+1}$ and we will use the same process as in the proof of Theorem \ref{T1}. Namely, we obtain the long exact sequences, see (\ref{eq3}), for every
$i=1,\ldots,p-1$:
\arraycolsep 0pt
{{
\begin{eqnarray}\cdots\rightarrow H^{2p-3,p-1}(E^i_{p,p+1})\to  H^{2p-2,p}(D^{i-1}_{p,p+1})\to \nonumber \\ \to  H^{2p-2,p}(D_{p,p+1}^i) \to H^{2p-2,p-1}(E^i_{p,p+1})\to \cdots\label{eq4}\end{eqnarray}
}}
\arraycolsep 6pt
For every $i=1,\ldots,p-1$ we can calculate explicitly the number of positive and negative crossings of $E^i_{p,p+1}$, and we can f\mbox{}ind explicitly the positive diagram to which $E^i_{p,p+1}$
is isotopic.\\

\indent One can easily see that the number of negative crossings of each $E^i_{p,p+1}$ is equal to $2p-2$ and hence the number of positive crossings is equal to $(p-1)(p+1)-i-(2p-2)=p^2-2p+1-i$. On the other hand, every $E^i_{p,p+1}$ for $i=1,\ldots,p-2$ is isotopic (by a sequence of R2 and R1 moves as explained in the proof of Lemma \ref{lem1}) to  the diagram  $D_{p-2,p-1}^{i-1}$, while $E^{p-1}_{p,p+1}$ is isotopic to the diagram $D_{p-2,p-1}^{p-3}\coprod U=D_{p-2,p-2}\coprod U$, where by $U$ we denote the unknot.  Hence we have that $\H^{l}(E_{p,p+1}^i)=0$ for $l<0$. Also, since $D_{p,q}^i$ is a positive braid knot with $p$ strands and $(p-1)q-i$ crossings, we have that  $\H^{0,(p-3)(p-2)-(i-1)\pm 1}(D_{p-2,p-1}^{i-1})=\Z$ and all other $\H^{0,j}(D_{p-2,p-1}^{i-1})$ are trivial. Hence, we have that the only nontrivial part of the zeroth homology group  of $E_{p,p+1}^i$ is given by 
 $\H^{0,(p-3)(p-2)-(i-1)\pm 1}(E_{p,p+1}^{i})=\Z$ for $i=1,\ldots,p-2$, and $\H^{0,(p-3)(p-2)-(p-3)\pm 1\pm 1}(E_{p,p+1}^{p-1})=\Z$. Thus,  we have that for every $i=1,\ldots,p-1$,  
$\H^{0,p^2-5p+4-i}(E_{p,p+1}^i)$ is trivial.\\

\indent Finally, since 
the number of negative crossings of $E_{p,p+1}^i$ is equal to $2p-2$, we have that $H^{2p-3}(E_{p,p+1}^i)$ is trivial. Furthermore, since the number of positive crossings  of $E_{p,p+1}^i$ is equal to $p^2-2p+1-i$ we have that 
$$H^{2p-2,p-1}(E_{p,p+1}^i)=\H^{0,p^2-5p+4-i}(E_{p,p+1}^i)$$
which is trivial for every $i=1,\ldots,p-1$. \\
\indent 
Hence from the long exact sequences (\ref{eq4}) we obtain 
$${H^{2p-2,p}(D_{p,p})}={H^{2p-2,p}(D_{p,p+1})},$$ as required. \kraj\\

\begin{conjecture} \label{con1}
The rank of the homology group $H^{2p-2,p}(D_{p,p})$ (and equivalently of $H^{2p-2,p}(D_{p,p+1})$) is nonzero.
\end{conjecture}

\indent As we saw, the validity of  Conjecture \ref{con1} implies that the homological width of the torus knot $T_{p,q}$ is at least $p$.\\
\indent Even though we don't (yet) have the proof of  Conjecture \ref{con1}, there is evidence that it is true. First of all, the computer program calculations show that the conjecture is true at least for $p \le 7$ (the calculations are mainly for knots, i.e. for $D_{p,p+1}$). Furthermore, Lee's variant $H_L^{i,j}$ of Khovanov homology (\cite{lee2}) for the $p$-component link $D_{p,p}$ has $2p$ generators in the homological degree $2p-2$. Also, as is well-known, there exist spectral sequences whose $E_{\infty}$-page is Lee's homology  and whose $E_2$-page is Khovanov homology (see \cite{ras}, \cite{turner}). So $H_L^{i,j}\subset H^{i,j}$ and hence $H^{2p-2}(D_{p,p})$ has at least $2p$ generators. So, we are left with proving that at least one of them has the $q$-grading equal to $p$.\\

\section{Stable $sl(n)$ homology of torus knots} \label{stable}

\indent Def\mbox{}ine the following normalization of the Poincar\'e polynomial of the homology of the torus knot:
$$P_{m,n}(t,q)=q^{-(m-1)n}P(T_{m,n})(t,q).$$ 
Then from the ``descending" properties of Theorem \ref{T3} we have the following:
\begin{theorem} For every $m \in \N$ there exists a stable homology polynomial $P^S_m$ given by:
$$P^S_m(t,q)=lim_{n\to \infty}P_{m,n}(t,q).$$
\label{stablet}
\end{theorem}

\indent Furthermore, as we have shown, the (normalized) Poincar\'e polynomial $P_{m,n}(t,q)$ of the torus knot $T_{m,n}$ coincides with the stable polynomial $P^S_m$, for all powers of $t$ up to $m+n-3$.\\
\indent Similar results are obtained at the conjectural level in \cite{dgr} (with a conjectural bound on the powers of $q$ for  agreement between the stable homology and the ef\mbox{}fective homology of any  particular torus knot). In \cite{dgr},  reduced homology (see e.g. \cite{pat}) is used, but the whole method and all proofs from the previous sections work in the same way for reduced homology. \\

\indent Also, in \cite{dgr} the existence of stable $sl(n)$ homology for torus knots is conjectured. However, in the course of proving the stability property in the $sl(2)$ case (Theorem \ref{T3}, formula (\ref{f1})) we basically used the long exact sequence in homology. Since the analogous long exact sequence exists for $sl(n)$ homology (it is again the mapping cone - see Remark \ref{rem2}), we can repeat the major part of the process. The long exact sequence in the case of $sl(n)$ homology is:
{\small{
\begin{equation}\cdots\rightarrow H^{i-1,j+1}(D_1)\to H^{i,j}(D)\to H^{i,j}(D_0)\to H^{i,j+1}(D_1)\to H^{i+1,j}(D)\to \cdots\label{lessln}\end{equation}
}}
where $D_i$, $i=0,1$ is obtained from $D$ after resolving the crossing $c$ into an $i$-resolution. Note that one of the diagrams $D_i$ is not a planar projection of a knot since it contains one thick edge. In the case that we are interested in (torus knots - positive knots), the diagram $D_1$ is the one which has one thick edge (for the details and notation see Sections \ref{sln} and \ref{kovroznot}).\\
\indent Again, we start from the diagram $D_{p,q}$ of the torus knot $T_{p,q}$, and we resolve the crossing $c_{p-1}=(p-1,1)$. We denote the diagram obtained by the 0-resolution by $D^1_{p,q}$, and the diagram obtained by the $1$-resolution by $\bar{E}^1_{p,q}$. Then we have the following long exact sequence:
{\small{
\begin{equation*}\cdots\rightarrow H_n^{i-1,j+1}(\bar{E}^1_{p,q})\to H^{i,j}_n(D_{p,q})\to H_n^{i,j}(D_{p,q}^1)\to H_n^{i,j+1}(\bar{E}^1_{p,q})\to \cdots\end{equation*}
}}
We continue the process, by resolving the crossings $c_l=(l,1)$, $l=p-2,\ldots,1$ of the diagram $D^{p-1-l}_{p,q}$ and we denote the 0- and 1-resolution obtained, by $D^{p-l}_{p,q}$ and $\bar{E}_{p,q}^{p-l}$, respectively. Then we have:
%{\small{
\begin{eqnarray*}
\cdots\rightarrow H_n^{i-1,j+1}(\bar{E}^l_{p,q})\to H^{i,j}_n(D^{l-1}_{p,q})\to H_n^{i,j}(D^l_{p,q})\to H_n^{i,j+1}(\bar{E}^l_{p,q})\to \cdots,\\
l=2,\ldots,p-1.
\end{eqnarray*}

%\begin{remark}
%Although it virtually seems that we are doing the sequence of the %resolutions in a dif\mbox{}ferent way than in the $sl(2)$ case, in fact %we have the same situation. Indeed, because of the symmetries of the %standard diagram of the torus knot, the resolutions that we are %perfoming above on the f\mbox{}irst full twist, can be also seen as %that
%they are done on the $q$-th (last) full twist, which together with one %general rotation gives the situation from the $sl(2)$ case. So, the only %dif\mbox{}ference from the $sl(2)$ case is the fact that now we are %having that two lower strands of the thick edge (i.e. of the diagram %$\bar{E}_i$), are the ones that are the ``neighbours" to each other %through the diagram (which was the case with the upper two strands in %$sl(2)$ case, see Lemma \ref{lem1}).  
%\end{remark}

Like in the $sl(2)$ case, we shall prove the following
\begin{lemma}\label{lem2} The homology group $H_n^i(\bar{E}^l_{p,q})$ is trivial for every $l<p$ and $i<p+q-3$.
\end{lemma}
This lemma, together with the above long exact sequences and the fact that $D_{p,q}^{p-1}=D_{p,q-1}$ gives
\begin{equation}
\label{spustsln}
H_n^{i,j}(D_{p,q-1})=H_n^{i,j}(D_{p,q}),\quad i<p+q-3.
\end{equation}
From this formula, we conclude 
the existence of the limit:
\begin{eqnarray*}
P^n_k(t,q)&=&\lim_{l\to\infty} \sum_{i,j\in \Z}t^iq^j \dim H_n^{i,j}(D_{k,l})=\\
&=&\lim_{l\to\infty} \sum_{i,j\in \Z}t^iq^j q^{(n-1)(k-1)l}\dim \H_n^{i,j}(T_{k,l})=\\
&=&\lim_{l\to\infty}q^{(n-1)(k-1)l}P^n(T_{k,l})(t,q),
\end{eqnarray*}
for every $k$, where $P^n(T_{k,l})(t,q)$ is the Poincar\'e polynomial of the chain complex assigned to $T_{k,l}$ by $sl(n)$-homology. In other words, we obtain
\begin{theorem}
There exists stable $sl(n)$ homology for torus knots.
\label{torussln}
\end{theorem}

\indent Thus, we are left with proving Lemma \ref{lem2}. 
%that $H_n^i(\bar{E}^l_{p,q})$ is %trivial for every $l<p$ and $i<p+q-3$. 
We will use more or less the same approach as in Lemma \ref{lem1}. Let $C_n(\bar{E}^l_{p,q})$ be the chain complex assigned by $sl(n)$-link homology (\cite{kovroz}) to  $\bar{E}^l_{p,q}$. Then $H_n^i(\bar{E}^l_{p,q})=H^i(C_n(\bar{E}^l_{p,q}))$. Note that in the $sl(n)$ case, a complex $C_n$ is assigned to  generalized regular diagrams, i.e. to  regular diagrams where we also allow trivalent vertices and thick edges.\\
\indent Since $\bar{E}^l_{p,q}$ has only positive crossings, its homology groups are trivial for negative homological degrees. However, we will show that 
\begin{equation}
C_n(\bar{E}^l_{p,q})\sim D_n[p+q-3],
\label{slnkompl}
\end{equation}
where $D_n$ is a complex such that all its chain groups $D_n^i$ are trivial for $i<0$ (in fact, we will def\mbox{}ine $D_n$ as the direct sum of the complexes of the form $C_n(D_{\Gamma}^i)$, where $D_{\Gamma}^i$'s are the generalized regular diagrams whose all crossings are positive). Here, by $\sim$ we denote a quasi-isomorphism, which implies that the two complexes have isomorphic homology groups. Then (\ref{slnkompl}) implies that $H_n^i(\bar{E}^l_{p,q})=H^i(C_n(\bar{E}^l_{p,q}))$ is trivial for $i<p+q-3$.\\
\indent Like in Lemma \ref{lem1}, we have that the two lower strands (thin edges) of the $\bar{E}_i$ part, 
%will have at least two intersections with each other 
will form at least two crossings with each other
(corresponding to an R1 move in the proof of Lemma \ref{lem1}) and both of them will have at least $p+q-5$ over- or undercrossings with the same strand (corresponding to an R2 move in the proof of Lemma \ref{lem1}). We will show that in each of these cases we can ``shift" up our complex by one homological degree, and thus obtain (\ref{slnkompl}). \\
\indent In order to prove this we will use the following fact (``cancellation principle" for chain complexes): if we quotient the chain complex $\C$ by an (arbitrary) acyclic subcomplex $\C'$ (i.e. a subcomplex with trivial homology), then the  quotient complex $\C/\C'$ is quasi-isomorphic to the complex $\C$, and so they have isomorphic homology groups (see e.g. Lemma 3.7 of \cite{bn}).\\

\paragraph{Untwisting an R1 move}

\indent First, let us work with the analog of the R1 move. Let $\bar{D}$ be the diagram that contains the following diagram as a subdiagram:

\begin{center}
\setlength{\unitlength}{4mm}
\begin{picture}(20,3.4) 
\linethickness{0.6pt}
\put(10,-0.2){
\put(0,0){\X}
\put(0,2){\GR}
}
\put(6.5,1.5){$D:$}
\end{picture}
\end{center}
Then (see Section \ref{kovroznot} and \cite{kovroz}) the chain complex $C_n(\bar{D})$ associated to $\bar{D}$ is the mapping cone of a certain homomorphism $f:C_n(\bar{D}_0)\to C_n(\bar{D}_1)\{-1\}$, where $\bar{D}_0$ and $\bar{D}_1$ are the 0- and 1-resolutions, respectively, of the crossing of $D$, i.e. they look the same as the diagram $\bar{D}$ except that its subdiagram $D$ is replaced by $D_0$ and $D_1$, respectively:

\begin{center}
\setlength{\unitlength}{4mm}
\begin{picture}(20,3.8) 
\linethickness{0.6pt}
\put(5,-0.2){
\put(0,0){\OF}
\put(0,2){\GR}
}
\put(14,-0.2){
\put(0,0){\GRF}
\put(0,2){\GR}
}
\put(2.5,1.5){$D_0:$}
\put(11.5,1.5){$D_1:$}
\end{picture}
\end{center}
In other words, the complex associated to $\bar{D}$ is the total complex of the complex given by $C_n(\bar{D}_0)\to C_n(\bar{D}_1)[1]\{-1\}$. Furthermore, we have that (see Section \ref{kovroznot}):
$$C_n(\bar{D}_1)\cong C_n(\bar{D}_0)\{1\} \oplus C_n(\bar{D}_0)\{-1\}.$$
We also have that the projection of $f$ to the f\mbox{}irst summand is an isomorphism (see \cite{kovroz}), and hence the complex $C_n(\bar{D})$ is quasi-isomorphic to $C_n(\bar{D}_0)[1]\{-2\}$ (by the cancellation principle). So, the last two complexes have isomorphic homology groups. Thus, we can ``untwist" the crossing involving two strands that are connected to the same thick edge (the analog of an R1 move in the $sl(2)$ case) by shifting the complex of the diagram obtained up by one in homological degree, as required. \\

\paragraph{Untwisting an R2 move}

\indent Hence, after untwisting  the two crossings of $\bar{E}^l_{p,q}$ that were resolved in the $sl(2)$ case by the R1 move, we are left with a diagram of the form:

\begin{equation*}
\setlength{\unitlength}{4mm}
\left.
{
\begin{picture}(5,6) 
\linethickness{0.6pt}
\put(0,0){
\put(0,4){\GR}
\qbezier(0,2)(2,4)(4,6)
\qbezier(0.5,1)(2.5,3)(4.5,5)
\qbezier(0,4)(0.4,3.6)(0.8,3.2)
\qbezier(1,4)(1.15,3.85)(1.3,3.7)
\qbezier(1.2,2.8)(1.38,2.62)(1.56,2.44)
\qbezier(1.7,3.3)(1.88,3.12)(2.06,2.94)
\qbezier(1.94,2.06)(2.47,1.53)(3,1)
\qbezier(2.44,2.56)(3.22,1.78)(4,1)
\qbezier(3.44,-3.94)(3.72,-4.22)(4,-4.5)
\qbezier(3.06,-3.56)(3,-3.5)(2.94,-3.44)
\qbezier(2,-2.5)(2.28,-2.78)(2.56,-3.06)
\qbezier(3,-3)(2.25,-3.75)(1.5,-4.5)
\qbezier(4,-3)(3.25,-3.75)(2.5,-4.5)
\qbezier(3,-3)(3,-2.8)(3,-2.6)
\qbezier(4,-3)(4,-2.8)(4,-2.6)
\put(0,1){
\qbezier(3,-0.4)(3,-0.2)(3,0)
\qbezier(4,-0.4)(4,-0.2)(4,0)
}
}
\put(3.45,-0.5){$\cdot$}
\put(3.45,-1.5){$\cdot$}
\put(3.45,-2.5){$\cdot$}
\end{picture}}
\right\}
p+q-5
\end{equation*}

%$$p+q-5 \quad R2 \quad moves$$

In other words, we have two neighbouring strands that are both connected to the same thick edge, and both go over or both go under $p+q-5$ strands. We will show that the complex corresponding to this diagram is quasi-isomorphic to the complex of diagrams whose crossings are all positive, shifted up in homological degree by $p+q-5$.\\
\indent Let $\bar{D}$ be a positive diagram that contains the following diagram as a subdiagram:

\begin{center}
\setlength{\unitlength}{4mm}
\begin{picture}(20,5.4) 
\linethickness{0.6pt}
\put(9,-0.2){
\put(0,0){\X}
\put(0,4){\GR}
\put(1,2){\X}
\qbezier(0,2)(0,3)(0,4)
\qbezier(2,6)(2,5)(2,4)
\qbezier(2,2)(2,1)(2,0)
\put(2,6){
\qbezier(0,0)(0.08,-0.08)(0.16,-0.16)
\qbezier(0,0)(-0.08,-0.08)(-0.16,-0.16)
}
}
\put(6.5,2.5){$D:$}
\end{picture}
\end{center}

\indent Denote by $D_0$, $D_1^1$, $D_1^2$ and $D_2$ the resolutions obtained from $D$ by resolving its two crossings, according to the following pictures:

\begin{center}
\setlength{\unitlength}{4mm}
\begin{picture}(20,5.8) 
\linethickness{0.6pt}
\put(5,-0.2){
\put(0,0){\OF}
\put(0,2){\OF}
\put(0,4){\GR}
\qbezier(2,6)(2,3)(2,0)
\put(2,6){
\qbezier(0,0)(0.08,-0.08)(0.16,-0.16)
\qbezier(0,0)(-0.08,-0.08)(-0.16,-0.16)
}
}
\put(14,-0.2){
\put(0,0){\GRF}
\put(0,4){\GR}
\put(0,2){\OF}
\qbezier(2,6)(2,3)(2,0)
\put(2,6){
\qbezier(0,0)(0.08,-0.08)(0.16,-0.16)
\qbezier(0,0)(-0.08,-0.08)(-0.16,-0.16)
}
}
\put(2.5,2.5){$D_0:$}
\put(11.5,2.5){$D_1^2:$}
\end{picture}
\end{center}

\begin{center}
\setlength{\unitlength}{4mm}
\begin{picture}(20,5.8) 
\linethickness{0.6pt}
\put(5,-0.2){
\put(0,0){\OF}
\put(1,2){\GRF}
\put(0,4){\GR}
\qbezier(0,2)(0,3)(0,4)
\qbezier(2,6)(2,5)(2,4)
\qbezier(2,2)(2,1)(2,0)
\put(2,6){
\qbezier(0,0)(0.08,-0.08)(0.16,-0.16)
\qbezier(0,0)(-0.08,-0.08)(-0.16,-0.16)
}

}
\put(14,-0.2){
\put(0,0){\GRF}
\put(0,4){\GR}
\put(1,2){\GRF}
\qbezier(0,2)(0,3)(0,4)
\qbezier(2,6)(2,5)(2,4)
\qbezier(2,2)(2,1)(2,0)
\put(2,6){
\qbezier(0,0)(0.08,-0.08)(0.16,-0.16)
\qbezier(0,0)(-0.08,-0.08)(-0.16,-0.16)
}
}
\put(2.5,2.5){$D_1^1:$}
\put(11.5,2.5){$D_2:$}
\end{picture}
\end{center}
Denote by $\bar{D}_0$, $\bar{D}_1^1$, $\bar{D}_1^2$ and $\bar{D}_2$ the diagrams obtained from $\bar{D}$, after replacing the subdiagram $D$  by $D_0$, $D_1^1$, $D_1^2$ and $D_2$, respectively.
Then we have that the complex $C_n(\bar{D})$ associated to the diagram $\bar{D}$ is the total complex of the following complex of complexes:

\begin{equation*}
\begin{array}{ccccc}
&&C_n(\bar{D}_1^1)[1]\{-1\}&&\\
&\nearrow&&\searrow &\\
C_n(\bar{D}_0)&&\bigoplus&&C_n(\bar{D}_2)[2]\{-2\}\\
&\searrow&&\nearrow &\\
&&C_n(\bar{D}_1^2)[1]\{-1\}&&
\end{array}
\end{equation*}
Like previously, we have that 
\begin{equation}
C_n(\bar{D}_1^2)\cong C_n(\bar{D}_0)\{1\} \oplus C_n(\bar{D}_0)\{-1\},
\label{dirsum1}
\end{equation} 
and the projection of the map from $C_n(\bar{D}_0)$ onto the f\mbox{}irst summand of $C_n(\bar{D}_1^2)[1]\{-1\}$ is an isomorphism. Hence, again we can quotient by an acyclic complex and obtain that $C_n(\bar{D})$ is quasi-isomorphic to the total complex of:
\begin{equation}
\begin{array}{ccc}
C_n(\bar{D}_1^1)[1]\{-1\}&&\\
&\searrow&\\
\bigoplus&&C_n(\bar{D}_2)[2]\{-2\} \\
&\nearrow&\\
C_n(\bar{D}_0)[1]\{-2\}&&
\end{array}
\label{kompl1}
\end{equation}
Also, we have that 
\begin{equation}
C_n(\bar{D}_2)\cong C_n(\bar{D}_0) \oplus C_n(\bar{D}_3),
\label{dirsum2}
\end{equation}
(``categorif\mbox{}ication" of the analogous version (\ref{axioma51}) of the f\mbox{}ifth axiom from Section \ref{homflynot}) and that the map from the second summand of (\ref{dirsum1}) to the f\mbox{}irst summand of (\ref{dirsum2}) is an isomorphism (see \cite{kovroz}). Here by $\bar{D}_3$ we denoted the diagram that is the same as $\bar{D}$ with the subdiagram $D$ replaced by the following diagram (see Section \ref{homflynot}):

\begin{center}
\setlength{\unitlength}{3.5mm}
\begin{picture}(22,5) 
\linethickness{0.6pt}
\put(9,0)
{\qbezier(0,0)(0.75,0.75)(1.5,1.5)
\qbezier(1.5,0)(1.5,0.75)(1.5,1.5)
\qbezier(3,0)(2.25,0.75)(1.5,1.5)
\qbezier(0,5)(0.75,4.25)(1.5,3.5)
\qbezier(1.5,5)(1.5,4.25)(1.5,3.5)
\qbezier(3,5)(2.25,4.25)(1.5,3.5)
\put(0.7,2.2){\small{3}}
\linethickness{3pt}
\qbezier(1.5,1.5)(1.5,2.5)(1.5,3.5)
\linethickness{0.6pt}
\put(1.5,0.7){
\qbezier(0,0)(-0.08,-0.08)(-0.16,-0.16)
\qbezier(0,0)(0.08,-0.08)(0.16,-0.16)
}
\put(1.5,4.5){
\qbezier(0,0)(-0.08,-0.08)(-0.16,-0.16)
\qbezier(0,0)(0.08,-0.08)(0.16,-0.16)
}
\put(0.7,0.7)
{\qbezier(0,0)(0,-0.1)(0,-0.2)
\qbezier(0,0)(-0.1,0)(-0.2,0)
}
\put(2.5,4.5)
{\qbezier(0,0)(0,-0.1)(0,-0.2)
\qbezier(0,0)(-0.1,0)(-0.2,0)
}
\put(2.3,0.7)
{\qbezier(0,0)(0,-0.1)(0,-0.2)
\qbezier(0,0)(0.1,0)(0.2,0)
}
\put(0.5,4.5)
{\qbezier(0,0)(0,-0.1)(0,-0.2)
\qbezier(0,0)(0.1,0)(0.2,0)
}
\put(-3,2.2){$D_3:$}
}

\end{picture}
\end{center}

Hence $C_n(\bar{D})$ is quasi-isomorphic to the total complex of 
\begin{equation}
C_n(\bar{D}_1^1)[1]\{-1\} \to C_n(\bar{D}_3)[2]\{-2\}
\label{zezanje}
\end{equation}
On the other hand the complex $C_n(\bar{D}_3)$ is quasi-isomorphic to the total complex of both of the following two complexes:
$$C_n(\bar{D}_0)[-1]\to C_n(\bar{D}_2),$$
and
$$C_n(\bar{D}'_0)[-1]\to C_n(\bar{D}'_2),$$
where $\bar{D}'_0$ and $\bar{D}'_2$ are the diagrams obtained from $\bar{D}$ after replacing $D$ by following two diagrams, respectively:

\begin{center}
\setlength{\unitlength}{4mm}
\begin{picture}(20,5.8) 
\linethickness{0.6pt}
\put(5,-0.2){
\put(1,0){\OF}
\put(1,4){\OF}
\put(1,2){\GRF}
\qbezier(0,6)(0,3)(0,0)
\put(0,6){
\qbezier(0,0)(0.08,-0.08)(0.16,-0.16)
\qbezier(0,0)(-0.08,-0.08)(-0.16,-0.16)
}
\put(1,6){
\qbezier(0,0)(0.08,-0.08)(0.16,-0.16)
\qbezier(0,0)(-0.08,-0.08)(-0.16,-0.16)
}
\put(2,6){
\qbezier(0,0)(0.08,-0.08)(0.16,-0.16)
\qbezier(0,0)(-0.08,-0.08)(-0.16,-0.16)
}
}
\put(2.5,2.5){$D'_0:$}
\put(14,-0.2){
\put(1,0){\GRF}
\put(1,4){\GR}
\put(0,2){\GRF}
\qbezier(2,2)(2,3)(2,4)
\qbezier(0,6)(0,5)(0,4)
\qbezier(0,2)(0,1)(0,0)
\put(0,6){
\qbezier(0,0)(0.08,-0.08)(0.16,-0.16)
\qbezier(0,0)(-0.08,-0.08)(-0.16,-0.16)
}
}
\put(11.5,2.5){$D'_2:$}
\end{picture}
\end{center}

Thus, the total complex of (\ref{kompl1}) (and hence $C_n(\bar{D})$) is quasi-isomorphic to the total complex of the following complex:

\begin{equation}
\begin{array}{ccc}
C_n(\bar{D}_1^1)[1]\{-1\}&&\\
&\searrow&\\
\bigoplus&&C_n(\bar{D}'_2)[2]\{-2\} \\
&\nearrow&\\
C_n(\bar{D}'_0)[1]\{-2\}&&
\end{array}
\label{kompl2}
\end{equation}

Thus, if we denote by $D_n$ the above complex shifted down in homological degree by 1, then we have that $C_n(\bar{D})\sim D_n[1]$, and all homology groups of $D_n$ are in nonnegative homological degrees, as required (the crossings in three diagrams appearing in (\ref{kompl2}) are all positive). We can now iterate the argument for each instance when an R2 move would occur in the $sl(2)$ case, since the two lower rightmost strands are both connected to the same thick edge in all three diagrams $D_1^1$, $D'_0$ and $D'_2$, and hence we can continue the process like with the initial diagram $D$.\\
\indent Completely analogously, we obtain the same result for the diagram $D$ of the form:

\begin{center}
\setlength{\unitlength}{4mm}
\begin{picture}(20,5.4) 
\linethickness{0.6pt}
\put(9,-0.2){
\put(1,0){\X}
\put(1,4){\GR}
\put(0,2){\X}
\qbezier(2,2)(2,3)(2,4)
\qbezier(0,6)(0,5)(0,4)
\qbezier(0,2)(0,1)(0,0)
\put(0,6){
\qbezier(0,0)(0.08,-0.08)(0.16,-0.16)
\qbezier(0,0)(-0.08,-0.08)(-0.16,-0.16)
}
}
\put(6.5,2.5){$D:$}
\end{picture}
\end{center}

Thus, we obtain the required shift in homological degree.

\appendix

\chapter{Planar graphs and alternating links}\label{a1}

In this appendix we will show that there is a bijection between planar graphs and alternating link diagrams (up to mirror image), and that the Jones polynomial of an alternating link is equal, up to a multiple, to the specialization of the dichromatic polynomial of the corresponding graph. Mainly, we follow \cite{kauf}, section II.8.\\
\indent In fact, there is a bijective correspondence between planar graphs and link shadows, i.e. the planar projection obtained by replacing each over- and undercrossing with an intersection. Also, every connected link shadow corresponds to two alternating links, one the mirror image of the another, and hence, up to mirror image, every connected planar diagram corresponds to an alternating link.\\
\indent If we have a planar graph $G$, then we assign to it the alternating link $L^G$, in the following way: we replace every edge of the graph $G$ by a crossing according to the following picture

\begin{center}
\setlength{\unitlength}{4mm}
\begin{picture}(20,4) 
\linethickness{0.6pt}
\put(5,0){
\put(0,2){\line(1,0){10}}
%\qbezier(0,2)(5,2)(10,2)
\linethickness{0.8pt}
\put(3.5,0.5){\qbezier(0,0)(1.5,1.5)(3,3)}
%{\line(1,1){3}}
%\qbezier(3.5,0.5)(5,2)(6.5,3.5)
\put(6.5,0.5){\qbezier(0,0)(-0.5,0.5)(-1,1)}
%{\line(-1,1){1}}
%(6,1)(5.5,1.5)
\put(3.5,3.5){\qbezier(0,0)(0.5,-0.5)(1,-1)}
%{\line(1,-1){1}}
\linethickness{2.5pt}
\put(0,2){\circle{0.2}}
\put(10,2){\circle{0.2}}
}
%\put(8,2.5){\line(-1,1){1}}
\end{picture}
\end{center}
and then we connect the strands near each vertex of the graph, in accordance with the following pattern:
\begin{center}
\setlength{\unitlength}{4mm}
\begin{picture}(20,6) 
\linethickness{0.5pt}

\qbezier(5,3)(10,3)(15,3)
\qbezier(10,0)(10,3)(10,6)
\linethickness{3pt}
\put(10,3){\circle{0.1}}
\put(10,3){\circle{0.15}}
\put(10,3){\circle{0.2}}
\linethickness{0.8pt}
\qbezier(9.7,4.7)(8.4,3.4)(7.1,2.1)
\qbezier(9.1,5.9)(10.4,4.6)(11.7,3.3)
\qbezier(12.9,3.9)(11.6,2.6)(10.3,1.3)
\qbezier(8.3,2.7)(9.6,1.4)(10.9,0.1)
\qbezier(7.1,3.9)(7.4,3.6)(7.7,3.3)
\qbezier(10.3,5.3)(10.6,5.6)(10.9,5.9)
\qbezier(12.3,2.7)(12.6,2.4)(12.9,2.1)
\qbezier(9.1,0.1)(9.4,0.4)(9.7,0.7)
\end{picture}
\end{center}
Obviously, in this way we obtain an alternating link.\\
\indent Conversely, from every alternating link diagram $L$, we can obtain its shadow by replacing every over- and undercrossing by an intersection (i.e. by a 4-valent vertex). Then we shade the 4-valent graph obtained in this way, in checkerboard fashion (we colour the regions in black and white such that neighbouring regions receive a dif\mbox{}ferent colour). Finally, we obtain the corresponding planar graph $G^L$ by assigning a vertex to every black region and every time two black regions touch, or one black region touches itself, at a crossing point of $L$, we connect the two corresponding vertices, or connect the vertex to itself, with an edge passing through that crossing point of $L$.\\

\indent Now, let us see the connection between the Jones polynomial for links and the dichromatic polynomial of graphs. First of all, denote the number of vertices of the graph $G$ by $N$, and the number of edges of $G$ by $m$. Fix an ordering of the  set of edges, $E(G)$, of the graph $G$. Then since the crossings of $L^G$ correspond to the edges of $G$, we obtain an ordering of the crossings of $L^G$, as well. As we saw in Chapter \ref{graf}, the dichromatic polynomial, $P_G(q,v)$ of a graph $G$ is given by:
\begin{equation}
P_G(q,v)=\sum_{\epsilon \in  \{0,1\}^m}{(-1)^{|\epsilon|}q^{|\epsilon|}v^{k(\epsilon)}},
\label{dpss}
\end{equation}
Here we denote the sum of the entries of $\epsilon$ by $|\epsilon|$, and  the number of connected components of $[G:s_{\epsilon}]$ (the spanning subgraph of $G$ whose set of edges is $s_{\epsilon}$) by 
$k(\epsilon)$. Recall that we denote by $s_{\epsilon}$ the subset of the set of edges of $G$ such that the $i$-th edge of $G$ belongs to $s_{\epsilon}$ if and only if $\epsilon_i=1$.
As we saw previously, to every edge $e$ of the graph $G$ we assign a crossing $c_e$ of $L^G$. In addition, the graph $G-e$ is naturally associated to a 0-resolution of $L^G$, i.e. the diagram that is obtained by resolving the crossing $c_e$ into a 0-resolution, and the graph $G/e$ is naturally associated to a 1-resolution of $L^G$, i.e. the diagram that is obtained by resolving the crossing $c_e$ into a 1-resolution. Hence,  every $\epsilon \in  \{0,1\}^m$ naturally corresponds to a total resolution of $L^G$, denoted by $L^G_{\epsilon}$, such that the $i$-th crossing is resolved into an $\epsilon_i$-resolution.\\
\indent As we saw in Section \ref{jonesnot}, the Jones polyomial of $L^G$, $J(L^G)(z)$, is, up to a multiple, equal to the Kauf\mbox{}fman bracket of $L^G$, $\langle L^G \rangle (z)$, which, according to formula \ref{kaufbr}, is given by:
\begin{equation}
J(L^G)(z) \sim \langle L^G \rangle(z)=\sum_{\epsilon\in\{0,1\}^m} {(-1)^{|\epsilon|}z^{|\epsilon|}(z+z^{-1})^{c(\epsilon)}}.
\label{jpss}
\end{equation}
Here by $\sim$ we denote equality up to a multiple, and by $c(\epsilon)$ we denote the number of circles in the total resolution $L^G_{\epsilon}$. Although $c(\epsilon)$ and $k(\epsilon)$ are not equal in general, there is a simple relation between them. Indeed from the Euler formula for planar graphs, one can easily obtain that:
\begin{equation}
k(\epsilon)=\frac{1}{2}\left( N - |\epsilon| + c(\epsilon) \right).
\label{smena}
\end{equation}
For more details see \cite{kauf}, section II.8. \\
\indent By substituting (\ref{smena}) into (\ref{dpss}) we obtain:
\begin{eqnarray*}
P_G(q,v)&=&\sum_{\epsilon \in  \{0,1\}^m}{(-1)^{|\epsilon|}q^{|\epsilon|}v^{k(\epsilon)}}=\\
&=& v^{N/2}\sum_{\epsilon \in  \{0,1\}^m}{(-1)^{|\epsilon|}(qv^{-1/2})^{|\epsilon|}(v^{1/2})^{c(\epsilon)}}.
\end{eqnarray*}
Now, if we take the specialization given by $v=q^2/(q-1)$ and if we denote $z=(q-1)^{1/2}$, then we have:
$$qv^{-1/2}=(q-1)^{1/2}=z,$$
and
$$v^{1/2}=\frac{q}{(q-1)^{1/2}}=(q-1)^{1/2}+\frac{1}{(q-1)^{1/2}}=z+z^{-1}.$$
Finally, by using (\ref{jpss}) we have:
\begin{equation*}
P_G(q,q^2/(q-1))\sim \sum_{\epsilon \in  \{0,1\}^m}{(-1)^{|\epsilon|}z^{|\epsilon|}(z+z^{-1})^{c(\epsilon)}}\sim J(L^G)(z),
\end{equation*}
which gives the required relation between this specialization of the dichromatic polynomial of a graph and the Jones polynomial of the corresponding alternating link.\\

\indent Also, since we know from Chapter \ref{graf} (formula (\ref{tutdich})), that the Tutte polynomial of a graph $G$ is, up to a multiple and a change of variables, equal to the dichromatic polynomial of a graph $G$, we have that 
$$T_G(x,y) \sim P_G(q,v),$$
where $q=1-y$ and $v=(x-1)(y-1)$, or equivalently $y=1-q$ and $x=1-v/q$. 
Hence, the specialization $v=q^2/(q-1)$ corresponds to the case $x=1-q/(q-1)=1/(1-q)$ and $y=1-q=1/x$.

\chapter{State-sum models for the HOMFLYPT polynomial}\label{ap2}

In this appendix we will start from the most general state-sum models involving thick edges (of the form (\ref{slnr1})-(\ref{slnr2})), and we will show that in this way we obtain two dif\mbox{}ferent diagrammatics for the HOMFLYPT polynomial of Sections \ref{homflynot} and \ref{coment}. Also, as we shall see, these are essentially the only two dif\mbox{}ferent possibilities.\\
\indent Let $L$ be an oriented link, and let $D$ be a braid presentation of $L$, i.e. $D$ is the closure of a braid $\sigma$, such that $D$ is a regular diagram of $L$. Denote the number of positive crossings of $D$ by $n_+$, the number of negative crossings by $n_-$ and the number of strands of the braid $\sigma$ by $s(D)$. We will def\mbox{}ine the value of an (intermediate) invariant $\langle D \rangle$ by f\mbox{}irst giving the value of the bracket on planar trivalent graphs with thick edges, and then we def\mbox{}ining its value on (generalized) diagrams via the following recursive relations:
\begin{eqnarray}
\langle D_+ \rangle&=&a\langle D_0 \rangle +b \langle D_{\Gamma} \rangle 
\label{apb1}\\
\langle D_- \rangle&=&c\langle D_0 \rangle +d \langle D_{\Gamma}\rangle. \label{apb2}
\end{eqnarray} 
Here, by $D_+$, $D_-$, $D_0$ and $D_{\Gamma}$ we denote the following crossings and resolutions:

\begin{center}
\setlength{\unitlength}{5mm}
\begin{picture}(25,3)
\linethickness{0.6pt}
\put(5,1){\XS}
\put(10,1){\YS}
\put(15,1){\OS}
\put(20,1){\GR}
\put(5.2,-0.2){$D_{+}$}
\put(10.2,-0.2){$D_{-}$}
\put(15.2,-0.2){$D_{0}$}
\put(20.2,-0.2){$D_{\Gamma}$}
\end{picture}
\end{center}
We will not require that the bracket is invariant under the (oriented) Reidemeister moves, but rather that it has the following simple behaviour:

\begin{equation}
\setlength{\unitlength}{4mm}
\begin{picture}(25,2)
\linethickness{0.7pt}
\put(-5,0){
\put(7.4,0.8){$\langle$}
\put(8,0){\LPP}
\put(9.8,0.8){$\rangle$}
\put(10.8,0.8){=}
\put(12,0.8){$x$}
\put(13,0.8){$\langle$}
\put(14,0){\qbezier(0,0)(0.5,1)(0,2)}
\put(15,0.8){$\rangle$}
\put(8,2){
\qbezier(0,0)(-0.03,-0.15)(-0.06,-0.3)
\qbezier(0,0)(0.14,-0.07)(0.28,-0.14)
}
\put(14,2){
\qbezier(0,0)(-0.03,-0.14)(-0.06,-0.28)
\qbezier(0,0)(0.14,-0.07)(0.28,-0.14)
}
}
\put(8,0){
\put(7.4,0.8){$\langle$}
\put(8,0){\LPM}
\put(9.8,0.8){$\rangle$}
\put(10.8,0.8){=}
\put(12,0.8){$y$}
\put(13,0.8){$\langle$}
\put(14,0){\qbezier(0,0)(0.5,1)(0,2)}
\put(15,0.8){$\rangle$}
\put(8,2){
\qbezier(0,0)(-0.03,-0.14)(-0.06,-0.28)
\qbezier(0,0)(0.14,-0.07)(0.28,-0.14)
}
\put(14,2){
\qbezier(0,0)(-0.03,-0.14)(-0.06,-0.28)
\qbezier(0,0)(0.14,-0.07)(0.28,-0.14)
}
}
\end{picture}
\label{r1xy}
\end{equation}

\begin{equation}
\setlength{\unitlength}{4mm}
\begin{picture}(25,4)
\linethickness{0.7pt}
\put(7.2,1.8){$\langle$}
\put(8,0){\X}
\put(8,2){\Y}
\put(9.6,1.8){$\rangle$}
\put(10.5,1.8){ {$=$}}
\put(12,1.8){$z$}
\put(13.2,1.8){$\langle$}
\put(14,0){\qbezier(0,0)(0.7,2)(0,4)\qbezier(1,0)(0.3,2)(1,4)}
\put(15.5,1.8){$\rangle$}
\put(8,4){
\qbezier(0,0)(-0.03,-0.14)(-0.06,-0.28)
\qbezier(0,0)(0.14,-0.07)(0.28,-0.14)
}
\put(9,4){
\qbezier(0,0)(0.03,-0.14)(0.06,-0.28)
\qbezier(0,0)(-0.14,-0.07)(-0.28,-0.14)
}
\put(14,4){
\qbezier(0,0)(-0.06,-0.15)(-0.12,-0.3)
\qbezier(0,0)(0.14,-0.07)(0.28,-0.14)
}
\put(15,4){
\qbezier(0,0)(0.06,-0.15)(0.12,-0.3)
\qbezier(0,0)(-0.14,-0.07)(-0.28,-0.14)
}
\end{picture}
\label{r2z}
\end{equation}
while it is invariant under the R3 move. We don't require any relations for the diagrams related by the Reidemeister 2B move, i.e. the move like in f\mbox{}igure (\ref{r2z}) but with one of the strands oriented downwards, since that situation cannot occur for a diagram which is the closure of a braid diagram. Here and throughout this appendix, by the above relations, we mean relations between two diagrams that look the same except near one crossing where they are as in the above pictures.  
We denote the value of the free circle by $U$:
\begin{equation}
\setlength{\unitlength}{3mm}
\begin{picture}(25,3.5) 
\linethickness{0.6pt}
\put(10,1.5){\oval(4,3)}
\put(10,3){\qbezier(0,0)(0.15,0.15)(0.3,0.3)\qbezier(0,0)(0.15,-0.15)(0.3,-0.3)}
\put(13,1.1){$=\quad U.$}

\end{picture}
\label{xyz}
\end{equation}

Then from (\ref{apb1})-(\ref{r1xy}) we have that the following relation holds:

\begin{equation} 
\setlength{\unitlength}{4mm}
\begin{picture}(25,3.5) 
\linethickness{0.6pt}
\put(8.5,0){
\put(1.8,2){\oval(1.6,2.5)}
\qbezier(1,1.5)(1,0.8)(0.5,0.2)
\qbezier(1,2.5)(1,3.2)(0.5,3.8)
\linethickness{2.5pt}
\qbezier(1,1.5)(1,2)(1,2.5)
\linethickness{0.6pt}
\put(2.6,2)
{\qbezier(0,0)(-0.12,0.12)(-0.24,0.24)
\qbezier(0,0)(0.12,0.12)(0.24,0.24)
}
\put(0.9,0.9)
{\qbezier(0,0)(0.03,-0.15)(0.06,-0.3)
\qbezier(0,0)(-0.13,-0.1)(-0.26,-0.2)
}
\put(0.8,3.32)
{\qbezier(0,0)(0,-0.16)(0,-0.32)
\qbezier(0,0)(0.14,-0.1)(0.28,-0.2)
}
\put(4.0,1.8){$=\quad B$}
\put(9.0,0){\qbezier(0,0.2)(0,2)(0,3.8)
\qbezier(0,2)(-0.12,1.88)(-0.24,1.76)
\qbezier(0,2)(0.12,1.88)(0.24,1.76)
}
}
\end{picture}
\end{equation}
i.e. the values of the bracket on these two diagrams are proportional, and we denote the coef\mbox{}f\mbox{}icient of proportionality by $B$. Hence, we have:
\begin{eqnarray*}
x&=&aU+bB,\\
y&=&cU+dB.
\end{eqnarray*}

Also, because of the recursive relations (\ref{apb1}) and (\ref{apb2}), the bracket will satisfy the relation (\ref{r2z}) (with $z=ac$), and it will be invariant under the R3 move,  if and only if the values of the bracket on trivalent graphs with thick edges satisfy the following axioms:

{\Large{
\begin{center}
\setlength{\unitlength}{5mm}
\begin{picture}(25,3.5) 
\linethickness{0.6pt}
\put(7,0){
\linethickness{2.5pt}\qbezier(0,0.5)(0,0.75)(0,1)\qbezier(0,3)(0,3.25)(0,3.5)
\linethickness{0.6pt}

\put(0,0.1){\qbezier(-0.5,2)(-0.6,1.85)(-0.7,1.7)
\qbezier(-0.5,2)(-0.4,1.85)(-0.3,1.7)}
\put(0,0.1){\qbezier(0.5,2)(0.6,1.85)(0.7,1.7)
\qbezier(0.5,2)(0.4,1.85)(0.3,1.7)}

\qbezier(0,1)(-0.5,1.6)(-0.5,2)
\qbezier(0,3)(-0.5,2.4)(-0.5,2)
\qbezier(0,1)(0.5,1.6)(0.5,2)
\qbezier(0,3)(0.5,2.4)(0.5,2)

%\qbezier(0,3)(0.15,2.95)(0.3,2.9)
%\qbezier(0,3)(-0.15,2.95)(-0.3,2.9)

%\qbezier(0,3)(0.05,2.85)(0.1,2.7)
%\qbezier(0,3)(-0.05,2.85)(-0.1,2.7)

\put(2.2,1.6){$=-\left(\frac{a}{b}+\frac{c}{d}\right)$}
\put(8.2,0){\linethickness{2.5pt}\qbezier(0,0.5)(0,2)(0,3.5)}
}
\linethickness{0.6pt}

\end{picture}
\end{center}
}}

{\Large{
\begin{center}
\setlength{\unitlength}{5mm}
\begin{picture}(22,6) 
\linethickness{0.6pt}
\put(-1,0){
\put(3,0){
\put(0,0){\GRF}
\put(0,4){\GR}
\put(1,2){\GRF}
\qbezier(2,0)(2,1)(2,2)
\qbezier(0,4)(0,3)(0,2)
\qbezier(2,6)(2,5)(2,4)
\put(2,6)
{\qbezier(0,0)(-0.08,-0.08)(-0.16,-0.16)
\qbezier(0,0)(0.08,-0.08)(0.16,-0.16)
}
\put(0,3)
{\qbezier(0,0)(-0.08,-0.08)(-0.16,-0.16)
\qbezier(0,0)(0.08,-0.08)(0.16,-0.16)
}
\put(2,1)
{\qbezier(0,0)(-0.08,-0.08)(-0.16,-0.16)
\qbezier(0,0)(0.08,-0.08)(0.16,-0.16)
}
\put(1,2)
{\qbezier(0.00,2.00)(0,1.90)(0,1.80)
\qbezier(0.00,2.00)(0.1,2)(0.20,2)
}
\put(0.7,-1.7)
{\qbezier(0.00,2.00)(0,1.90)(0,1.80)
\qbezier(0.00,2.00)(0.1,2)(0.20,2)
}
\put(-0.7,-1.7)
{\qbezier(1.00,2.00)(1,1.9)(1,1.80)
\qbezier(1.00,2.00)(0.9,2)(0.8,2)
}
\put(0,0)
{\qbezier(1.00,2.00)(1,1.9)(1,1.80)
\qbezier(1.00,2.00)(0.9,2)(0.8,2)
}
}
\put(6,2.9){$-\quad\frac{ac}{bd}$}

\put(10,0){
\put(0,2){\GRF}
\qbezier(1,0)(1,1)(1,2)
\qbezier(0,0)(0,1)(0,2)
\qbezier(0,6)(0,5)(0,4)
\qbezier(1,6)(1,5)(1,4)
\qbezier(2,0)(2,3)(2,6)
\put(0,6)
{\qbezier(0,0)(-0.08,-0.08)(-0.16,-0.16)
\qbezier(0,0)(0.08,-0.08)(0.16,-0.16)
}
\put(1,6)
{\qbezier(0,0)(-0.08,-0.08)(-0.16,-0.16)
\qbezier(0,0)(0.08,-0.08)(0.16,-0.16)
}
\put(2,6)
{\qbezier(0,0)(-0.08,-0.08)(-0.16,-0.16)
\qbezier(0,0)(0.08,-0.08)(0.16,-0.16)
}
\put(0,1)
{\qbezier(0,0)(-0.08,-0.08)(-0.16,-0.16)
\qbezier(0,0)(0.08,-0.08)(0.16,-0.16)
}
\put(1,1)
{\qbezier(0,0)(-0.08,-0.08)(-0.16,-0.16)
\qbezier(0,0)(0.08,-0.08)(0.16,-0.16)
}

}
\put(13,2.85){$=$}
\put(15,0){
\put(0,2){\GRF}
\put(1,4){\GR}
\put(1,0){\GRF}
\qbezier(0,0)(0,1)(0,2)
\qbezier(2,4)(2,3)(2,2)
\qbezier(0,6)(0,5)(0,4)
\put(0,6)
{\qbezier(0,0)(-0.08,-0.08)(-0.16,-0.16)
\qbezier(0,0)(0.08,-0.08)(0.16,-0.16)
}
\put(0,1)
{\qbezier(0,0)(-0.08,-0.08)(-0.16,-0.16)
\qbezier(0,0)(0.08,-0.08)(0.16,-0.16)
}
\put(2,3)
{\qbezier(0,0)(-0.08,-0.08)(-0.16,-0.16)
\qbezier(0,0)(0.08,-0.08)(0.16,-0.16)
}
\put(1,0)
{\qbezier(0.00,2.00)(0,1.90)(0,1.80)
\qbezier(0.00,2.00)(0.1,2)(0.20,2)
}
\put(1.7,-1.7)
{\qbezier(0.00,2.00)(0,1.90)(0,1.80)
\qbezier(0.00,2.00)(0.1,2)(0.20,2)
}
\put(0.3,-1.7)
{\qbezier(1.00,2.00)(1,1.9)(1,1.80)
\qbezier(1.00,2.00)(0.9,2)(0.8,2)
}
\put(0,2)
{\qbezier(1.00,2.00)(1,1.9)(1,1.80)
\qbezier(1.00,2.00)(0.9,2)(0.8,2)
}

}

\put(18,2.9){$-\quad\frac{ac}{bd}$}
}
\put(20,0){
\put(1,2){\GRF}
\qbezier(1,0)(1,1)(1,2)
\qbezier(2,0)(2,1)(2,2)
\qbezier(2,6)(2,5)(2,4)
\qbezier(1,6)(1,5)(1,4)
\qbezier(0,0)(0,3)(0,6)
\put(0,6)
{\qbezier(0,0)(-0.08,-0.08)(-0.16,-0.16)
\qbezier(0,0)(0.08,-0.08)(0.16,-0.16)
}
\put(1,6)
{\qbezier(0,0)(-0.08,-0.08)(-0.16,-0.16)
\qbezier(0,0)(0.08,-0.08)(0.16,-0.16)
}
\put(2,6)
{\qbezier(0,0)(-0.08,-0.08)(-0.16,-0.16)
\qbezier(0,0)(0.08,-0.08)(0.16,-0.16)
}
\put(1,1)
{\qbezier(0,0)(-0.08,-0.08)(-0.16,-0.16)
\qbezier(0,0)(0.08,-0.08)(0.16,-0.16)
}
\put(2,1)
{\qbezier(0,0)(-0.08,-0.08)(-0.16,-0.16)
\qbezier(0,0)(0.08,-0.08)(0.16,-0.16)
}

}

\end{picture}
\end{center}
}}

%Note that we have $z=ac$.\\
\indent We will obtain the invariant $I(D)$ as an overall multiple of the bracket:
\begin{equation}
I(D)={\alpha}^{n_+}{\beta}^{n_-}{\gamma}^{s(D)}\langle D \rangle.
\label{invi}
\end{equation}
We determine the coef\mbox{}f\mbox{}icients $\alpha$, $\beta$ and $\gamma$ from the requirement that they must compensate the values $x$, $y$ and $z$ from the behaviour of the bracket under the Reidemeister moves
%(formulas (\ref{r1xy}) and (\ref{r2z}))
. Thus, from (\ref{invi}), (\ref{r1xy}) and (\ref{r2z}) we obtain:

\begin{center}
\setlength{\unitlength}{4mm}
\begin{picture}(25,2)
\linethickness{0.7pt}
\put(-6,0)
{
\put(7,0.8){$I($}
\put(8,0){\LPP}
\put(9.8,0.8){$)$}
\put(10.8,0.8){=}
\put(11.8,0.8){$\alpha\gamma x$}
\put(14,0.8){$I($}
\put(15.2,0){\qbezier(0,0)(0.5,1)(0,2)}
\put(16,0.8){$)$}
\put(8,2){
\qbezier(0,0)(-0.03,-0.15)(-0.06,-0.3)
\qbezier(0,0)(0.14,-0.07)(0.28,-0.14)
}
\put(15.2,2){
\qbezier(0,0)(-0.03,-0.14)(-0.06,-0.28)
\qbezier(0,0)(0.14,-0.07)(0.28,-0.14)
}

}
\put(7,0){
\put(7,0.8){$I($}
\put(8.2,0){\LPM}
\put(9.8,0.8){$)$}
\put(10.8,0.8){=}
\put(11.8,0.8){$\beta\gamma y$}
\put(14,0.8){$I($}
\put(15.2,0){\qbezier(0,0)(0.5,1)(0,2)}
\put(16,0.8){$)$}
\put(8.2,2){
\qbezier(0,0)(-0.03,-0.15)(-0.06,-0.3)
\qbezier(0,0)(0.14,-0.07)(0.28,-0.14)
}
\put(15.2,2){
\qbezier(0,0)(-0.03,-0.14)(-0.06,-0.28)
\qbezier(0,0)(0.14,-0.07)(0.28,-0.14)
}

}
\end{picture}
\end{center}

%\glup

\begin{center}

\setlength{\unitlength}{3mm}
\begin{picture}(25,4)
\linethickness{0.7pt}
\put(6.5,1.8){$I($}
\put(8,0){\X}
\put(8,2){\Y}
\put(9.6,1.8){$)$}
\put(10.5,1.8){=}
\put(12,1.8){$\alpha\beta z$}
\put(15,1.8){$I($}
\put(16.2,0){\qbezier(0,0)(0.7,2)(0,4)\qbezier(1,0)(0.3,2)(1,4)}
\put(17.8,1.8){$)$}
\put(8,4){
\qbezier(0,0)(-0.03,-0.14)(-0.06,-0.28)
\qbezier(0,0)(0.14,-0.07)(0.28,-0.14)
}
\put(9,4){
\qbezier(0,0)(0.03,-0.14)(0.06,-0.28)
\qbezier(0,0)(-0.14,-0.07)(-0.28,-0.14)
}
\put(16.2,4){
\qbezier(0,0)(-0.06,-0.15)(-0.12,-0.3)
\qbezier(0,0)(0.14,-0.07)(0.28,-0.14)
}
\put(17.2,4){
\qbezier(0,0)(0.06,-0.15)(0.12,-0.3)
\qbezier(0,0)(-0.14,-0.07)(-0.28,-0.14)
}

\end{picture}

\end{center}

So, we def\mbox{}ine $\alpha$, $\beta$ and $\gamma$ such that the above three coef\mbox{}f\mbox{}icients are equal to 1, i.e.:
\begin{eqnarray}
\alpha&=&{\left(\frac{y}{xz}\right)}^{1/2}={\left(\frac{cU+dB}{ac(aU+bB)}\right)}^{1/2},\\
\beta&=&{\left(\frac{x}{yz}\right)}^{1/2}={\left(\frac{aU+bB}{ac(cU+dB)}\right)}^{1/2},\\
\gamma&=&{\left(\frac{z}{xy}\right)}^{1/2}={\left(\frac{ac}{(aU+bB)(cU+dB)}\right)}^{1/2}.
\end{eqnarray}

Note that:
\begin{equation}
(a\alpha) (c\beta)=ac/z=1.
\label{jednako}
\end{equation}

\indent Finally, let us see what skein relation the invariant $I$ satisf\mbox{}ies. First of all, from the recursive relations (\ref{apb1}) and (\ref{apb2}), the bracket $\langle,\rangle$ satisf\mbox{}ies the following relation:

\begin{equation*}
d\langle D_+ \rangle - b\langle D_- \rangle = (ad-bc) \langle D_0 \rangle.
\end{equation*} 

Thus, from the def\mbox{}inition of $I(D)$ (\ref{invi}), we obtain: 

\begin{equation*}
\frac{d}{\alpha} I(D_+) - \frac{b}{\beta} I(D_-) = (ad-bc) I(D_0).
\end{equation*} 

After dividing both sides of the above equality by $(abcd)^{1/2}$, denoting $q={(ad/bc)}^{1/2}$, and bearing in mind (\ref{jednako}), we get:

\begin{equation}
(q\beta c)\,\, I(D_+) - (q\beta c)^{-1}\,\, I(D_-) = (q-q^{-1})\,\, I(D_0).
\label{skeinkon}
\end{equation} 

Thus, we have shown that $I(D)$ satisf\mbox{}ies the skein relation of the HOMFLYPT polynomial. 
\begin{remark} \label{apbr1}
Note that $q$ and 
$$\beta c=\left(\frac{U+\frac{b}{a}B}{U+\frac{d}{c}B}\right)^{1/2},$$
only depend on the values $U,$ $B$ and the quotients $b/a$ and $d/c$. Also note that we would obtain the same invariant, if we def\mbox{}ined the skein relations for the bracket (\ref{apb1}) and (\ref{apb2}) such that $a$ and $c$ are replaced by 1, and $b$ and $d$ by $b/a$ and $d/c$, respectively, and then def\mbox{}ined the invariant $I(D)$ by (\ref{invi}) with $\alpha$ replaced by $a\alpha$ and $\beta$ replaced by $c\beta$. 
\end{remark}
\indent So, from now on, we will assume that $a=c=1$, and hence $d=q^2 b$. In order to obtain the $sl(n)$ specialization of the HOMFLYPT polynomial, we must have (see (\ref{skeinkon})):
\begin{equation}
q\beta c=q^{\pm n}
\label{plusminus}
\end{equation}
In the case when we have the plus sign, we obtain the following relation:
$$U+bB=q^{2n-2}(U+bq^2B),$$
or equivalently:
\begin{equation*}
-q^{-1}[n-1]U=[n]bB.
\end{equation*}
Here, by $[n]$ we have denoted the quantum integer $n$, i.e. $[n]=(q^n-q^{-n})/(q-q^{-1})$. Hence, we have a natural solution:
$$U=[n],\quad B=[n-1],\quad b=-q^{-1},\quad d=q.$$
 In this case, the parameter $\gamma$ in the def\mbox{}inition of the invariant $I$ is equal to 1, and if we add the following requirement 
\begin{center}
\setlength{\unitlength}{5mm}
\begin{picture}(35,3.5) 
\linethickness{0.6pt}
\put(3,0)
{
\put(2,2){\oval(2,3)}
\qbezier(1,1.5)(1,0.8)(0.5,0.2)
\qbezier(1,2.5)(1,3.2)(0.5,3.8)
\qbezier(3,1.5)(3,0.8)(3.5,0.2)
\qbezier(3,2.5)(3,3.2)(3.5,3.8)
\linethickness{2.5pt}
\qbezier(1,1.5)(1,2)(1,2.5)
\qbezier(3,1.5)(3,2)(3,2.5)
\linethickness{0.6pt}
%\put(2.6,2)
%{\qbezier(0,0)(-0.12,0.12)(-0.24,0.24)
%\qbezier(0,0)(0.12,0.12)(0.24,0.24)
%}
\put(0.9,0.9)
{\qbezier(0,0)(0.03,-0.15)(0.06,-0.3)
\qbezier(0,0)(-0.13,-0.1)(-0.26,-0.2)
}
\put(3.1,3.1)
{\qbezier(0,0)(-0.03,0.15)(-0.06,0.3)
\qbezier(0,0)(0.13,0.1)(0.26,0.2)
}
\put(0.8,3.32)
{\qbezier(0,0)(0,-0.16)(0,-0.32)
\qbezier(0,0)(0.14,-0.1)(0.28,-0.2)
}
\put(3.2,0.68)
{\qbezier(0,0)(0,0.16)(0,0.32)
\qbezier(0,0)(-0.14,0.1)(-0.28,0.2)
}
\put(1.9,0.5)
{\qbezier(0,0)(0.14,0.1)(0.28,0.2)
\qbezier(0,0)(0.14,-0.1)(0.28,-0.2)
}
\put(2.1,3.5){
\qbezier(0,0)(-0.14,0.1)(-0.28,0.2)
\qbezier(0,0)(-0.14,-0.1)(-0.28,-0.2)
}
}
\put(8.5,1.8){$=$}
\put(10,0)
{\qbezier(0.5,0.2)(2,2.5)(3.5,0.2)
\qbezier(0.5,3.8)(2,1.5)(3.5,3.8)
\put(2.1,1.35)
{\qbezier(0,0)(-0.14,0.1)(-0.28,0.2)
\qbezier(0,0)(-0.14,-0.1)(-0.28,-0.2)
}
\put(1.9,2.65)
{\qbezier(0,0)(0.14,0.1)(0.28,0.2)
\qbezier(0,0)(0.14,-0.1)(0.28,-0.2)
}
}
\put(14.5,1.8){$+\,\,[n-2]$}
\put(18.5,0)
{\qbezier(0.5,0.2)(2.2,2)(0.5,3.8)
\qbezier(3.5,0.2)(1.8,2)(3.5,3.8)
\put(1.35,2.1)
{\qbezier(0,0)(-0.1,-0.14)(-0.2,-0.28)
\qbezier(0,0)(0.1,-0.14)(0.2,-0.28)
}
\put(2.65,1.9)
{\qbezier(0,0)(-0.1,0.14)(-0.2,0.28)
\qbezier(0,0)(0.1,0.14)(0.2,0.28)
}
}
\end{picture}
\end{center}
This can be easily seen to imply invariance under the Reidemeister 2B move (the move with diagram as in the picture (\ref{r2z}) when the orientations of the two strands are opposite). Hence, this model works for arbitrary regular diagrams, and not only for the closures of braid diagrams, and this is the model used for the def\mbox{}inition of the HOMFLYPT polynomial in Section \ref{homflynot} of Chapter \ref{not}.\\
\indent In the case of the minus sign in (\ref{plusminus}), we have:
$$U+bB=q^{-2n-2}(U+bq^2B),$$
which gives
$$U(1+q^2+q^4+\ldots+q^{2n})=-bq^2B(1+q^2+\dots+q^{2(n-1)}).$$
The last equation has the solution
$$U=1+q^2+\ldots+q^{2(n-1)},\quad B=1+q^2+\ldots+q^{2n},\quad b=-q^{-2},\quad d=-1.$$
In this model, the value of $\gamma$ is not equal to 1, and hence the invariant $I(D)$ can be def\mbox{}ined only by starting from diagrams that are the closures of braid diagrams. Also, following Remark \ref{apbr1}, if we def\mbox{}ine $a=q^2$ and $c=q^{-2}$, and multiply $b$ and $d$ accordingly, we obtain the model used in Section \ref{coment} of Chapter \ref{Koszul} for the categorif\mbox{}ication of the $n$-specializations of the HOMFLYPT polynomial.

\end{document}